\documentclass[preprint,review,11pt,authoryear]{elsarticle}

\usepackage{amsmath,amssymb}
\usepackage{amsthm}
\usepackage{bbm, bm}
\usepackage{graphicx}
\usepackage{color}
\definecolor{cornell-red}{RGB}{179,27,27}
\usepackage{subcaption}
\usepackage{mathtools}
\usepackage{dsfont}
\usepackage{color}
\usepackage{natbib}
\usepackage{multirow}
\usepackage{pgfplots}
\usepackage{caption,subcaption}
\usepackage[left=1in,right=1in,top=1in,bottom=1in]{geometry}
\usepackage{indentfirst}
\usepackage{enumitem}
\usepackage[lined,boxed,commentsnumbered, ruled]{algorithm2e}

\usepackage{xcolor}
\usepackage[colorlinks = true,linkcolor = blue,urlcolor  = blue,citecolor = blue,anchorcolor = blue]{hyperref}

\usepackage{makecell, hhline}  

\usepackage{natbib}

\usepackage{multibib}
\newcites{ec}{References}

\graphicspath{ {Figures/} }

\newtheorem{theorem}{Theorem}
\newtheorem{lemma}{Lemma}
\newtheorem{proposition}{Proposition}

\theoremstyle{definition}

\newtheorem{example}{Example}[section]
\newtheorem{assumption}{Assumption}

\makeatletter  
\def\th@remark{%
  \thm@headfont{\bfseries}%
  \normalfont 
  \thm@preskip\topsep \divide\thm@preskip\tw@
  \thm@postskip\thm@preskip
}
\makeatother

\theoremstyle{remark}
\newtheorem{remark}{\textbf{Remark}}

\numberwithin{table}{section}

\newcommand{\ra}[1]{\renewcommand{\arraystretch}{#1}}  

\DeclareMathOperator*{\argmin}{arg\,min}

\DeclareMathOperator{\Var}{Var}

\DeclareMathOperator*{\cl}{cl}
\DeclareMathOperator*{\intr}{int}

\DeclareMathOperator{\Pro}{Prob} 

\newcommand{\B}{\mathbb{B}}
\newcommand{\E}{\mathbb{E}}
\newcommand{\Prob}{\mathbb{P}}
\newcommand{\N}{\mathbb{N}}

\newcommand{\R}{\mathbb{R}}

\newcommand{\convp}{\overset{p}{\to}}
\newcommand{\as}{\overset{a.s.}{\to}}
\newcommand{\Probh}{\widehat{\mathbb{P}}}

\newcommand{\calB}{\mathcal{B}}
\newcommand{\calC}{\mathcal{C}}
\newcommand{\calD}{\mathcal{D}}

\newcommand{\calF}{\mathcal{F}}

\newcommand{\calI}{\mathcal{I}}

\newcommand{\calN}{\mathcal{N}}

\newcommand{\calP}{\mathcal{P}}

\newcommand{\calS}{\mathcal{S}}

\newcommand{\calU}{\mathcal{U}}
\newcommand{\calV}{\mathcal{V}}

\newcommand{\calX}{\mathcal{X}}
\newcommand{\calY}{\mathcal{Y}}
\newcommand{\calZ}{\mathcal{Z}}

\newcommand{\one}{\bm{1}}
\newcommand{\zero}{\bm{0}}

\newcommand{\tp}{\top}
\newcommand{\norm}[1]{\left\lVert#1\right\rVert}
\newcommand{\norms}[1]{\lVert#1\rVert}

\newcommand{\hb}{\pmb{h}}

\newcommand{\qb}{\pmb{q}}

\renewcommand{\sb}{\pmb{s}}

\newcommand{\ub}{\pmb{u}}
\newcommand{\vb}{\pmb{v}}
\newcommand{\wb}{\pmb{w}}
\newcommand{\xb}{\pmb{x}}
\newcommand{\yb}{\pmb{y}}
\newcommand{\zb}{\pmb{z}}

\newcommand{\Cb}{\pmb{C}}

\newcommand{\Tb}{\pmb{T}}

\newcommand{\Wb}{\pmb{W}}

\newcommand{\Zb}{\pmb{Z}}

\newcommand{\xbbar}{\overline{\xb}}

\newcommand{\xbt}{\widetilde{\xb}}
\newcommand{\xbhat}{\widehat{\xb}}

\newcommand{\vhat}{\widehat{v}}
\newcommand{\Qh}{\widehat{Q}}
\newcommand{\Fh}{\widehat{F}}

\newcommand{\calXhat}{\widehat{\calX}}
\newcommand{\Ch}{\widehat{C}}
\newcommand{\mh}{\widehat{m}}
\newcommand{\calUhat}{\widehat{\calU}}
\newcommand{\rbar}{\overline{r}}
\newcommand{\ubar}{\overline{u}}

\newcommand{\pib}{\pmb{\pi}}

\newcommand{\sigmab}{\pmb{\sigmab}}

\newcommand{\xib}{\pmb{\xi}}

\newcommand{\Pit}{\widetilde{\Pi}}
\newcommand{\Pih}{\widehat{\Pi}}
\newcommand{\epst}{\widetilde{\varepsilon}}

\newcommand{\Lbar}{\overline{L}}
\newcommand{\Ubar}{\overline{U}}

\newcommand{\co}{c^{\textup{\tiny o}}}
\newcommand{\cu}{c^{\textup{\tiny u}}}
\newcommand{\cw}{c^{\textup{\tiny w}}}

\newcommand{\itt}{\widetilde{i}}

\usepackage{titlesec}

\titlespacing{\section}{0pt}{0ex}{0ex}
\titlespacing{\subsection}{0pt}{0ex}{0ex}
\titlespacing{\subsubsection}{0pt}{0ex}{0ex}

\makeatletter
\def\ps@pprintTitle{%
  \let\@oddhead\@empty
  \let\@evenhead\@empty
  \let\@oddfoot\@empty
  \let\@evenfoot\@oddfoot
}
\makeatother

\newcommand{\MYT}[1]{%
  \nonumnote{\textit{Man Yiu Tsang,\enspace}#1}}
\newcommand{\TS}[1]{%
  \nonumnote{\textit{Tony Sit (corresponding author),\enspace}#1}}
\newcommand{\HYW}[1]{%
  \nonumnote{\textit{Hoi Ying Wong,\enspace}#1}}

\begin{document}

\begin{frontmatter}

\title{Contextual Quantile Minimization for Two-Stage Stochastic Programs}

\author{Man Yiu Tsang}
\author{Tony Sit} 
\author{Hoi Ying Wong\vspace{-5mm}} 

\MYT{Department of Industrial, Manufacturing, and Systems Engineering, Texas Tech University, Lubbock, TX, USA}
\TS{Department of Statistics and Data Science, The Chinese University of Hong Kong, Shatin, N.T., Hong Kong}
\HYW{Department of Statistics and Data Science, The Chinese University of Hong Kong, Shatin, N.T., Hong Kong}





\begin{abstract}

\noindent  Contextual stochastic optimization is an advanced methodology to model uncertainty in the presence of contextual information during decision planning processes. Although classical methodologies focus on minimizing the expectation of a random loss, in many applications, risk-averse decision-makers may be interested in minimizing a specific quantile as a more prudent alternative. In this paper, we propose a new risk-averse contextual stochastic optimization problem with a quantile objective for general two-stage problems. Given historical data on the model's random parameters and contextual information, we model the conditional quantile by replacing the conditional expectation in its variational characterization with a generic estimator. Under two sets of mild regularity conditions, we derive the asymptotic almost-sure convergence and convergence in probability of the optimal solution and the optimal value of the associated optimization problem to their true counterparts. Optimization problems with a quantile objective is usually non-convex, which are generally regarded as challenging to solve. To address the computational difficulties, we propose a new stochastic inexact constraint generation method with convergence guarantee. Finally, through numerical experiments on a single-server appointment scheduling problem, we study the computational performance of our proposed solution method as well as operational performance of our proposed methodology. Our results demonstrate the importance of incorporating useful contextual information and decision-maker's risk attitude into the optimization model.


\end{abstract}
\end{frontmatter}

\section{Introduction} \label{sec:introduction}

Decision-making processes often involve uncertainty. For instance, in financial portfolio optimization problems where investors devise their investment strategies, returns from assets are changing over time; in supply chain management where suppliers determine their order quantities, customers’ demands are unknown. Due to the presence of uncertainty, classical deterministic modeling approaches cannot be directly applied. In contrast, stochastic programming has received much attention as a methodology to model uncertainty in decision-making processes over the past decades. With mature theories (see, e.g., \citealp{Shapiro_et_al:2014}), it has been successfully applied to a wide range of real-life applications, for example, in the fields of finance, health care, and supply chain management, to name but a few \citep{Grass_Fischer:2016, Grieco_et_al:2021, Li_Grossmann:2021}.

Unlike classical stochastic programming approaches that aim at minimizing the expected cost, a recent line of research focuses on incorporating contextual information (a.k.a. side information and covariates). {The emergence of this new research direction is partly due to the large amount of data that a decision-maker can access and collect nowadays.} Specifically, the contextual stochastic programming approach aims at minimizing the \textit{conditional} expected cost given the contextual information. Therefore, they are expected to produce a better decision due to improvement in the prediction of random parameters; see, e.g., \cite{Ban_Rudin:2019, Bertsimas_Kallus:2020}.

However, in existing literature, formulation of contextual stochastic programming is often based on a risk-neutral framework, which uses expected loss in the objective. {As pointed out by \cite{Li_et_al:2022}, in practice, decision-makers may be interested in minimizing quantiles of the random loss, in particular quantiles in the upper tail parts that correspond to severe losses. Thus, as argued in \cite{Wang_et_al:2018}, using the quantile objective could be more attractive than the expectation objective when the decision-maker is risk-averse and the sacrifice on expected loss is comparatively small. In such a situation, potential extreme losses could be protected with a modest increase in expected loss. The use of quantile objective has been employed in various application domains. In the financial sector, quantile, commonly known as value-at-risk (VaR), is widely adopted due to regulatory requirements (e.g., Basel III and Solvency II). In optimal treatment regime problems, the decision-maker may want to optimize the (lower) quantiles of the treatment effect with a little sacrifice on expected outcomes \citep{Beyerlein:2014, Wang_et_al:2018}. Other applications that adopt a quantile objective include facility location \citep{Chen_et_al:2006}, supply chain network design \citep{Soleimani_et_al:2014}, and the appointment scheduling problem \citep{Sang_et_al:2021}.}

The use of quantile objective usually leads to a non-convex optimization problem. Many existing works consider an alternative risk measure, conditional value-at-risk (CVaR), which is essentially the average losses beyond a given quantile. Admitting a linear reformulation \citep{Rockafellar_Uryasev:2000}, CVaR has been a popular choice of risk measure in different optimization models \citep{Filippi_et_al:2020}. However, as shown in \cite{Gneiting:2011} and \cite{Ziegel:2016}, quantile is elicitable but CVaR is not, where elicitability describes how accurately a risk measure can be forecasted. Moreover, in \cite{Jiang_Kou:2021}, through asymptotic expansions for relative errors of quantile and CVaR, they showed that simulating quantile is easier than CVaR. {As discussed earlier, the quantile objective is of interest in many applications. It is also more readily to be employed in practice since it directly measures the quantile of the given loss function.} Therefore, in this paper, we propose the risk-averse contextual stochastic optimization problem with a quantile objective, which takes into account decision-makers’ risk preferences. Specifically, we propose a generic framework for modeling conditional quantiles in the presence of contextual information, as well as the corresponding theoretical guarantees.

Besides the modeling perspective, we propose a novel stochastic inexact column generation (SiCG) method to solve linear two-stage stochastic programs with a quantile objective. As discussed in \cite{Pavlikov_et_al:2017}, optimization models involving quantile are challenging in general. Although a mixed-integer program (MIP) reformulation is possible, such an MIP typically involves a large number of variables and constraints, as well as big-$M$ parameters. Hence, directly solving the MIP reformulation could be computationally expensive. The proposed SiCG method, as a decomposition algorithm, allows solving a relaxed problem at each iteration inexactly. Thus, our SiCG method could provide a more computationally efficient solution approach when compared with the direct MIP reformulation. To examine the performance of our proposed methodology and solution method, we conduct numerical experiments on a classical appointment scheduling problem. 

\subsection{Contributions}

The contributions of this paper are summarized into the following three aspects, namely methodology, computational algorithm, and application.

\begin{itemize}[topsep=1pt,itemsep=-1pt,leftmargin=5mm]
    \item \textit{Methodology}. We propose a risk-averse contextual stochastic optimization problem with a quantile objective, which has not been studied in the exiting literature. To estimate the conditional quantile from a given set of data, we leverage its variational characterization as the minimizer of the conditional expectation of the checked loss function. Specifically, we obtain the estimated conditional quantile by replacing the conditional expectation with a generic estimator. Under different sets of regularity conditions, we derive the asymptotic almost-sure convergence and convergence in probability of the optimal value and solutions of the associated optimization model to their true counterparts.
    
    \item \textit{Computational Algorithm}. Unlike existing literature that focuses on an expectation objective, our proposed model, which involves a quantile objective, poses fundamental optimization challenges. To address this, we propose a novel stochastic inexact column generation (SiCG) method to solve general linear stochastic programs with a quantile objective. To the best of our knowledge, our proposed SiCG method is the first decomposition method with convergence guarantee for solving two-stage linear stochastic programs with a quantile objective. This new decomposition method allows inexact solves of a relaxed optimization problem in each iteration. With the introduced stochasticity, we verify that our SiCG converges in a finite number of solutions almost surely. Although motivated by the two-stage setting, where the random cost function is the optimal value of a linear program, our proposed SiCG method could be readily extended to solve problems with a general decomposable cost function.
    
    \item \textit{Application}. We conduct numerical experiments on a classical single-server appointment scheduling problem (ASP) {as an example to demonstrate the performance our proposed methodology and algorithm}. Our computational results demonstrate the computational efficiency of our proposed SiCG when compared with the standard MIP reformulation. We also derive managerial insights in ASP when contextual information is present in the planning process with the use of a quantile objective. While being demonstrated via a single-serve ASP, our proposed approach can be applied to general two-stage stochastic programs.
\end{itemize}

\subsection{Organization of the Paper}

The paper is organized as follows. In Section~\ref{sec:literature_review}, we review literature related to our work. In Section~\ref{sec:CQM}, we introduce the proposed risk-averse optimization problem with a quantile objective and the associated modeling methodology. In Section~\ref{sec:convergence_CQM}, we discuss the convergence properties of our proposed model. Next, in Section~\ref{sec:SiCG}, we propose the new SiCG algorithm to solve general linear stochastic programs with a quantile objective. Experiments results investigating the operational and computational performance of our proposed modeling and solution techniques on the appointment scheduling problem are discussed in Section~\ref{sec:numerical_expt}. Finally, we conclude in Section~\ref{sec:conclusion}. We relegate all the proofs to \ref{appdx:math_proof}.

\section{Relevant Literature} \label{sec:literature_review}

In this section, we review literature related to our work in two different perspectives: recent modeling approaches on contextual stochastic optimization and computational methods for solving stochastic programs with a quantile objective.

\subsection{Contextual Stochastic Optimization}

{Contextual stochastic optimization has received growing interests since predictors could be used to support decision-making processes. Thus, data-driven decisions with better performance could be obtained. We refer readers to \cite{Sadana_et_al:2025} for a recent survey on contextual stochastic optimization methodologies. While there are various strategies to model the conditional expectation objective, we here review three different existing methodologies.}

The first approach is to integrate prediction and optimization. Motivated by the fact that the loss function is independent of the optimization problem in classical statistical estimation, this approach formulates the loss function using the optimal solution of the optimization problem. This integration results in a challenging estimation problem, and recent research has been devoted to addressing the associated computational issues. \cite{Donti_et_al:2017} assumed that the optimization problem is strongly convex with stochastic inequality and deterministic equality constraints. For this setting, a gradient-based numerical method was proposed. In a simplified setting, \cite{Elmachtoub_Grigas:2022} considered a linear optimization problem and proposed a convex loss function which results in either a tractable linear reformulation or a convex optimization problem, depending on the polyhedrality of the constraint set. \cite{Qi_et_al:2025} generalized their framework to address convex optimization problems. Demonstrating that the estimation problem is non-convex and non-differentiable, they approximated the loss function with some classes of differentiable functions so that gradient-based methods could be applied. \cite{Munoz_et_al:2022} considered a bilevel framework for convex optimization with both deterministic inequality and equality constraints. In various applications, they reformulated the model into a single-level, mixed-integer linear or quadratic problem, which could be solved via off-the-shelf solvers. Recently, there are studies applying implicit differentiation approaches to learn the unknown prediction model parameters \citep{Costa_Iyengar:2023, Mckenzie_et_al:2023, Sun_et_al:2022}.

An alternative methodology is the residuals-based approach popularized by \cite{Kannan_et_al:2025}, which can be found also in some earlier works, e.g., \cite{Deng_Sen:2022}. Instead of designing a new loss function, they employed existing statistical models such as the ordinary least-squares and LASSO for prediction and adopted a bootstrap method for the residuals. \cite{Kannan_et_al:2023} extended the residuals-based model by proposing a distributionally robust optimization counterpart. The advantage of the residuals-based approach is that the tractability of the optimization model is preserved, and only an additional statistical model fitting and bootstrap procedure are required. However, theoretical guarantees such as asymptotic consistency are derived under the assumption that the model class is correctly specified. Despite the fast convergence rate of such an approach, it is difficult to specify a correct functional form of the model in practice.

In contrast to the parametric approach adopted by \cite{Kannan_et_al:2025}, \cite{Bertsimas_Kallus:2020} proposed a non-parametric approximation approach by recognizing that conditional expectations could be estimated non-parametrically. Unlike the sample average approximation (SAA) method in the stochastic programming approach where scenarios are equally weighted, scenarios are re-weighted according to the proximity to the given contextual information. The weights could then be computed via non-parametric techniques such as kernel regression and $k$-nearest neighborhood. Since only scenario probabilities are adjusted, tractability of the optimization model is preserved. Hence, from the computational perspective, it is particularly interesting to explore the feasibility of extending this approach to a quantile objective, where the original problem itself is already challenging. This method has been further investigated with different variants including variance-based regularization \citep{Srivastava_et_al:2021}, robustification \citep{Bertsimas_et_al:2021, Le_et_al:2021} and in multi-period settings \citep{Bertsimas_et_al:2019}. Recently, instead of using non-parametric methods which only extract local information, \cite{Bertsimas_Koduri:2022} proposed two global approximation methods using a reproducing kernel Hilbert space approach that admits tractable reformulations. One of the methods extends the linear decision rule studied in \cite{Ban_Rudin:2019}. With the use of global approximations, they discussed and demonstrated how the curse of dimensionality in covariates encountered in the non-parametric approach could be alleviated.

Different from the above literature, we consider contextual risk-averse stochastic optimization models with a quantile objective. As pointed out in the recent survey \citep{Sadana_et_al:2025}, contextual stochastic optimization model in the risk averse setting remains as an active future research direction. However, existing approaches focus on learning the uncertainty set in robust optimization via contextual information \citep{Chenreddy_et_al:2022, Sun_et_al:2023, Wang_et_al:2023}. Recently, \cite{Tao_et_al:2025} proposed a new contextual optimization framework that considered simultaneously problem data uncertainty and contextual uncertainty. To obtain the optimal policy as a function of contextual information, they restricted the hypothesis space to a reproducing kernel Hilbert space and adopted the SAA approach. Different from \cite{Tao_et_al:2025}, we do not consider contextual uncertainty. Instead, we focus on settings where contextual information (e.g., customers' features, demographics, etc.) is known or can be accurately predicted in advance. Moreover, our framework is flexible, allowing for generic conditional quantile estimators, which includes the SAA estimator as a special case. Notably, \cite{Rahimian_Pagnoncelli:2023} and \cite{Rahimian_Pagnoncelli:2024} studied stochastic programs with contextual chance constraints and expected-value constraints, respectively. While quantile minimization problems can be cast as a chance-constrained problem, our work is different from \cite{Rahimian_Pagnoncelli:2023} in several aspects. First, while \cite{Rahimian_Pagnoncelli:2023} focused on deriving finite-sample guarantees (e.g., feasibility guarantee) for specific classes of local regression approaches, we derive general conditions under which asymptotic convergence of our proposed generic framework holds. Second, we propose a new decomposition algorithm for solving our contextual quantile minimization problem as discussed next.

\subsection{Computational Methods}
We focus on reviewing literature that investigates algorithms for solving stochastic optimization models with a quantile objective. For heuristic and exact solution approaches in the portfolio optimization context, we refer readers to \cite{Babat_et_al:2018, Cetinkaya_Thiele:2015, Feng_et_al:2015, Gaivoronski_Pflug:2005} and references therein. It is well known that quantile minimization problems can be reformulated as MIPs (see, e.g., \citealp{Benati_Rizzi:2007, Pavlikov_et_al:2017}). To improve the computational efficiency of the MIP reformulation, existing literature proposed ways to search for tight big-$M$ parameters involved and develop valid inequalities. \cite{Qiu_et_al:2014} discussed various techniques on cutting planes and branching for a $k$-violation linear program. \cite{Pavlikov_et_al:2017} focused on the MIP reformulation with a VaR objective in a generic optimization context. They discussed several valid inequalities and variable bounds that could improve the computational performance. 

Decomposition methods are alternative solution approaches to the MIP reformulations. \cite{Ivanov_Naumov:2012} studied a quantile minimization problem with a polyhedral loss (i.e., as a maximum of affine functions) and proposed a Benders' decomposition (BD) algorithm with convergence guarantee. Their proposed approach requires the knowledge of all coefficients in the polyhedral loss function, precluding its direct application to our two-stage setting. \cite{Naumov_Bobylev:2012} extended the idea of \cite{Ivanov_Naumov:2012} to study a general two-stage linear stochastic optimization problem with a quantile objective via a confidence set formulation (which results in a minimax problem). The proposed BD algorithm requires enumerating all vertices of the dual feasible set (as a polyhedron) of the recourse function. \cite{Zhenevskaya_Naumov:2018} alleviated such a disadvantage to enumerate all vertices and proposed a nested algorithm to solve the problem. However, the latter two algorithms output a feasible solution only without theoretical convergence to the true solution. To the best of our knowledge, our proposed SiCG is the first decomposition method with convergence guarantee for solving two-stage linear stochastic programs with a quantile objective based on the MIP formulation. In particular, our approach iteratively generates the unknown dual feasible vertices. Moreover, we introduce an inexactness mechanism that effectively accelerates convergence by allowing for inexact solutions to relaxations of the MIP formulation, a feature that would not be possible if the MIP formulation is solved directly as in the classical approach.

Finally, we remark that two-stage optimization problems with a quantile objective is related to chance-constrained optimization problems; see \cite{Bai_et_al:2021}, \cite{Liu_et_al:2016}, \cite{Luedtke:2014}, \cite{Zeng_et_al:2014} and so forth for state-of-the-art decomposition methods on such problems. However, these algorithms are designed for chance constraints indicating feasibility of the decisions, where the feasibility is explicitly characterized by linear constraints. Recently, \cite{Cattaruzza_et_al:2024} developed exact and heuristic solution approaches for solving a class stochastic optimization problems where the objective is a convex combination of the expected linear cost and a sum of quantiles of a linear cost. This is different from our setting where the (non-linear) random cost is characterized by a minimization problem. Thus, these algorithms could not be directly used to solve our problem of interest.

\section{Contextual Quantile Minimization Problem} \label{sec:CQM}

In this section, we first introduce the proposed risk-averse contextual stochastic optimization problem with a quantile objective, which we denote the contextual quantile minimization problem, in Section~\ref{subsec:CQM_prob} and its sample counterpart in Section~\ref{subsec:sample_quantile}.

\subsection{The Contextual Quantile Minimization Problem}  \label{subsec:CQM_prob}

Let $(\Omega,\calF,\calP)$ be a probability space. We denote $\xb\in\R^n$ and $\yb\in\R^m$ as the first-stage and second-stage decisions respectively. The random parameter in the optimization model is given by $\xib\in\Xi$, where $\Xi\subseteq\R^\ell$ is the support of $\xib$. Accordingly, we let $\calX$ be the set of first-stage feasible decisions and $\calY(\xb,\xib)=\{\yb\in\R^m \mid \Tb\xb +  \Wb\yb + \Cb\xib \geq \hb,\, \yb\geq \zero\}$ be the set of second-stage feasible decisions characterized by linear constraints. For any random variable $S$, define $Q_\tau(S)=\inf\{t \mid \Prob(S \leq t)\geq \tau\}$ as the quantile function, where $\tau\in(0,1)$. That is, $Q_\tau(S)$ computes the $(\tau\times 100\%)$th percentile of the distribution of $S$. Note that we suppress the dependence on the distribution $\Prob$ for notational simplicity. With this notation, we define the two-stage risk-averse stochastic optimization problem with a quantile objective as
\begin{equation} \label{eqn:QM}
    v^\star= \min_{\xb\in\calX} \, \, Q_\tau\big(f(\xb,\xib)\big)
\end{equation}
with the linear second-stage cost 
\begin{equation} \label{eqn:second_stage}
    f(\xb,\xib) = \min_{\yb\in\calY(\xb,\xib)} \qb^\tp \yb,
\end{equation}
where $\qb\in\R^m$ is the coefficient vector in the second-stage problem.  We assume that $f(\xb,\xib)$ is positive and has an optimal solution for any $\xb\in\calX$ and $\xib\in\Xi$, which is satisfied in practice, for example, if $f(\xb,\xib)$ represents the implementation cost of decision $\xb$. In \eqref{eqn:QM}, $\tau$ captures the risk aversion of the decision maker. Typical choices of $\tau$ would be values close to $1$ (such as $0.90$ and $0.95$) that correspond to extreme losses.



In practice, decision-makers may possess data relevant to the unknown parameter $\xib$. That is, some useful contextual information is present that provide additional knowledge on the distribution of $\xib$. Mathematically, let $\Zb$ be the random vector representing the contextual information with support $\calZ\subseteq \R^p$ and distribution $\Prob_Z$. Then, for a given contextual information $\zb\in\calZ$, we define the contextual quantile minimization problem as
\begin{equation} \label{eqn:CQM}
    v^\star(\zb)= \min_{\xb\in\calX} \, \, \Big\{Q_\tau(\xb;\zb):=Q_\tau\big(f(\xb,\xib)\, \big| \, \Zb=\zb\big)\Big\},
\end{equation}
where the conditional quantile of $f(\xb,\xib)$ given $\Zb=\zb$ is computed as $Q_\tau(f(\xb,\xib)\mid\Zb=\zb)=\inf\{t\in\R \mid \Prob(f(\xb,\xib) \leq t \mid \Zb =\zb)\geq \tau\}$. We denote $\calX^\star(\zb)$ as the set of optimal solutions to \eqref{eqn:CQM}. 


\subsection{Sample Counterpart of the Contextual Quantile Minimization Problem}  \label{subsec:sample_quantile}

In practice, the true distribution $\Prob$ of $\xib$ is unknown, but decision-makers may possess a set of historical observations of the uncertainty $\xib$ and the corresponding contextual information $\Zb$. Note that for any $\xb\in\calX$ and $\zb\in\calZ$, the quantile $Q_\tau(\xb;\zb)$ admits the following variational characterization:
\begin{equation} \label{eqn:cond_quantile_min_characterization}
    Q_\tau(\xb;\zb) = \argmin_{u\in\R} \,\, \E_{\Prob(\zb)}\Big[\rho_\tau\big(f(\xb;\xib)-u\big)\Big],
\end{equation}
where $\rho_\tau(u) = u(\tau-\one(u<0))$ is the check loss function \citep{Koenker_Hallock:2001} and $\Prob(\zb)$ is the unknown true conditional distribution of $\xib$ given $\Zb=\zb$.  If the minimizer is not unique, we pick the one with the smallest value. From \eqref{eqn:cond_quantile_min_characterization}, we observe a direct relationship between the conditional quantile and the conditional expectation. Thus, we could adopt any statistical methods that approximate conditional expectation, e.g., local regression approaches \citep{Ban_Rudin:2019, Bertsimas_Kallus:2020}, to produce an estimate of the conditional quantile. Specifically, let $\calD_N=\{(\Zb^1,\xib^1),\dots,(\Zb^N,\xib^N)\}$ be a set of samples of $(\Zb,\xib)$. With the set of sample data, we obtain a conditional quantile estimate by replacing the unknown conditional expectation in \eqref{eqn:cond_quantile_min_characterization} with an estimate:%
\begin{equation} \label{eqn:cond_quantile_min_characterization_est}
    \Qh^N_\tau(\xb;\zb) = \argmin_{u\in\R} \,\, \E_{\Probh_N(\zb)}\Big[\rho_\tau\big(f(\xb,\xib) - u\big)\Big],
\end{equation}
where $\Probh_N(\zb)$ is an approximation of the unknown true conditional distribution $\Prob(\zb)$. In Example~\ref{eg:local_regression} below, we provide an example of $\Probh_N(\zb)$ via the popular local regression approaches.  Using the estimator $\Qh^N_\tau$, we define the sample counterpart of the contextual quantile minimization problem as
\begin{equation} \label{eqn:CQM_sample}
    \vhat^N(\zb)= \min_{\xb\in\calX} \, \,\Qh^N_\tau(\xb;\zb).
\end{equation}
We denote by $\calXhat^N(\zb)$ the set of optimal solutions to \eqref{eqn:CQM_sample}.

\begin{example}[Local Regression Estimate] \label{eg:local_regression}
     Consider the local regression estimator of $\Prob(\zb)$; see, e.g., \cite{Walk:2010}. Let $\Probh_N(\zb)=\sum_{i=1}^N w_{i,N}(\zb)\delta_{\xib^i}$ for some weights $\{w_{i,N}(\zb)\}_{i=1}^N$ satisfying $w_{i,N}(\zb)\geq 0$ and $\sum_{i=1}^N w_{i,N}(\zb)=1$, where $\delta_{\xib}$ is the Dirac measure on $\xib$. The weights $\{w_{i,N}(\zb)\}_{i=1}^N$ typically captures the proximity of the new contextual information $\zb$ to the data $\Zb^i$, with a larger value of $w_{i,N}$ implying that $\zb$ is closer to $\Zb^i$. Then, problem~\eqref{eqn:cond_quantile_min_characterization_est} reduces to 
    \begin{equation} \label{eqn:cond_quantile_min_local_reg}
        \Qh^N_\tau(\xb;\zb) = \argmin_u \,\, \sum_{i=1}^N w_{i,N}(\zb) \rho_\tau\big(f(\xb,\xib^i) - u\big).
    \end{equation}
    That is, unlike the classical SAA approach where each scenario has a weight of $1/N$, the local regression method assigns a new weight to each scenario based on the proximity of $\zb$ to the data $\Zb^i$. Note that by \cite{Laksaci_et_al:2009}, 
    \begin{equation} \label{eqn:cond_quantile_est}
        \Qh^N_\tau(\xb;\zb) =\inf \Big\{ t\in\R \,\Big|\, \Fh^N(t\mid\zb) \geq \tau \Big\},
    \end{equation}
    where $\Fh^N(t\mid\zb)=\sum_{i=1}^N  w_{i,N}(\zb) \one(f(\xb,\xib^i) \leq t)$ is the estimated conditional cumulative distribution function. 
\end{example}

\section{Convergence Properties} \label{sec:convergence_CQM} 

In this section, we discuss the asymptotic convergence properties of model~\eqref{eqn:CQM_sample}. Note that our proposed problem involves two layers of minimization: finding an optimal solution $\xb\in\calX$ as in \eqref{eqn:CQM_sample} and searching for the quantile for a given $\xb\in\calX$ that corresponds to the optimal solution of another minimization problem \eqref{eqn:cond_quantile_min_characterization_est}. In other words, problem~\eqref{eqn:CQM_sample} can be viewed as a form of bilevel optimization problem. Thus, different from existing literature that considers problems with a single layer of minimization (see, e.g., \citealp{Bertsimas_Kallus:2020, Kannan_et_al:2025}), our analysis involves two layers of convergence: (a) the convergence of the sample quantile $\Qh_\tau^N$, and eventually (b) the convergence of the optimal solution and the set of optimal solutions to \eqref{eqn:CQM_sample}. We highlight that existing contextual stochastic bilevel optimization literature focuses on solving \eqref{eqn:CQM_sample} directly assuming that new samples $\xib$ can be generated through simulations, e.g., Monte Carlo techniques \citep{Bouscary_et_al:2025, Hu_et_al:2023}, and the asymptotic convergence properties of \eqref{eqn:CQM_sample} remains unclear.

In Section~\ref{subsec:convg_assumptions}, we first provide some assumptions that lay the foundation of our convergence analysis. In Section~\ref{subsec:convg_as}, we analyze the almost-sure convergence of the optimal value $\vhat^N(\zb)$ and the set of optimal solutions $\calXhat^N(\zb)$ of problem~\eqref{eqn:CQM_sample} to their true counterparts. Next, in Section~\ref{subsec:convg_prob}, we discuss conditions under which $\vhat^N(\zb)$ and  $\calXhat^N(\zb)$ converge to their counterparts in probability. In the following, we use  ``a.e.'' to denote almost everywhere, ``a.s.'' to denote almost surely, ``$\as$'' to denote almost sure convergence, and ``$\convp$'' to denote convergence in probability. For an event $A$, we define $\one_A$ as the indicator function taking value one if $A$ occurs and taking value zero otherwise. For two sets $\calS_1$ and $\calS_2$, we define the distance $D(\calS_1,\calS_2)=\sup_{\sb_1\in \calS_1}\inf_{\sb_2\in \calS_2}\norms{\sb_1-\sb_2}$.

\subsection{Assumptions}  \label{subsec:convg_assumptions}

We first impose the following regularity conditions related to the feasible set $\calX$, random cost function $f$, and its quantile.

\begin{assumption} \label{assumption:compact_X}
    The feasible set $\calX$ is compact.
\end{assumption}

\begin{assumption} \label{assumption:regularity_cond_1}
    For any $\xb\in\calX$ and $\Prob_Z$-a.e. $\zb\in\calZ$,
    \begin{enumerate}[topsep=1pt,itemsep=-1pt] 
        \item[(a)] the random variable $f(\xb,\xib)$ is integrable, i.e., $\E_{\Prob(\zb)}|f(\xb,\xib)|<\infty$;
        \item[(b)] the quantile $Q_\tau\big(f(\xb,\xib)\mid \Zb=\zb \big)<\infty$ is unique.
    \end{enumerate}
\end{assumption}

Assumption~\ref{assumption:compact_X} is standard in the stochastic optimization literature \citep{Duchi_et_al:2021, Shapiro_et_al:2014, Tsang_Shehadeh:2025}. Assumption~\ref{assumption:regularity_cond_1}(a) is a basic assumption that ensures the integrability of the random variable $f(\xb,\xib)$. Assumption~\ref{assumption:regularity_cond_1}(b) is standard in the related literature that develop convergence of conditional quantile estimates \cite[see, e.g.,][]{Laksaci_et_al:2009,Rachdi_et_al:2021}. This assumption holds, for example, when the cumulative distribution function of the random variable $f(\xb,\xib)$ is continuous.

Next, we impose additional continuity assumptions on the quantile function $Q_\tau(\xb;\zb)$ and its sample counterpart $\Qh^N_\tau(\xb;\zb)$.

\begin{assumption} \label{assumption:Holder_cont}
    The quantile function $Q_\tau(\xb;\zb)$ and its sample counterpart $\Qh^N_\tau(\xb;\zb)$ satisfy the following conditions.
    \begin{enumerate}[topsep=1pt,itemsep=-1pt] 
        \item[(a)] For $\Prob_Z$-a.e. $\zb\in\calZ$, the quantile function $Q_\tau(\xb;\zb)$ is $\alpha$-H\"{o}lder continuous on $\calX$ for some $\alpha>0$; i.e., there exists $L_1>0$ such that $\big|Q_\tau(\xb_1;\zb)-Q_\tau(\xb_2;\zb)\big|\leq L_1\norms{\xb_1-\xb_2}^\alpha$ for any $\{\xb_1,\xb_2\}\subseteq\calX$;
        \item[(b)] For $\Prob_Z$-a.e. $\zb\in\calZ$, the estimated quantile function $\Qh^N_\tau(\xb;\zb)$ is $\beta$-H\"{o}lder continuous on $\calX$ for some $\beta>0$ independent of the data $\calD_N$; i.e., there exists $L_2>0$ independent of $\calD_N$ such that $\big|\Qh^N_\tau(\xb_1;\zb)-\Qh^N_\tau(\xb_2;\zb)\big|\leq L_2\norms{\xb_1-\xb_2}^\beta$ for any $\{\xb_1,\xb_2\}\subseteq\calX$.
    \end{enumerate}
\end{assumption}

Assumption~\ref{assumption:Holder_cont} ensures that quantile function $Q_\tau(\xb;\zb)$ and its sample counterpart $\Qh^N_\tau(\xb;\zb)$, i.e., the objective functions of problems~\eqref{eqn:CQM} and \eqref{eqn:CQM_sample}, respectively, are H\"{o}lder continuous, which is standard in the existing stochastic literature when analyzing asymptotic convergence of stochastic optimization models \citep{Bertsimas_Kallus:2020, Bertsimas_McCord:2019, Pichler_Xu:2022, Sun_Xu:2016}. Assumption~\ref{assumption:Holder_cont}(a) holds under mild assumptions on the distribution $\Prob(\zb)$; we refer readers to \cite{Henrion_Romisch:2004} for detailed discussions. Assumption~\ref{assumption:Holder_cont}(b) ensures that the estimate quantile function $\Qh^N_\tau(\xb;\zb)$ is Lipschitz continuous with Lipschitz constant $L_2$  independent of the data $\calD_N$. This holds, for example, if we adopt the local regression estimation (see Example~\ref{eg:local_regression}). We first illustrate, using Example \ref{eg:Lip_quantile}, that the resulting sample quantile function $\Qh_\tau^N(\xb;\zb)$ is Lipschitz continuous, i.e., $1$-H\"{o}lder continuous.

\begin{example} \label{eg:Lip_quantile}
Consider the following two max-affine functions.
\begin{align*}
    f(x_1,x_2) &= \max\{2x_1-3x_2+3,\,3x_1-2x_2+2,\,x_1-3x_2+5\} \\
    g(x_1,x_2) &= \max\{2x_1-x_2+2,\,2x_1-2x_2+3,\,x_1-3x_2+5\}
\end{align*}
Assume that each of the function corresponds to one scenario, each with probability $0.5$.  Figure~\ref{fig:eg_quantile} shows the resulting quantile function with $\tau\leq0.5$ (i.e., the minimum of $f$ and $g$), which is non-convex and is composed of several pieces of affine functions. Figure \ref{fig:eg_partition} shows the partition of the domain where each piece of affine functions belongs to. The domain is partitioned into regions with boundaries characterized by hyperplanes (i.e., lines in $\R^2$). Note that the intersections are either caused by the intersections of affine functions from the same function or from two different functions. It is not difficult to see, from Figure \ref{fig:eg_quantile}, that the estimated quantile function is Lipschitz continuous for some large enough Lipschitz constant.
\begin{figure}
     \centering
     \begin{subfigure}[b]{0.45\textwidth}
         \centering
         \includegraphics[width=\textwidth]{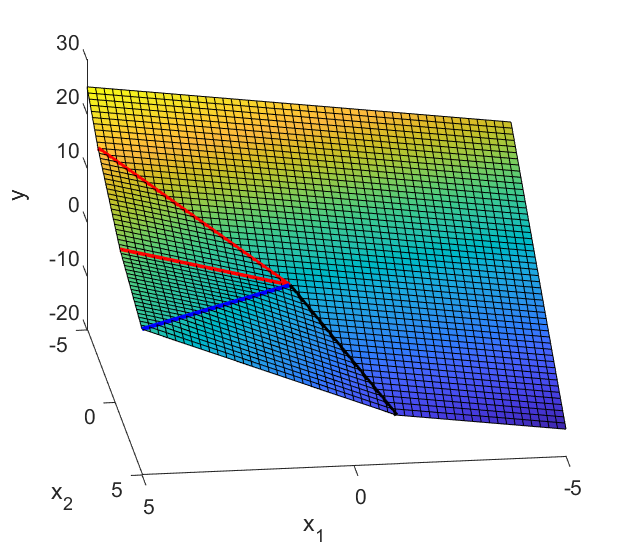}
         \caption{$\Qh^N_\tau$ with $\tau\leq0.5$}
         \label{fig:eg_quantile}
     \end{subfigure}
     \hfill
     \begin{subfigure}[b]{0.45\textwidth}
         \centering
         \includegraphics[width=\textwidth]{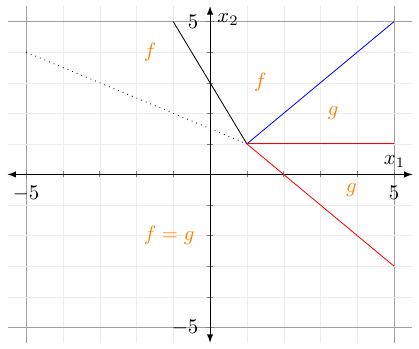}
         \caption{Area of active hyperplanes}
         \label{fig:eg_partition}
     \end{subfigure}
     \caption{Example of the estimated quantile function and its domain partition} \label{fig:eg_quantile_partition}     
\end{figure}
\end{example}

From Example \ref{eg:Lip_quantile}, we observe that $\Qh_\tau^N(\xb;\zb)$ is a piecewise affine function. Indeed, we can verify that this is true for general cost function \eqref{eqn:second_stage} when using the local regression approach. With the use of a recent result on Lipschitz continuity of piecewise Lipschitz functions by \cite{Leobacher_Steinicke:2022}, we formally prove the Lipschitz continuity of $\Qh_\tau^N(\xb;\zb)$ in Theorem \ref{thm:Lip_quantile}.

\begin{theorem} \label{thm:Lip_quantile}

Let $\Pit$ be the set of extreme points of the dual feasible set of \eqref{eqn:second_stage}. Then, for any $\zb\in\calZ$, the estimated quantile function $\Qh_\tau^N(\xb;\zb)$ using local regression approaches is Lipschitz continuous on $\calX$ with Lipschitz constant $L=\max_{\pib\in\Pit}\norms{\Tb^\tp\pib}<\infty$.
\end{theorem}

\subsection{Almost-Sure Convergence} \label{subsec:convg_as}

In this section, we establish the a.s.~convergence of the optimal value $\vhat^N(\zb)$ and the set of optimal solutions $\calXhat^N(\zb)$ of problem~\eqref{eqn:CQM_sample} to their true counterparts. To this end, in Theorem~\ref{thm:convergence_as_est_quantile}, we first prove the a.s.~convergence of the estimated quantile function $\Qh_\tau^N(\xb;\zb)$ to its true counterpart $Q_\tau(\xb;\zb)$. 

\begin{theorem} \label{thm:convergence_as_est_quantile}
    In addition to Assumption \ref{assumption:regularity_cond_1}, suppose that the following hold for $\Prob_Z$-a.e. $\zb\in\calZ$.
    \begin{itemize}[topsep=1pt,itemsep=-1pt] 
        \item [(a)] There exists a compact set $B_{\zb}\subseteq\Xi$ such that $\Prob(\xib\in B_{\zb}\mid \Zb=\zb)>0$.
    
        \item [(b)] For any $\xb\in\calX$ and $u\in\R$, we have $\E_{\Probh_N(\zb)}\big[\rho_\tau\big(f(\xb,\xib) - u\big)\big] \as \E_{\Prob(\zb)}\big[\rho_\tau\big(f(\xb,\xib) - u\big)\big]$.

        \item [(c)] For any measurable set $A\subseteq\Xi$, we have $\E_{\Probh_N(\zb)}(\one_A)\as\E_{\Prob(\zb)}(\one_A)$.
    \end{itemize}
    Then, for $\Prob_Z$-a.e. $\zb\in\calZ$, we have $\Qh^N_{\tau}(\xb;\zb)\as  Q_\tau(\xb;\zb)$ as $N\rightarrow\infty$ for any $\xb\in\calX$.
    %
\end{theorem}

We comment on the additional assumptions in Theorem~\ref{thm:convergence_as_est_quantile}. Assumption~(a) is a mild assumption that ensures that the conditional distribution $\Prob(\zb)$ has a non-zero measure on some compact subset $B_{\zb}$ of $\Xi$. This holds for most of the common probability distributions. Assumptions~(b) and (c) ensure the a.s.~convergence of the estimated quantile. They hold for popular classes of estimator of conditional expectations. We provide some examples below.

\begin{example}[Kernel Regression]
    Consider the local regression estimator discussed in Example~\ref{eg:local_regression} with weights $w_{i,N}(\zb)$ given by
    $$w_{i,N}(\zb)=\frac{K\big((\Zb^i-\zb)/h_N\big)}{\sum_{j=1}^N K\big((\Zb^j-\zb)/h_N\big)}$$
    for all $i\in\{1,\dots,N\}$ with bandwidth $h_N>0$, where $K:\R^d\rightarrow\R$ is a symmetric kernel satisfying $C_1 H(\norms{\vb})\leq K(\vb)\leq C_2H(\norms{\vb})$ for some $0<C_1<C_2<\infty$ and non-decreasing function $H:\R_+\rightarrow\R_+$ with $\lim_{t\downarrow 0} H(t)>0$ and $\lim_{t\rightarrow\infty} t^p H(t)=0$. This includes the following kernel functions: (Naive) $K(\vb)=\one(\norms{\vb}\leq 1)$, (Epanechnikov) $K(\vb)=(1-\norms{\vb}^2)\one(\norms{\vb}\leq 1)$, and (Tricubic) $K(\vb) = (1-\norms{\vb}^3)^3\one(\norms{\vb}\leq 1)$. Suppose that $\E_\Prob\big[|f(\xb,\xib)-u| \max\{ \log|f(\xb,\xib)-u|,0\}\big]<\infty$ for all $u\in\R$. This holds, for example, if $\Var(f(\xb,\xib))<\infty$, and in particular, if $f(\xb,\xib)$ is a bounded random variable. Such a finite variance or boundedness assumption holds in many application domains such as facility location, scheduling, and unit commitment problems; see \cite{Jiang_et_al:2019, Wang_et_al:2020, Zhou_et_al:2019}, to name but a few. It then follows from Theorem~3 of \cite{Walk:2010} that assumptions~(b) and (c) hold if one of the following conditions is satisfied.
    \begin{itemize}[leftmargin=5mm]
        \item The data $\calD_N=\{(\Zb^i,\xib^i)\}_{i=1}^N$ comes from a $\rho$-mixing process with $\rho(k)=O(k^{-\gamma})$ for some $\gamma>0$ with bandwidth $h_N=CN^{-\delta}$ for some $C>0$ and $0<\delta<2\gamma/(p+2p\gamma)$. This includes the special case with $\rho=0$ where the data $\calD_N$ are independent and identically distributed (i.i.d.).

        \item The data $\calD_N=\{(\Zb^i,\xib^i)\}_{i=1}^N$ comes from a $\alpha$-mixing process with $\alpha(k)=O(k^{-\gamma})$ for some $\gamma>1$ with bandwidth $h_N=CN^{-\delta}$ for some $C>0$ and $0<\delta<2(\gamma-1)/(3p+2p\gamma)$. 
    \end{itemize}
\end{example}

\begin{example}[$k$NN Regression]
    Consider the local regression estimator discussed in Example~\ref{eg:local_regression} with weights $w_{i,N}(\zb)$ given by
    $$w_{i,N}(\zb) = \frac{1}{k}\one(\Zb^i \text{ is a $k$NN of } \zb),$$
    where $k=k_N=\min\big\{ \lceil CN^\delta \rceil, N-1 \big\}$ for some $C>0$ and $\delta\in(0,1)$. Here, $\Zb^i$ is a $k$NN of $\zb$ if $i\in N_k(\zb)=\big\{i'\in\{1,\dots,N\}\mid \sum_{j=1}^N \one(\norms{\zb-\zb_{i'}}\geq \norms{\zb-\zb_j}\}) \leq k \big\}$. Then, it follows from Theorem~5 of \cite{Walk:2010} that assumptions (b) and (c) hold if the data $\calD_N=\{(\Zb^i,\xib^i)\}_{i=1}^N$ are i.i.d.
\end{example}

A key step toward establishing the desired a.s.~convergence of the optimal value and the set of optimal solutions is the uniform a.s.~convergence of estimated quantile function $\Qh_\tau^N(\xb;\zb)$ to its true counterpart $Q_\tau(\xb;\zb)$ over $\xb\in\calX$. Lemma~\ref{lem:unif_convergence_as_obj} below establishes such a uniform convergence result.

\begin{lemma} \label{lem:unif_convergence_as_obj}
    Suppose that Assumptions~\ref{assumption:compact_X}--\ref{assumption:Holder_cont} and the assumptions in Theorem~\ref{thm:convergence_as_est_quantile} hold. Then, for $\Prob_Z$-a.e.~$\zb\in\calZ$, we have $\Qh_\tau^N(\xb;\zb)\as Q_\tau(\xb;\zb)$ uniformly over $\xb\in\calX$.
\end{lemma}

With Lemma~\ref{lem:unif_convergence_as_obj}, we are now ready to prove the a.s.~convergence of our problem~\eqref{eqn:CQM_sample}. Theorem~\ref{thm:convergence_as_opt_val_sol} establishes that, under the stated assumptions, the optimal value $\vhat^N(\zb)$ and the set of optimal solutions $\calXhat^N(\zb)$ to problem~\eqref{eqn:CQM_sample} converge a.s.~to the true optimal value $v^\star(\zb)$ and the set of optimal solutions $\calX^\star(\zb)$ to \eqref{eqn:CQM}, respectively.

\begin{theorem} \label{thm:convergence_as_opt_val_sol}
    Under the same assumptions in Lemma~\ref{lem:unif_convergence_as_obj}, for $\Prob_Z$-a.e.~$\zb\in\calZ$, as $N\rightarrow\infty$, we have (i)  $\vhat^N(\zb) \as  v^\star(\zb)$ (ii) $\sup_{\xb\in\calXhat^N(\zb)} Q_\tau(\xb^N;\zb)  \as v^\star(\zb)$, and (iii) $D(\calXhat^N(\zb),\calX^\star(\zb))\as 0$. 
\end{theorem}

\subsection{Convergence in Probability} \label{subsec:convg_prob}

In this section, we derive the convergence of the optimal value $\vhat^N(\zb)$ and the set of optimal solutions $\calXhat^N(\zb)$ of problem~\eqref{eqn:CQM_sample} to their true counterparts in probability under a weaker set of assumptions. As pointed out in \cite{Kannan_et_al:2025}, this is useful in part because statistics literature typically focus on establishing the convergence in probability of the estimated conditional distribution $\Probh_N(\zb)$. This leads to a set of conditions that could be more readily verified compared with those in Section~\ref{subsec:convg_as}. 

Again, we start by establishing the convergence of the estimated quantile function $\Qh_\tau^N(\xb;\zb)$ to its true counterpart $Q_\tau(\xb;\zb)$ in probability. 

\begin{theorem} \label{thm:convergence_prob_est_quantile}
    In addition to Assumption \ref{assumption:regularity_cond_1}, suppose that the following hold for $\Prob_Z$-a.e.~$\zb\in\calZ$.
    \begin{itemize}[topsep=1pt,itemsep=-1pt] 
        \item [(a)] There exists a compact set $B_{\zb}\subseteq\Xi$ such that $\Prob(\xib\in B_{\zb}\mid \Zb=\zb)>0$.
    
        \item [(b)] For any $\xb\in\calX$ and $u\in\R$, we have $\E_{\Probh_N(\zb)}\big[\rho_\tau\big(f(\xb,\xib) - u\big)\big] \convp \E_{\Prob(\zb)}\big[\rho_\tau\big(f(\xb,\xib) - u\big)\big]$.

        \item [(c)] For any measurable set $A\subseteq\Xi$, we have $\E_{\Probh_N(\zb)}(\one_A)\convp\E_{\Prob(\zb)}(\one_A)$.
    \end{itemize}
    Then, for $\Prob_Z$-a.e. $\zb\in\calZ$, we have $\Qh^N_{\tau}(\xb;\zb)\convp Q_\tau(\xb;\zb)$ as $N\rightarrow\infty$ for any $\xb\in\calX$.
    %
\end{theorem}

\begin{remark}
    Different from the existing works that analyze the convergence in probability of these models under compactness assumption (see, e.g., \citealp{Bertsimas_McCord:2019, Kannan_et_al:2025}), the estimated quantile $\Qh^N_\tau(\xb;\zb)$ and the true quantile $Q_\tau(\xb;\zb)$ respectively correspond to optimization models \eqref{eqn:cond_quantile_min_characterization} and \eqref{eqn:cond_quantile_min_characterization_est} with an unbounded feasible set $\R$. Therefore, our analysis extends those in the existing literature by relaxing the compactness assumption and replacing it with some (weaker) uniform coerciveness conditions of the objective function, which is $\rho_\tau$ in our setting.    
\end{remark}

Next, in Lemma~\ref{lem:unif_convergence_prob_obj}, we prove the uniform convergence of the estimated quantile function $\Qh^N_\tau(\xb;zb)$ to its true counterpart $Q_\tau(\xb;\zb)$ over $\xb\in\calX$. Our results extend the one in \cite{Bertsimas_McCord:2019} by considering a more general H\"{o}lder continuity assumption; see Assumption~\ref{assumption:Holder_cont}.

\begin{lemma} \label{lem:unif_convergence_prob_obj}
    Suppose that Assumptions~\ref{assumption:compact_X}--\ref{assumption:Holder_cont} and the assumptions in Theorem~\ref{thm:convergence_prob_est_quantile} hold. Then, for $\Prob_Z$-a.e.~$\zb\in\calZ$, we have $\Qh_\tau^N(\xb;\zb)\convp Q_\tau(\xb;\zb)$ uniformly over $\xb\in\calX$.
\end{lemma}

With Lemma~\ref{lem:unif_convergence_prob_obj}, we are now ready to prove the desired convergence results for problem~\eqref{eqn:CQM_sample}. Theorem~\ref{thm:convergence_prob_opt_val_sol} establishes that, under the stated assumptions, the optimal value $\vhat^N(\zb)$ and the set of optimal solutions $\calXhat^N(\zb)$ to problem~\eqref{eqn:CQM_sample} converge to the true optimal value $v^\star(\zb)$ and the set of optimal solutions $\calX^\star(\zb)$ to \eqref{eqn:CQM}, respectively, in probability.

\begin{theorem} \label{thm:convergence_prob_opt_val_sol}
    Under the same assumptions in Lemma~\ref{lem:unif_convergence_prob_obj}, for $\Prob_Z$-a.e.~$\zb\in\calZ$, as $N\rightarrow\infty$, we have (i)  $\vhat^N(\zb) \convp  v^\star(\zb)$ (ii) $\sup_{\xb\in\calXhat^N(\zb)} Q_\tau(\xb^N;\zb)  \convp v^\star(\zb)$, and (iii) $D(\calXhat^N(\zb),\calX^\star(\zb))\convp 0$. 
\end{theorem}

\section{Stochastic Inexact Constraint Generation Method} \label{sec:SiCG}

Note that the solvability of \eqref{eqn:CQM_sample} depends on the approach adopted to estimate the conditional quantile. In this section, we discuss solution approaches that are designed to search for the optimal solution to the quantile minimization problem~\eqref{eqn:CQM_sample} with the local regression estimator that often leads to tractable reformulations; see Example~\ref{eg:local_regression}.
\begin{equation} \label{prob:quantile_min}
    \vhat^N(\zb) = \min_{\xb\in\calX} \Bigg\{ \Qh_\tau^N(\xb;\zb)=\inf \bigg\{ t\in\R \,\bigg|\, \sum_{i=1}^N  w_{i,N}(\zb) \one(f(\xb,\xib^i) \leq t) \geq \tau \bigg\}\Bigg\}.
\end{equation}%
Nevertheless, our proposed algorithm can be applied to solve any quantile minimization problem with a discrete distribution. Since we focus on solving \eqref{prob:quantile_min} with a given predictor $\zb$ and a historical data set, for notational simplicity, we let $w_i=w_{i,N}(\zb)>0$ for all $i\in\{1,\dots,N\}$ to suppress the dependence on $\zb$. In particular, if $w_i=1/N$ for all $i\in\{1,\dots,N\}$, then \eqref{prob:quantile_min} reduces to the classical SAA of the two-stage quantile minimization problem without any contextual information, and our proposed algorithm can also be applied. In Section \ref{subsec:sol_MIP_reformulation}, we discuss a direct MIP reformulation of \eqref{prob:quantile_min} as a standard approach to quantile minimization problems. {Recognizing that such a reformulation, with a potentially large number of additional variables and constraints, is challenging to solve,} we introduce our proposed decomposition method and investigate its convergence properties in Section~\ref{subsec:sol_SiCG}.

\subsection{MILP reformulation} \label{subsec:sol_MIP_reformulation}

First, in view of \eqref{eqn:cond_quantile_est}, we introduce binary variable $z_i$ taking value $1$ if $f(\xb,\xib^i)\leq t$, and $0$ otherwise, for $i\in\{1,\dots,N\}$. Then, \eqref{prob:quantile_min} is equivalent to
\begin{subequations} 
\begin{align}
\underset{\xb\in\calX,\,t\in\R,\,\vb\in\{0,1\}^N}{\text{minimize}\,} \quad
&  t     \label{model:MILP_quan_min_obj} \\    
\text{subject to} \hspace{5mm} \quad
&  M(1-v_i) \geq f(\xb,\xib_i)-t,\quad\forall i\in\{1,\dots,N\},     \label{model:MILP_quan_min_con1} \\    
&  \textstyle\sum_{i=1}^N w_i v_i \geq \tau,   \label{model:MILP_quan_min_con2}
\end{align} \label{model:MILP_quan_min}%
\end{subequations} 
where $M$ is a large enough constant. Indeed, when $f(\xb,\xib^i)> t$, constraints \eqref{model:MILP_quan_min_con1} ensures that $v_i=0$; otherwise, $v_i$ is unrestricted by constraints \eqref{model:MILP_quan_min_con1}. However, constraints \eqref{model:MILP_quan_min_con2} and the minimization nature of our objective function ensure that $v_i=1$ in the latter case.  As pointed out in \cite{Pavlikov_et_al:2017}, the choice 
$$M=\max_{\xb\in\calX} \max_{i=1,\dots,N} f(\xb,\xib^i)-\min_{\xb\in\calX} \min_{i=1,\dots,N} f(\xb,\xib^i)$$
is sufficient to ensure that \eqref{prob:quantile_min} and \eqref{model:MILP_quan_min} are equivalent. We note that finding a tight big-$M$ parameter depends on the actual application, which is not the main focus of this paper. Finally, using \eqref{eqn:second_stage}, we can reformulate \eqref{model:MILP_quan_min} into a mixed-integer linear program (MILP) shown in Proposition \ref{prop:MILP_quan_min_LP_form}.

\begin{proposition} \label{prop:MILP_quan_min_LP_form}
Problem \eqref{model:MILP_quan_min} is equivalent to
\allowdisplaybreaks
\begin{subequations} 
\begin{align}
\underset{\xb\in\calX,\,t\in\R,\,\vb\in\{0,1\}^N,\,\yb^i\in\R^N_+}{\textup{minimize}\,} \quad
&  t     \label{model:MILP_quan_min_LP_form_obj} \\    
\textup{subject to} \hspace{9.5mm} \quad
&  M(1-v_i) \geq \qb^\tp \yb^i -t,\quad\forall i\in\{1,\dots,N\},     \label{model:MILP_quan_min_LP_form_con1} \\
&  \Tb\xb + \Wb\yb^i + \Cb\xib^i \geq \hb,\quad\forall i\in\{1,\dots,N\},     \label{model:MILP_quan_min_LP_form_con2} \\
&  \textstyle\sum_{i=1}^N w_i v_i \geq \tau,   \label{model:MILP_quan_min_LP_form_con3} 
\end{align} \label{model:MILP_quan_min_LP_form}%
\end{subequations} 
and for any optimal solution $(\xb^\star,t^\star,\vb^\star,\yb^\star)$ to \eqref{model:MILP_quan_min_LP_form},  $(\xb^\star,t^\star,\vb^\star)$ is an optimal solution to \eqref{model:MILP_quan_min}.
\end{proposition}

\subsection{Stochastic Inexact Constraint Generation (SiCG)} \label{subsec:sol_SiCG}

The direct MILP reformulation in Proposition \ref{prop:MILP_quan_min_LP_form} involves $mN$ additional variables $\{\yb_i\}_{i=1}^N$ and $d_h N$ additional constraints in \eqref{model:MILP_quan_min_LP_form_con2}, where $d_h$ is the dimension of $\hb$. Although the model size scales linearly with $N$, the MILP \eqref{model:MILP_quan_min_LP_form} remains challenging to solve when $N$ is fairly large. In this section, we propose a novel stochastic inexact constraint generation algorithm, which we call SiCG, to solve \eqref{model:MILP_quan_min}. Unlike \eqref{model:MILP_quan_min_LP_form} that considers a primal perspective (i.e., with primal variables and constraints of $f(\xb,\xib)$), our proposed SiCG considers a dual perspective as detailed below.

Note that we can rewrite the second-stage cost using its dual:
\begin{equation} \label{eqn:second_stage_dual}
    f(\xb,\xib)=\max_{\pib}\Big\{ (\hb-\Tb\xb-\Cb\xib)^\tp\pib \,\Big|\, \Wb^\tp \pib\leq \qb ,\, \pib\geq\zero \Big\}.
\end{equation}
Denote the set of dual feasible solution as $\Pi=\{\pib\mid \Wb^\tp\pib\leq\qb,\, \pib\geq 0\}$. Then, \eqref{model:MILP_quan_min} is equivalent to
\begin{subequations} 
\begin{align}
\underset{\xb\in\calX,\,t\in\R,\,\vb\in\{0,1\}^N}{\text{minimize}\,} \quad
&  t     \label{model:MILP_quan_min_dual_obj} \\    
\text{subject to} \hspace{5mm} \quad
&  M(1-v_i) \geq (\hb-\Tb\xb-\Cb\xib^i)^\tp\pib-t,\quad\forall \pib\in\Pi,\,i\in\{1,\dots,N\},     \label{model:MILP_quan_min_dual_con1} \\    
&  \textstyle\sum_{i=1}^N w_i v_i \geq \tau.   \label{model:MILP_quan_min_dual_con2}
\end{align} \label{model:MILP_quan_min_dual}%
\end{subequations} 
Note that \eqref{model:MILP_quan_min_dual} does not require any additional variables as opposed to \eqref{model:MILP_quan_min_LP_form}. However, constraints \eqref{model:MILP_quan_min_dual_con1} involve an infinite number of dual variables $\pib\in\Pi$. As a result, \eqref{model:MILP_quan_min_dual} is a semi-infinite program that cannot be solved directly. This motivates us to propose a constraint-generation based algorithm that solves \eqref{model:MILP_quan_min_dual} with a gradually increasing number of constraints \eqref{model:MILP_quan_min_dual_con1} with respect to some $\pib\in\Pi$.

The proposed SiCG algorithm is different from classical constraint generation algorithms (see, e.g., \citealp{Jiang_et_al:2012, Minguez_et_al:2021, Thiele_et_al:2009}) in the following two aspects. First, instead of requiring master problems (described next) to be solved to optimality, master problems are solved to an adaptively updated MIP gap or time limit (see also Remark \ref{rem:time_limit}). Such a feature is especially attractive in quantile minimization problems since the number of binary variables scale with the number of scenarios $N$ (see \eqref{model:MILP_quan_min_LP_form}). To ensure the convergence of the algorithm, we employ a backtracking routine proposed in \cite{Tsang_et_al:2023} to adaptively reduce the inexactness (i.e., the MIP gap) in solving master problems. We note that in \cite{Tsang_et_al:2023}, they considered a two-stage robust optimization problem, which is different from our quantile minimization problem. In particular, unlike robust optimization problems where a newscnario can be identified by solving a subproblem, a different procedure is needed to identify new dual solutions in our quantile minimization problem. Such a difference also leads to our second differentiation point with classical constraints generation algorithms, which we elaborate next.

Second, when searching for a new dual solution $\pib\in\Pi$, one may identify a set of potential new dual solutions, which may not be a singleton. This is mainly due to the use of quantile objective, which does not naturally correspond to a subproblem as in the robust optimization context, i.e., maximizing $Q(\xb,\xib)$ over a set of $\xib\in\Xi$ (see, e.g., \citealp{Jiang_et_al:2012,Thiele_et_al:2009}). As a result, we need to design a new subproblem for identifying new dual solutions. Specifically, instead of deterministically select (or compute) a new dual solution, we randomly select one new dual solution from a set of admissible dual solutions at each iteration. This distinguishes our proposed SiCG from classical constraint generation methods.

Algorithm \ref{algo:SiCG} summarizes the three steps of SiCG. First, in step 1, we solve the master problem \eqref{model:MP_quan_min_SiCG}, which is essentially \eqref{model:MILP_quan_min_dual} with only a finite subset $\Pih$ of the dual variables in $\Pi$ to a MIP gap of $\varepsilon_{MP}^j$ at iteration $j$. Note that \eqref{model:MP_quan_min_SiCG_con3} enforces a lower bound $\Lbar$ on the objective value, which would be useful to accelerate the lower bound improvement (see discussions in \citealp{Tsang_et_al:2023}). A lower bound $L^j\geq \Lbar$ and an upper bound $U^j$ from solving \eqref{model:MP_quan_min_SiCG} are recorded.  Moreover, we update $\ell$ that keeps track of the validity of the lower bound $L^j$ to the optimal value $\vhat^N(\zb)$ in \eqref{prob:quantile_min}, and the value of $\Lbar$ that is used as a lower bound of the next master problem. Next, in step 2, given a feasible solution $\xb^j\in\calX$, we evaluate the true objective value $\Qh_\tau^N(\xb^j;\zb)$. Thus, this provides an upper bound to the optimal value $\vhat^N(\zb)$ in \eqref{prob:quantile_min}.

Finally, in step 3, depending on the value of the \textit{inexact} gap $(\Ubar-U^j)/\Ubar$ (since $U^j$ may not be a valid lower bound), we proceed to either the exploitation or the exploration step. When the inexact gap is large, i.e., $\upsilon^\star_j$ from master problem with current $\Pih$ provides a loose lower bound to $\vhat^N(\zb)$, we proceed to the exploration step. In this step, as in classical constraint generation methods, we seek for a new dual variable and enlarge the set $\Pih$, i.e., the subset of dual variables used in the master problem. To obtain a new dual solution, we solve the dual of $f(\xb^j,\xib^i)$ and record its optimal solution $\pib^\star_i$  for all $i\in\{1,\dots,N\}$. Then, we randomly select $\pib^\star$ from $\{\pib^\star_i\}_{i=1}^N$ such that $\pib^\star$ is not in the current $\Pih$ and enlarge $\Pih$ with this new dual solution. We assume that each of these admissible dual solutions has a positive probability of being selected. As we show in the proof of Theorem~\ref{thm:SiCG_convergence}, we can always find such a new dual solution $\pib^\star$. Our convergence result does not depend on the randomization, as long as a new dual solution is added in the exploration step. For example, one could uniformly select $\pi^\star$ from the set of potential dual solutions. On the other hand, if the inexact gap is small, instead of enlarging the set of dual solutions, we proceed to the exploitation step. In this step, we restore the lower bound of the next master problem as the most current valid lower bound $L^\ell$. Then, we reduce the MIP gap of the master problem as a way to reduce the inexactness from solving the subsequent master problems. We then return to step~1 to re-solve the master problem, and the algorithm terminates when the actual gap $(\Ubar-L^\ell)/\Ubar$ is within the pre-specified tolerance $\varepsilon$.
\LinesNumbered  
\IncMargin{1em}
\begin{algorithm}[t!]
\ra{0.5}
\DontPrintSemicolon\small 
\SetKwInOut{Initialization}{Initialization}
\Initialization{Set $\Lbar=0$, $\Ubar=\infty$, $\varepsilon>0$, $\epst\in(0,\varepsilon/(1+\varepsilon))$,  $\{\varepsilon_{MP}^j>0\}_{j\in\N}$, $\alpha<1$, $\Pih=\emptyset$, $j=1$, $\ell=0$, $L^0=0$.}
\While{$(\Ubar-L^\ell)/\Ubar > \varepsilon$}{
\textbf{Step 1: Master Problem}\\
Solve the following master problem to a relative MIP gap of $\varepsilon_{MP}^j$:
\begin{subequations} 
\begin{align}
v^\star_j = \underset{\xb\in\calX,\,t\in\R,\,\vb\in\{0,1\}^N}{\text{minimize}\,} \quad
&  t     \label{model:MP_quan_min_obj_SiCG} \\    
\text{subject to} \hspace{5mm}\quad
&  M(1-v_i) \geq \pib^\tp(\hb-\Tb\xb-\Cb\xib^i)-t,\quad\forall \pib\in\Pih,\, i\in\{1,\dots,N\},     \label{model:MP_quan_min_SiCG_con1} \\    
&  \textstyle\sum_{i=1}^N w_i v_i \geq \tau,   \label{model:MP_quan_min_SiCG_con2}     \\
&  t \geq \Lbar   \label{model:MP_quan_min_SiCG_con3}    
\end{align} \label{model:MP_quan_min_SiCG}%
\end{subequations} 
\hspace{-2mm}and obtain the best feasible solution $(\xb^j,t^j,\vb^j)$. \\
Record a lower bound $L^j\geq \Lbar$ and an upper bound $U^j=t^j$ of $v^\star_j$. If $L^j > \Lbar$, set $\ell =j$.\\
Update $\Lbar \leftarrow U^j$. \\
\BlankLine
\textbf{Step 2: Subproblem} \\
Evaluate $\Qh_\tau^N(\xb^j;\zb)$ as in \eqref{eqn:cond_quantile_est} and update $\Ubar \leftarrow \min\{\Ubar,\, \Qh_\tau^N(\xb^j;\zb)\}$. \\
\BlankLine
\textbf{Step 3: Backtracking -- Exploitation or Exploration Routine} \\
\eIf {$(\Ubar-U^j)/\Ubar < \epst$}  {   
    Set $j=\ell$ and $\Lbar=L^\ell$. Update $\varepsilon_{MP}^j \leftarrow \alpha \varepsilon_{MP}^j$ for $j\geq \ell$. \tcc*[r]{Exploitation} 
} { 
    Obtain $\pib_i$, dual solutions of $f(\xb^j,\xib^i)$, for all $i\in\{1,\dots,N\}$. \tcc*[r]{Exploration} 
    Randomly select $\pib^\star\in\big\{\pib_i,\, i\in\{1,\dots,N\} \mid \pib_i\not\in\Pih\big\}$.\\
    Enlarge the scenario set $\Pih \leftarrow \Pih\cup\{\pib^\star\}$. \\
}
}
Return the solution with best objective $\Ubar$.
\BlankLine
\caption{The SiCG method for the quantile minimization problem}\label{algo:SiCG}
\end{algorithm}\DecMargin{1em}


In Theorem \ref{thm:SiCG_convergence}, we provide the a.s. finite convergence of SiCG. We remark that modifications in the convergence proof of a related algorithm in \cite{Tsang_et_al:2023} are needed to accommodate the difference in the studied problem (i.e., the use of quantile objective) and the stochasticity involved in the exploration step.
\begin{theorem} \label{thm:SiCG_convergence}
If $\epst<\varepsilon/(1+\varepsilon)$, then SiCG terminates in a finite number of iterations almost surely (w.r.t. the probability measure of selecting a new dual solution $\pib^\star$).
\end{theorem}

\begin{remark} \label{rem:time_limit}
Similar to \cite{Tsang_et_al:2023}, we can consider a time-limit variant of SiCG that imposes a time limit $\kappa$ when solving the master problem. Then, at each exploitation step, we update the time limit $\kappa\leftarrow \kappa+\beta$ for some $\beta>0$. This could be useful if the initial master problems are challenging but the lower bounds obtained are loose. We also implemented this variant in our numerical experiments.
\end{remark}

\begin{remark} 
Although motivated by the two-stage setting, where the random cost function is the optimal value of a linear program, our proposed SiCG method could be readily extended to solve problems with a general decomposable cost function of the form $f(\xb,\xib)=\sup_{\pib\in\Pi}\{h(\xb,\xib;\pib)\}$. For example, if $f(\xb,\xib)$ is a convex in $\xb$, we can write $f(\xb,\xib)=\max_{\pib\in\calX}\{f(\pib,\xib)+g(\pib,\xib)^\tp(\xb-\pib)\}$, where $g(\pib,\xib)$ is a subgradient of $f(\cdot,\xib)$ at $\pib\in\calX$. One could then replace the term $\pib^\tp(\hb-\Tb\xb-\Cb\xib^i)$ in \eqref{model:MP_quan_min_SiCG_con2} of the master problem~\eqref{model:MP_quan_min_SiCG} with $h(\xb,\xib;\pib)$. We leave further theoretical and numerical investigations on such generic cost functions for future work.
\end{remark}

\section{Numerical Experiments} \label{sec:numerical_expt}

In this section, we conduct numerical experiments on a single-server appointment scheduling problem (ASP). In this problem,  a set of appointments needs to be scheduled on a single server.  The objective is to determine the start time for each appointment such that the costs of appointment waiting times, as well as server idle time and overtime, are minimized. We describe the problem and experiment settings in Section \ref{subsec:expt_setting}. Then, in Section \ref{subsec:sol_time}, we compare the solution times of our proposed SiCG with those of the MILP reformulation.  Next, in Section \ref{subsec:opt_schedule}, we study the optimal schedule from classical ASP models where contextual information is not considered or a mean objective is employed. Finally, in Section \ref{subsec:OS_performance}, we provide managerial insights via out-of-sample simulations and discuss when significant gains could be obtained from our proposed model.

\subsection{Experiment Setting} \label{subsec:expt_setting}

We consider $n$ appointments that arrive at a single server with a fixed order. The server then processes each appointment in order upon arrival, and the service duration $s_i$ of appointment $i$ is random. The goal is to decide the amount of time $x_i$ for appointment $i$ that minimizes an objective comprising of costs associated with the server's idle time and overtime and appointment waiting time. The server has a working hour of $T$, and any work beyond $T$ is considered as overtime. 

As in classical ASP, the set of first-stage feasible decision is
$$\calX=\bigg\{\xb\in\R^n\,\bigg|\, \sum_{i=1}^n x_i = T,\, x_i\geq 0\bigg\},$$
where we require the sum of appointment times $x_i$ to be $T$. Let $u_i$ and $w_i$ be the idle and waiting time associated with appointment $i$ for $i\in\{1,\dots,n\}$ respectively, and $w_{n+1}$ be the overtime. Denote $\cu_i$, $\cw_i$, and $\co$ be the per unit penalty of idle time, waiting time, and overtime respectively. Under the assumption that $\cu_{i+1}-\cu_i \leq \cw_{i+1}$ for all $i\in\{1,\dots,n-1\}$ \citep{Ge_et_al:2014, Jiang_et_al:2017}, we can compute the weighted sum of idling, waiting, and overtime cost as a linear program (LP):
\begin{subequations} 
\begin{align}
f(\xb,\sb)=\underset{\wb,\,\ub}{\text{minimize}\,} \quad
&  \sum_{i=1}^n (\cu_i u_i + \cw_i w_i ) + \co w_{n+1}     \label{model:appt_sch_obj} \\    
\text{subject to} \quad
&  w_i-u_{i-1}=s_{i-1}+w_{i-1}-x_{i-1},\quad\forall i\in\{2,\dots,n+1\},     \label{model:appt_sch_con1} \\    
&  w_1 = 0,  \label{model:appt_sch_con2}    \\
&  w_i \geq 0,\, u_i\geq 0,\, w_{n+1}\geq 0,\quad\forall i\in\{1,\dots,n\}.  \label{model:appt_sch_con3}    
\end{align} \label{model:ASP}%
\end{subequations} 
Objective \eqref{model:appt_sch_obj} is to minimize a linear weighted sum of idling, waiting and overtime costs. Constraints \eqref{model:appt_sch_con1} yield the waiting, overtime or idle time. Since the first appointment always arrives at time $0$, constraint \eqref{model:appt_sch_con2} ensures that the associated waiting time is also zero.

In Section \ref{subsec:sol_time}, to compare the computational performance between our proposed SiCG and the MILP reformulation, we focus on solving the generated sample average approximation (SAA) instances, i.e., solving the ASP with a quantile objective using the empirical distribution from the generated data of size $N$ (see, e.g., \citealp{Shapiro_et_al:2014}).  We generate the data based on settings employed in the existing literature (see, e.g., \citealp{Jiang_et_al:2017, Mak_et_al:2015}, and references therein). Specifically, each appointment follows a lognormal distribution with mean $\mu_i\sim U[36,44]$ (where $U[a,b]$ represents a uniform distribution on $[a,b]$) and standard deviation $\sigma=\nu\mu_i$. Note that $\nu$ controls the variability of individual service duration. Next, we set the working hour $T=\sum_{i=1}^n \mu_i + R\sqrt{\sum_{i=1}^n \sigma_i^2}$, where $R$ adjusts the length of the server's regular working hour. In our experiments, we consider $\nu\in\{0.2,0.5\}$ and $R\in\{0.5,1.0\}$ that represent different extent of service duration variability and length of the server's working hour. For the cost structure, we set $\cu:\cw:\co=0.5:1:10$. 

In Sections \ref{subsec:opt_schedule}--\ref{subsec:OS_performance}, to derive insights of our proposed methodology, we generate the instances under the presence of contextual information. In particular, each appointment is associated with a given characteristic $z\sim U[-15,15]$. The service duration follows a lognormal distribution with mean $\mu+z$ and $\sigma=\nu\mu$, where we fix $\mu=40$. Hence, the given characteristic would provide additional information on the actual service duration distribution. The corresponding working hour is set to $T=n\mu + R\sqrt{n}\sigma$. For illustrative purposes, we focus on the case where there are $n=6$ appointments arriving at the server. We investigate the operational performances of our proposed methodology using three different characteristics (i.e., predictors) of these $6$ appointments: (a) $(0,0,0,0,0,0)^\tp$, (b) $(-15,-9,-3,3,9,15)^\tp$, and (c) $(15,9,3,-3,-9,-15)^\tp$. Case (a) corresponds to the case where all service durations follow the same distribution with mean $40$, while cases (b) and (c) correspond to increasing and decreasing means of service duration, respectively. {While ASP with a quantile objective or contextual information under an expectation objective has been investigated separately \citep{Sadghiani_et_al:2021, Sang_et_al:2021}, and our work is the first to combine these two elements for ASP and to study the potential benefits.}

In our experiments, we implement a tailored SiCG based on the problem structure of \eqref{model:ASP}. Specifically, instead of solving $f(\xb^j,\xib^i)$ for dual solutions $\pib_i$ in the exploration, we can generate a vertex from the dual feasible set of \eqref{model:ASP} in a less computationally expensive manner by leveraging the structure of the dual of \eqref{model:ASP}. We provide the details of the implemented SiCG in \ref{appdx:SiCG_ASP}. All the experiments are conducted on a computer with AMD Ryzen 7 5800X 3.8 GHz CPU and 16 Gb memory.

\subsection{Computational Time} \label{subsec:sol_time}
In this section, we compare the computational time between our proposed SiCG and the MILP reformulation, where we denote the latter as MILP in the following. For SiCG parameters, with reference to \cite{Tsang_et_al:2023}, we set $\epst=1.5\%$, $\alpha=0.5$, and initialize the master problem MIP gap as $\varepsilon_{MP}=5\%$. Moreover, we impose a time limit of $\kappa=30$s to the master problem, subject to an increase of $\beta=60$s in the exploitation step (see Remark \ref{rem:time_limit}). We set the termination criterion as a relative gap of $\varepsilon=2\%$ for both SiCG and MILP. For a given number of appointments $n\in\{6,8\}$ and a set of parameters $(v,R)$, we generate $30$ sets of scenarios with size $N\in\{200,500,1000\}$ as described in Section \ref{subsec:expt_setting}. To investigate the computational performance of the algorithms, we consider two different quantile levels $\tau\in\{0.95,0.90\}$. For presentation brevity, we focus on $n=6$ and relegate similar results for $n=8$ to \ref{appdx:expt_compt_time}.

\begin{figure}[t] 
\centering
\includegraphics[scale=0.75]{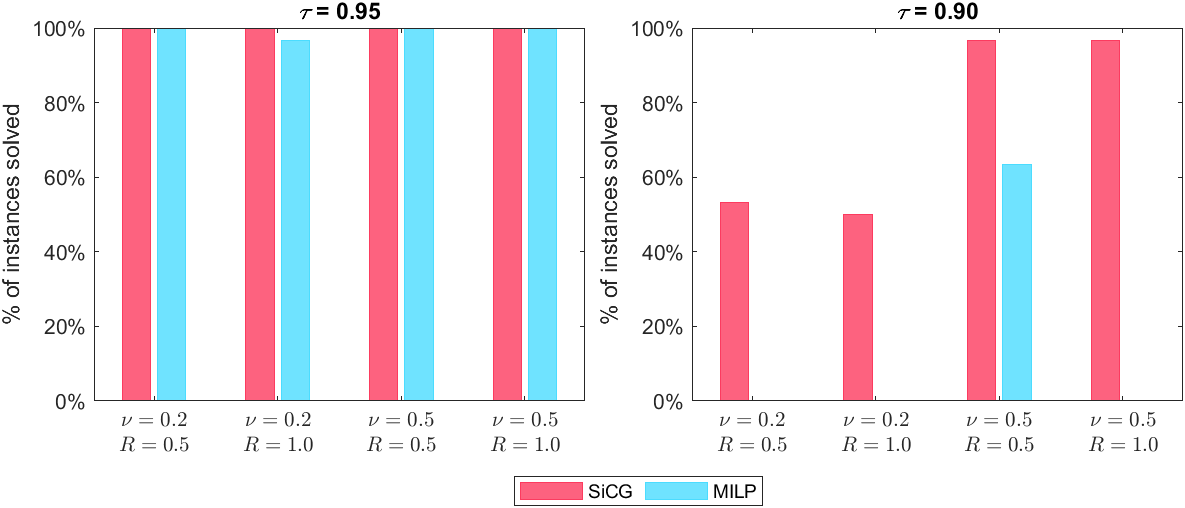}
\caption{Percentage of instances (over $30$ replications) solved within $1$ hour time limit using SiCG and MILP for $\tau\in\{0.95,0.90\}$ when $n=6$ and $N=1000$.} \label{fig:num_solved_n6}
\end{figure}

First, in Figure \ref{fig:num_solved_n6}, we show the percentage of instances (over 30 replications) that could be solved within a $1$-hour time limit when $N=1,000$, where we note that all the instances could be solved for a smaller $N$. When $\tau=0.95$, both SiCG and MILP can solve similar numbers of instances. However, when $\tau=0.90$, the number of instances that could be solved by SiCG is significantly greater than that by MILP.  Except under the setting of $\nu=0.5$ and $R=0.5$, MILP could not solve any instances within $1$ hour. In addition, we note that the relative MIP gap reported by the solver when terminated after $1$ hour by MILP is huge. For example, more than half of the terminating MIP gaps are greater than $20\%$ and some of them are even greater than $50\%$. These results illustrate that our proposed SiCG could efficiently solve instances, especially with a smaller $\tau$ (e.g., when $\tau=0.90$).

In Tables \ref{table:sol_time_n6_tau95}--\ref{table:sol_time_n6_tau90}, we present the lower quantile, median, and upper quantile of the solution times (over $30$ generated instances) under different numbers of scenarios and settings when $\tau=0.95$ and $\tau=0.90$, respectively. Because most of the instances when $N=1,000$ could not be solved by MILP, we did not present the solution times in the tables. It is observable that the solution time increases with the number of scenarios $N$. Moreover, instances with scenarios generated by a smaller variance (i.e., $\nu=0.2$) appear to be more challenging. This may be due to the difficulty in determining which scenarios are tail events when the scenarios are closer to each other. Despite the challenges of the non-convex ASP problem with a quantile objective, our proposed SiCG can solve most instances within $70$s and $1,600$s under $\tau=0.95$ and $\tau=0.90$, respectively. In addition, when $N=500$, the solution time of our proposed SiCG is significantly shorter than that of MILP. These results further demonstrate that our proposed SiCG is more computationally efficient than MILP.
\begin{table}[t]\centering \small
\ra{1.0}  
\caption{Solution time (in seconds) when $n=6$ and $\tau=0.95$ under different number of scenarios $N\in\{200,500\}$. Q1, Q2, and Q3, are the lower quantile, median, and upper quantile of solution times over $30$ generated instances.} \label{table:sol_time_n6_tau95}
\begin{tabular}{ll||rrr|rrr} \Xhline{1.0pt}
\multicolumn{2}{c||}{$N = 200$} & \multicolumn{3}{c|}{SiCG} & \multicolumn{3}{c}{MILP} \\
\multicolumn{2}{c||}{}          & Q1      & Q2     & Q3     & Q1     & Q2     & Q3     \\ \hline
$\nu = 0.2$     & $R = 0.5$    & 1.50    & 1.71   & 2.32   & 1.35   & 1.53   & 1.77   \\
$\nu = 0.2$     & $R = 1.0$    & 1.86    & 2.19   & 3.09   & 1.61   & 1.95   & 2.31   \\
$\nu = 0.5$     & $R = 0.5$    & 0.84    & 1.01   & 1.36   & 1.01   & 1.13   & 1.29   \\
$\nu = 0.5$     & $R = 1.0$    & 1.54    & 1.67   & 2.25   & 1.23   & 1.43   & 1.72   \\ \Xhline{1.0pt}
\multicolumn{2}{c||}{$N = 500$} & \multicolumn{3}{c|}{SiCG} & \multicolumn{3}{c}{MILP} \\
\multicolumn{2}{c||}{}          & Q1      & Q2     & Q3     & Q1     & Q2     & Q3     \\ \hline
$\nu = 0.2$     & $R = 0.5$    & 10.59   & 12.70  & 20.10  & 28.06  & 44.32  & 52.29  \\
$\nu = 0.2$     & $R = 1.0$    & 10.66   & 15.97  & 21.02  & 48.18  & 56.01  & 66.97  \\
$\nu = 0.5$     & $R = 0.5$    & 5.98    & 8.81   & 12.06  & 7.91   & 10.44  & 13.01  \\
$\nu = 0.5$     & $R = 1.0$    & 6.91    & 9.52   & 13.44  & 15.17  & 23.82  & 40.20  \\
\Xhline{1.0pt}
\end{tabular}
\end{table}
\begin{table}[t]\centering \small
\ra{1.0}  
\caption{Solution time (in seconds) when $n=6$ and $\tau=0.90$ under different number of scenarios $N\in\{200,500\}$. Q1, Q2, and Q3, are the lower quantile, median, and upper quantile of solution times over $30$ generated instances.} \label{table:sol_time_n6_tau90}
\begin{tabular}{ll||rrr|rrr} \Xhline{1.0pt}
\multicolumn{2}{c||}{$N = 200$} & \multicolumn{3}{c|}{SiCG} & \multicolumn{3}{c}{MILP}  \\
\multicolumn{2}{c||}{}          & Q1      & Q2     & Q3     & Q1     & Q2     & Q3      \\ \hline
$\nu = 0.2$     & $R = 0.5$    & 2.59    & 4.02   & 6.76   & 3.31   & 4.83   & 6.16    \\
$\nu = 0.2$     & $R = 1.0$    & 2.77    & 5.74   & 9.65   & 4.87   & 5.88   & 8.83    \\
$\nu = 0.5$     & $R = 0.5$    & 1.95    & 2.80   & 4.35   & 2.50   & 3.24   & 4.27    \\
$\nu = 0.5$     & $R = 1.0$    & 4.08    & 4.95   & 6.27   & 3.20   & 3.93   & 5.19    \\ \Xhline{1.0pt}
\multicolumn{2}{c||}{$N = 500$} & \multicolumn{3}{c|}{SiCG} & \multicolumn{3}{c}{MILP}  \\
\multicolumn{2}{c||}{}          & Q1      & Q2     & Q3     & Q1     & Q2     & Q3      \\ \hline
$\nu = 0.2$     & $R = 0.5$    & 71.87   & 143.93 & 217.13 & 324.17 & 460.69 & 730.16  \\
$\nu = 0.2$     & $R = 1.0$    & 122.01  & 167.06 & 266.19 & 660.85 & 857.32 & 1565.23 \\
$\nu = 0.5$     & $R = 0.5$    & 35.26   & 43.02  & 55.88  & 67.87  & 81.59  & 238.84  \\
$\nu = 0.5$     & $R = 1.0$    & 33.16   & 71.62  & 128.96 & 100.79 & 341.11 & 492.07  \\
\Xhline{1.0pt}
\end{tabular}
\end{table}

\subsection{Optimal Schedules} \label{subsec:opt_schedule}

In this section, we investigate the operational performances of our proposed methodology. We generate $N = 10,000$ pairs of historical data, and each pair contains the appointment characteristic and service duration. As a result, because we have $n$ appointments, the total number of potential scenarios amounts to $N^n$. Such a large number of scenarios cannot be used directly in optimization models from the computational perspective. Thus, we adopt the following sub-sampling technique to select $N'=1,000$ scenarios as the input to the optimization model. In \ref{appdx:expt_subsample}, we examine the effect of the random sub-sampling and observe that it does not significantly impact the operational performance.
\begin{itemize}[nosep,leftmargin=5mm]
    \item \textit{Sample Average Approximation (SAA)}. We obtain a random sub-sample of size $N'$. Then, we solve the ASP with a quantile objective using the empirical distribution of this sub-sample.
    \item \textit{Contextual Stochastic Optimization (CSO)}. In our proposed methodology, we employ the kernel $k(z;h)=\one(|z|\leq h)$ with bandwidth $h$ to each individual appointment and obtain a subset of historical data for each appointment, say, of size $N_i$. That is, we discard historical data with characteristic that is not close to the one we are planning for. Essentially, we are applying the kernel $k(\zb)=\one(\norm{\zb}_\infty \leq h)$ to the set of $N^n$ scenarios. Then, we obtain a random sub-sample of size $N'$ within the $\prod_{i=1}^n N_i$ scenarios for the optimization model. 
\end{itemize}

Note that because $N^n$ scenarios are involved, using the $k$NN regression approach is computationally difficult since it requires evaluating distances between the given predictor and the $N^n$ scenario. Thus, we do not consider this approach in our experiments. We generate $20$ sets of historical data and solve for both SAA and CSO models. We solve both models using SiCG to a relative gap of $5\%$ or terminate with a time limit $1,800$s. For the bandwidth $h$ in CSO, we pick $h=1.0$ since the solutions and performances using other tested bandwidths are similar (see \ref{appdx:expt_bandwidth}). As a benchmark, we also generate $20$ sets of scenarios generated from the true conditional distribution as described in Section \ref{subsec:expt_setting} and solve the optimization model. This benchmark allows us to investigate the effect of contextual information that is incorporated in our proposed methodology. In the experiments, we choose $\tau=0.95$ that corresponds to a risk-averse decision-maker who is more concerned about huge costs. 

\begin{figure}[t] 
\centering 
\includegraphics[scale=0.62]{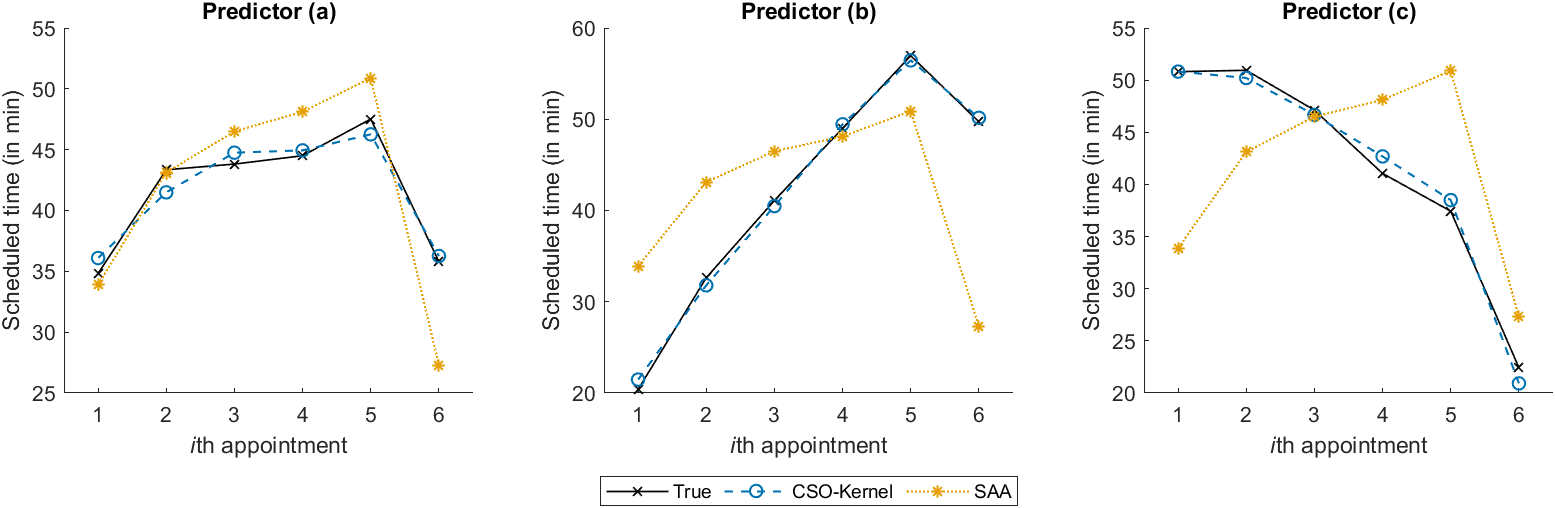}
\caption{Optimal appointment schedules from three different models when $\nu=0.2$ and $R=0.5$ under a quantile objective.} \label{fig:opt_schedule_1_quantile}
\end{figure}

First, we investigate the appointment time structure. For presentation brevity, we focus on the setting $\nu=0.2$ and $R=0.5$. Similar observations and discussions for the remaining settings can also be found in \ref{appdx:expt_opt_schedule}. Figure \ref{fig:opt_schedule_1_quantile} shows the optimal time assigned to each appointment, based on the average of the $20$ optimal solutions from $20$ generated sets of historical data, for three different predictors (see Section \ref{subsec:expt_setting}). Note that SAA is independent of the predictor information. From Figure \ref{fig:opt_schedule_1_quantile}, we observe that the optimal schedule from CSO is very close to the one obtained when simulating from the true distribution (which is unknown in practice). Moreover, CSO could capture the underlying conditional distribution as shown in the optimal schedules of predictors (b) and (c), where the mean service durations of the appointments vary. In contrast, SAA is not able to capture such an important information, and would eventually lead to a poorer operational performance.

\begin{figure}[t] 
\centering 
\includegraphics[scale=0.62]{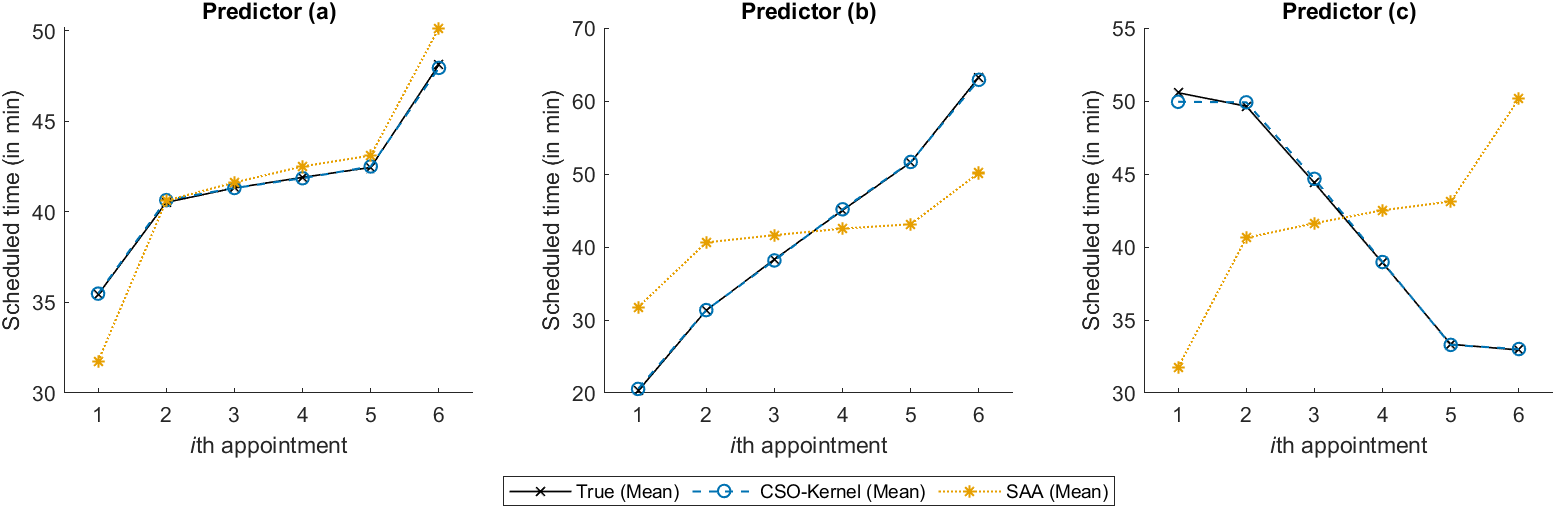}
\caption{Optimal appointment schedules from three different models when $\nu=0.2$ and $R=0.5$ under an expectation objective.} \label{fig:opt_schedule_1_mean}
\end{figure}

Next, we also study the difference in the optimal schedules using a quantile objective and an expectation objective. Figure \ref{fig:opt_schedule_1_mean} shows the optimal solutions when using a qunatile objective and an expectation objective. We observe that the optimal schedule using an expectation objective allocates less time to each appointment, except the last one (because we require that the sum of scheduled times to be always $T$). For instance, under predictor (a) where each appointment has a mean of $40$ minutes, each appointment is allocated with about $41$ minutes when using an expectation objective, as opposed to around $45$ minutes when using a quantile objective. This illustrates that using a quantile objective could lead to a more risk-averse schedule, in the sense that more time is assigned to each appointment to protect against unexpected long service durations, which may lead to huge waiting costs and thus, operational costs.

\subsection{Out-of-Sample Performance} \label{subsec:OS_performance}

To investigate the out-of-sample performance, for each given optimal solution (see Section \ref{subsec:opt_schedule}), we generate a new set of scenarios of size $N''=10,000$. Then, we implement the obtained optimal solutions on this new set of $N''$ scenarios, where we can compute the out-of-sample operational costs. 

\begin{figure}[t] 
\centering 
\includegraphics[scale=0.65]{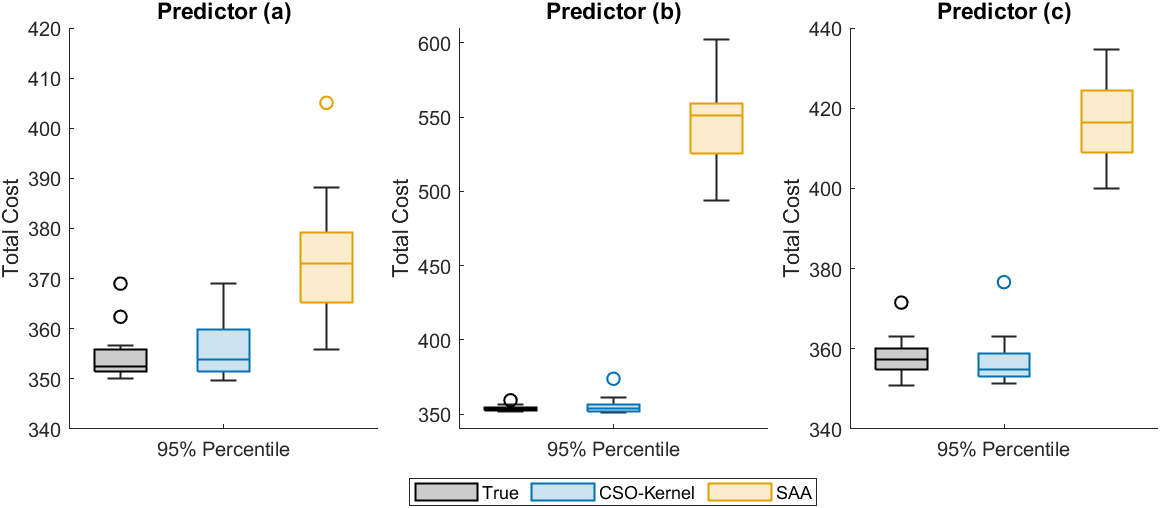}
\caption{Out-of-sample 95\% percentiles of the total cost over $20$ optimal schedules under a quantile objective from three different models when $\nu=0.2$ and $R=0.5$.} \label{fig:OS_95Per_1}
\end{figure}

First, we study the performance when generating the out-of-sample scenarios using the same distribution for in-sample data. Figure \ref{fig:OS_95Per_1} plots the out-of-sample 95\% percentiles of the total cost (i.e., the weighted sum of idle, waiting, and overtime costs) when $\nu=0.2$ and $R=0.5$. Note that a risk-averse decision-maker would be more concerned about tail behaviors than the mean behaviors. Since the optimal schedule from CSO is close to the one obtained when simulating from the true distribution, they both have similar out-of-sample performance. We observe that SAA yields a higher total cost, especially under predictors (b) and (c) where the distributions of the appointments are heterogeneous as opposed to predictor (a). These results illustrate that CSO, which is able to capture useful information that describes the underlying distribution, could be more preferable than the classical SAA.


For presentation brevity, in Table \ref{table:OS_95Per}, we summarize the out-of-sample 95\% percentiles of the total cost for all tested combinations of $\nu$ and $R$ (see \ref{appdx:expt_OS} for figures). As pointed out earlier, CSO and the one under true distribution give similar performance. When comparing CSO with SAA, we observe that CSO performs better than SAA in general, especially when $\nu=0.2$ and under predictors (b) and (c). Recall that $\nu$ controls the service duration variance. When $\nu$ is relatively large, the effect of the predictor (which in this case, a change of mean) is weakened. As a result, the benefits of CSO become less pronounced when $\nu=0.5$. Nevertheless, results from Table~\ref{table:OS_95Per} demonstrate that CSO is still more preferable than SAA in general and could yield significant improvement in many cases.
\begin{table}[t]\centering \small
\ra{1.0}  
\caption{Out-of-sample 95\% percentiles of the total cost under a quantile objective from three different models. Q1, Q2, and Q3, are the lower quantile, median, and upper quantile over $20$ optimal solutions.} \label{table:OS_95Per}
\begin{tabular}{ll||rrr|rrr|rrr} \Xhline{1.0pt}
\multicolumn{2}{c||}{Predictor (a)} & \multicolumn{3}{c|}{True} & \multicolumn{3}{c|}{CSO-Kernel} & \multicolumn{3}{c}{SAA}  \\ 
\multicolumn{2}{l||}{}              & Q1      & Q2     & Q3     & Q1        & Q2       & Q3       & Q1     & Q2     & Q3     \\ \hline
$\nu = 0.2$       & $R = 0.5$      & 351.6   & 352.4  & 355.8  & 351.5     & 353.9    & 359.8    & 365.2  & 373.1  & 379.2  \\
$\nu = 0.2$       & $R = 1.0$      & 253.2   & 255.2  & 258.1  & 253.0     & 255.7    & 261.8    & 267.4  & 275.1  & 279.0  \\
$\nu = 0.5$       & $R = 0.5$      & 979.4   & 984.9  & 990.8  & 980.7     & 989.7    & 998.0    & 978.9  & 986.5  & 998.2  \\
$\nu = 0.5$       & $R = 1.0$      & 732.0   & 735.4  & 747.5  & 738.6     & 747.0    & 757.4    & 735.1  & 743.0  & 751.6  \\  \Xhline{1.0pt}
\multicolumn{2}{c||}{Predictor (b)} & \multicolumn{3}{c|}{True} & \multicolumn{3}{c|}{CSO-Kernel} & \multicolumn{3}{c}{SAA}  \\
\multicolumn{2}{l||}{}              & Q1      & Q2     & Q3     & Q1        & Q2       & Q3       & Q1     & Q2     & Q3     \\ \hline
$\nu = 0.2$       & $R = 0.5$      & 352.4   & 353.6  & 354.7  & 352.0     & 353.5    & 356.4    & 525.8  & 551.2  & 559.4  \\
$\nu = 0.2$       & $R = 1.0$      & 254.3   & 255.5  & 256.1  & 253.9     & 255.5    & 259.3    & 428.6  & 446.9  & 454.9  \\
$\nu = 0.5$       & $R = 0.5$      & 972.4   & 975.6  & 983.6  & 975.3     & 981.3    & 987.4    & 1044.7 & 1074.3 & 1093.2 \\
$\nu = 0.5$       & $R = 1.0$      & 729.3   & 732.4  & 739.3  & 730.8     & 735.6    & 739.0    & 801.9  & 830.9  & 846.7  \\  \Xhline{1.0pt}
\multicolumn{2}{c||}{Predictor (c)} & \multicolumn{3}{c|}{True } & \multicolumn{3}{c|}{CSO-Kernel} & \multicolumn{3}{c}{SAA}  \\
\multicolumn{2}{l||}{}              & Q1      & Q2     & Q3     & Q1        & Q2       & Q3       & Q1     & Q2     & Q3     \\ \hline
$\nu = 0.2$       & $R = 0.5$      & 354.8   & 357.3  & 360.1  & 353.0     & 354.9    & 358.9    & 409.0  & 416.3  & 424.5  \\
$\nu = 0.2$       & $R = 1.0$      & 256.2   & 260.2  & 262.4  & 255.4     & 256.9    & 260.3    & 313.4  & 320.8  & 324.7  \\
$\nu = 0.5$       & $R = 0.5$      & 977.1   & 982.4  & 991.6  & 980.3     & 986.7    & 995.1    & 1001.0 & 1019.9 & 1029.4 \\
$\nu = 0.5$       & $R = 1.0$      & 732.4   & 740.1  & 744.8  & 736.1     & 743.3    & 751.1    & 759.4  & 773.6  & 783.0  \\
\Xhline{1.0pt}
\end{tabular}
\end{table}

\begin{figure}[t] 
\centering 
\includegraphics[scale=0.65]{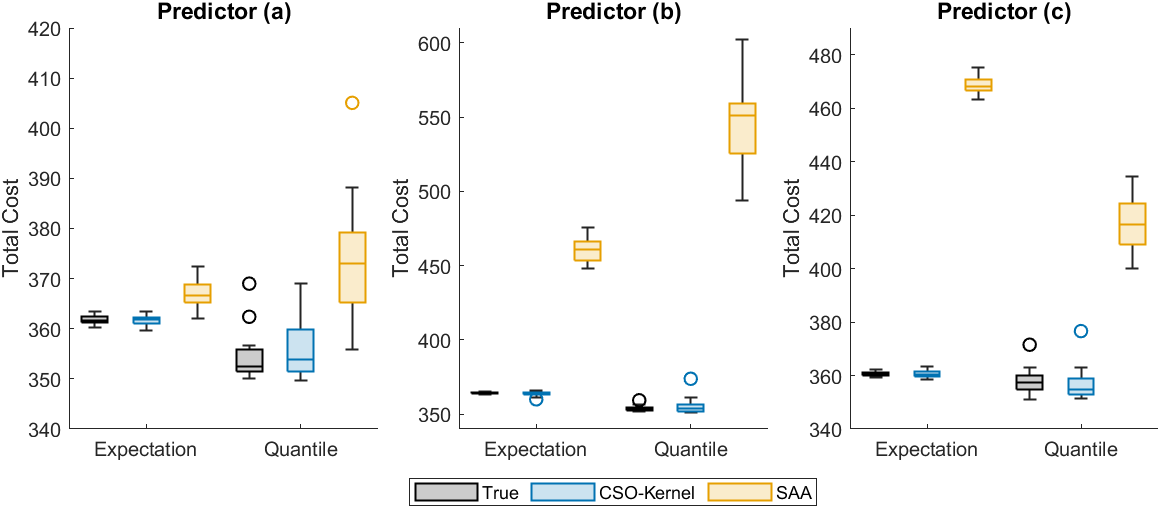}
\caption{Out-of-sample 95\% percentiles of the total cost over $20$ optimal schedules under an expectation or quantile objective from three different models when $\nu=0.2$ and $R=0.5$.} \label{fig:OS_with_mean_95Per_1}
\end{figure}

Next, we compare the results with optimal schedules obtained from using an expectation objective. Figure \ref{fig:OS_with_mean_95Per_1} shows the out-of-sample 95\% percentiles of the total cost with both expectation and quantile objectives when $\nu=0.2$ and $R=0.5$. Note that the variability of the performance is larger under a quantile objective because in such optimization problems, only scenarios corresponding to the tail events would be used in determining the optimal schedule. From Figure \ref{fig:OS_with_mean_95Per_1}, we  observe  that the out-of-sample 95\% percentiles using the quantile objective are generally smaller than those using the expectation objective. As we could expect from the optimal schedules  (see Section \ref{subsec:opt_schedule}), the waiting time using the quantile objective is significantly shorter than that using the expectation objective because more time is assigned to each appointment.  However,  as a result,  the overtime using the quantile objective is slightly longer because the last appointment is allocated with a shorter time. Thus, these results illustrate that using a quantile objective could offer effective protection against possible appointment delays. In \ref{appdx:expt_OS}, we present similar results for the remaining settings. We point out that when $R = 1.0$, i.e., the length of the working hour becomes longer, the respective performances using the two objectives become similar due to the increased flexibility in time scheduling.

Finally, we investigate the out-of-sample performance under mis-specified distributions. This may represent the situation where the actual service duration behaves differently from the historical data. We consider the following two different sets of distribution perturbations in generating out-of-sample scenarios: (I) increasing the service duration variance by a factor of $1.5$; (II) increasing the service duration mean by a factor of $1.2$. These two changes correspond to situations where the actual service duration is more volatile or longer than what  historical data would suggest. We compare the out-of-sample performance among our proposed methodology, SAA, and models with a mean objective.

Table \ref{table:OS_misspec_95Per_1} shows the mean and standard deviation of percentage change in the 95\% percentile of total cost when changing the implementation of the CSO solution based on a quantile objective to the SAA solution or the CSO solution based on an expectation objective under $\nu=0.2$ and $R=0.5$. A positive percentage change refers to a smaller 95\% percentile in the CSO solution based on a quantile objective. In \ref{appdx:expt_OS}, we present similar results for the remaining combinations of $(\nu,R)$. First, we observe that the advantage of using CSO over SAA becomes less significant when there is a change in distribution. This is because the predictors become less informative when the out-of-sample scenarios come from a perturbed distribution. Although CSO may be similar to SAA in performance in some cases (when the distribution is perturbed), we still observe a better performance of CSO under predictors (b) and (c).  Next, when compared CSO with the one using an expectation objective, our CSO based on a quantile objective performs well even when the distribution is perturbed. In particular, under Set II where the mean of appointment duration is increased, the improvement in the 95\% percentile becomes larger. This illustrates the advantage of using a quantile objective that yields a more risk-averse schedule than using an expectation objective even under mis-specified distributions. These results demonstrate the benefits of our proposed approach that integrates contextual information into a quantile objective.
\begin{table}[t]\centering \small
\ra{1.0}  
\caption{Mean and standard deviation of the percentage change in the 95\% percentile of total cost when changing the implementation of the CSO solution based on a quantile objective to the SAA solution or the CSO solution based on an expectation objective over $20$ optimal solutions under $\nu=0.2$ and $R=0.5$. (The `Original' column corresponds to the out-of-sample scenarios generated under the true distribution as for the in-sample scenarios).} \label{table:OS_misspec_95Per_1}
\begin{tabular}{l||rr|rr|rr} \Xhline{1.0pt}
Predictor (a)       & \multicolumn{2}{c|}{Original} & \multicolumn{2}{c|}{Set I} & \multicolumn{2}{c}{Set II} \\
                    & Mean          & Std          & Mean         & Std        & Mean         & Std         \\ \hline
CSO v.s. SAA        & 4.7\%         & 3.3\%        & -0.1\%       & 1.3\%      & -1.3\%       & 1.8\%       \\
CSO v.s. CSO (Mean) & 1.5\%         & 1.8\%        & 1.0\%        & 0.8\%      & 2.9\%        & 1.5\%       \\ \Xhline{1.0pt}
Predictor (b)       & \multicolumn{2}{c|}{Original} & \multicolumn{2}{c|}{Set I} & \multicolumn{2}{c}{Set II} \\
                    & Mean          & Std          & Mean         & Std        & Mean         & Std         \\ \hline
CSO v.s. SAA        & 53.0\%        & 7.8\%        & 22.5\%       & 4.2\%      & -0.5\%       & 3.1\%       \\
CSO v.s. CSO (Mean) & 2.4\%         & 1.6\%        & 1.5\%        & 0.9\%      & 2.9\%        & 1.0\%       \\ \Xhline{1.0pt}
Predictor (c)       & \multicolumn{2}{c|}{Original} & \multicolumn{2}{c|}{Set I} & \multicolumn{2}{c}{Set II} \\
                    & Mean          & Std          & Mean         & Std        & Mean         & Std         \\ \hline
CSO v.s. SAA        & 16.9\%        & 3.6\%        & 6.7\%        & 1.9\%      & 9.6\%        & 1.9\%       \\
CSO v.s. CSO (Mean) & 1.0\%         & 1.5\%        & 0.3\%        & 0.6\%      & 2.6\%        & 1.4\%       \\ 
\Xhline{1.0pt}
\end{tabular}
\end{table}

\section{Conclusion} \label{sec:conclusion}
In this paper, we propose a new risk-averse contextual stochastic optimization problem with a quantile objective for general two-stage stochastic programs. The use of quantile in the objective allows decision-makers to take their risk aversiveness into account. Moreover, the use of quantile objective that corresponds to the extreme losses is intuitive for decision-makers to adopt as compared with other risk measures in many application domains such as portfolio optimization and optimal treatment problems. Given a set of data consisting of the unknown random parameters' historical realizations and contextual information, we model the conditional quantile by replacing the conditional expectation of the checked loss function in its variational characterization with a generic estimator. Under two sets of mild regularity conditions, we derive the asymptotic almost-sure convergence and convergence in probability of the optimal value and solutions of the corresponding optimization problem to their true counterparts. Importantly, we establish the convergence in probability under a set of weaker assumptions, which enables the use of a larger class of conditional expectation estimators developed in the statistics literature.

While the use of quantile objective usually leads to non-convex problem which is challenging to solve, we propose a new stochastic inexact constraint generation (SiCG) method and prove its convergence. In our experiments on a single-server appointment scheduling problem, we demonstrate that our proposed SiCG has superior computational advantages over the standard mixed-integer linear program reformulations. Moreover, we also derive managerial insights of our proposed methodology. In particular, we illustrate the importance of incorporating potential contextual information into the optimization model. In addition, our results also demonstrate that using a quantile objective leads to a more risk-averse schedule, which may be preferred by risk-averse decision-makers. In general, our model provides an attractive alternative to existing models from a data-driven and risk-averse perspective.

This paper combines contextual information and risk preferences into a single optimization model. It would be interesting to extend our work by considering other risk measures, such as conditional value-at-risk and generic coherent risk measures. Moreover, as in many existing studies on contextual stochastic optimization, the set of contextual information is pre-specified.  It would be interesting to investigate how decision-makers could select a subset of meaningful predictors from a possibly large number of predictors in an optimization context. Finally, we would like to apply  our proposed methodology to different practical problems and investigate its potential advantages in other application domains.


\appendix
\newpage

\section{Mathematical Proofs} \label{appdx:math_proof}

\subsection{Proof of Theorem \ref{thm:Lip_quantile}} \label{appdx:thm_Lip_quantile}

Without loss of generality, assume that $\calX$ is connected. Otherwise, one could apply the argument on a connected set $\calC$ such that $\calX\subseteq\calC$. First, note that since $f(\xb,\xib)$ has an optimal solution, we can write $f(\xb,\xib)=\max_{\pib\in\Pit}\big\{ (\hb-\Tb\xb-\Cb\xib)^\tp\pib \big\}$ for some finitely countable set $\Pit$ (i.e., extreme point of the dual feasible set $\{\pib \geq \zero \mid \Wb^\tp \pib\leq \qb\}$). That is, $f(\xb,\xib)$ is a max-affine function for each $\xib\in\Xi$. Also, recall from \eqref{eqn:cond_quantile_est} that we have $\Qh_\tau^N(\xb;\zb)=\inf\{t \mid \Fh^n(t\mid \zb) \geq \tau\}$. Therefore, for a given $\xb\in\calX$, we have $\Qh_\tau^N(\xb;\zb)=f(\xb,\xib^i)$ for some $i\in\{1,\dots,N\}$.

Now, consider a given $\xbt\in\calX$ and let $\itt$ be the smallest index such that $\Qh^N_\tau(\xbt;\zb)=f(\xbt,\xib^{\itt})$. Denote the interior and closure of a set $A$ by $\intr(A)$ and $\cl(A)$, respectively. Construct the set
\begin{equation*}
    \calN(\xbt)=\cl\left\{\xb\in\calX \,\left|\,
    \begin{array}{l}
        \itt\in\argmin\big\{ i\in\{1,\dots,N\}\mid \Qh^N_\tau(\xb;\zb)=f(\xb,\xib^i)\big\},   \\
        \exists \text{ continuous } \gamma:[0,1]\rightarrow\calX \text{ s.t. } \gamma(0)=\xb,\, \gamma(1)=\xbt, \text{ and} \\
        \itt\in\argmin\big\{i\in\{1,\dots,N\} \,\big|\, \Qh_\tau^N(\gamma(t);\zb)=f(\gamma(t),\xib^i) \big\},\, \forall t\in[0,1]
    \end{array}
    \right. \right\}.
\end{equation*}
Note that the first constraint ensures that $\itt$ is the smallest index such that the estimated quantile equals the second-stage cost at scenario $\xib^i$. The second constraint says that $\xb\in\calN(\xbt)$ is always path-connected to $\xbt$. In addition, points along the path must satisfy the same condition as in the first constraint. By construction, the set $\calN(\xbt)$ is closed. We want to show that the boundaries of $\calN(\xbt)$, denoted as $\partial\calN(\xbt)$, can be characterized by the boundary of $\calX$ and hyperplanes. (Note that $\calN(\xbt)$ is not necessarily convex as shown in Example \ref{eg:Lip_quantile}). 

Note that by definition, a point $\xb\in\partial N(\xbt)$ if $\xb\in N(\xbt)$ but $\xb\not\in\intr(\calN(\xbt))$. Also, note that $\xb\not\in\intr(\calN(\xbt))$ is equivalent to one of the following: there does not exist $r>0$ such that for all $\xbhat\in\calB_r(\xb)$ (i.e., a ball centered at $\xb$ with radius $r$) such that (a) $\calB_r(\xb)\subseteq \calX$, or (b) $\itt$ is not the smallest index such that $\Qh^N_\tau(\xbhat;\zb)=f(\xbhat,\xib^i)$. For case (a), it is trivial that the boundaries are characterized by that of $\calX$. For case (b), we first define the set
$$\calI(\xb)=\big\{i\in\{1,\dots,N\}\mid \exists r>0 \text{ s.t. } f(\xbhat,\xib^{i})=f(\xbhat,\xi^{\itt}),\, \forall \xbhat\in \calB_r(\xb)\big\}.$$
That is, for all $i\in\calI(\xb)$ (including $\itt$), the functions $f(\xb,\xib^i)$ overlap with each other in a neighborhood of $\xb$. We claim that there exists $j\not\in\calI(\xb)$ such that $f(\xb,\xib^{\itt})=f(\xb,\xib^j)$. Suppose on the contrary, that for all $j\not\in\calI(\xb)$, we have $f(\xb,\xib^{\itt})\ne f(\xb,\xib^j)$.
\begin{itemize}[topsep=1pt,itemsep=-1pt]
\item If $j\not\in\calI(\xb)$, then trivially, by our assumption $f(\xb,\xib^{\itt})\ne f(\xb,\xib^j)$.
\item If $j\in\calI(\xb)$, then by definition of $\calI(\xb)$, $f(\xb,\xib^j)$ and $f(\xb,\xib^{\itt})$ are identical in a neighborhood of $\xb$.
\end{itemize}
These two cases imply that there exist a neighborhood of $\xb$ such that $\itt$ is still the smallest index such that $\Qh^N_\tau(\xb;\zb)=f(\xb,\xib^i)$, which contradicts the condition (b). Therefore, we can conclude that the boundary is a result of the intersection between two non-overlapping functions $f(\xb,\xib^{\itt})$ and $f(\xb,\xib^j)$ (i.e., not identical in a neighborhood of $\xb$). Since $f(\xb,\xib)$ is a max-affine function, intersections between the two functions $f(\cdot,\xib^i)$ and $f(\cdot,\xib^j)$ can be characterized by hyperplanes (i.e., affine equalities). This shows the boundaries of $\calN(\xb)$ in case (b) can be characterized by hyperplanes.

Next, we claim that there is at most finitely many distinct $\calN(\xb)$ on $\calX$. Suppose, on the contrary, that the number of distinct $\calN(\xb)$ is not finite, i.e., we have a sequence of distinct sets $\{\calN_k\}_{k\in\N}$. Pick one $\xb_k$ in each distinct set $\calN_k$. Since $\xb_k\in \calN_k$, we let $i_k$ be the smallest index such that $\Qh_\tau^N(\xb_k;\zb) =f(\xb_k,\xib^{i_k})$. Since $i_k\in\{1,\dots,N\}$ for all $k\in\N$, there exists $i\in\{1,\dots,N\}$ such that $i$ appears infinitely often in $\{i_k\}_{k\in\N}$. This implies that $f(\xb,\xib^i)$ intersects with other functions $f(\xb,\xib^j)$ for $j\ne i$ infinitely many times. However, $f(\xb,\xib^j)$ is max-affine with a finite number of pieces, say $p$. Then, for each piece of affine function in $f(\xb,\xib^i)$, it intersects with $f(\xb,\xib^j)$ for $j\ne i$ at most $p(N-1)$ times. Therefore, since $f(\xb,\xib^i)$ has only $p$ pieces of affine functions, there is at most $p^2(N-1)$ intersections. This results in a contradiction that $f(\xb,\xib^i)$ intersects with other functions $f(\xb,\xib^j)$ for $j\ne i$ infinitely many times.

Finally, let the sequence of distinctive $\calN(\xb)$ be $\{\calN_k\}_{k=1}^K$. Note that $f(\xb,\xib^i)$ is Lipschitz continuous since it is a max-affine function with a finite number of pieces. Therefore, $\Qh_\tau^N(\xb;\zb)$ is Lipschitz continuous on the entire domain $\calX$ (probably) except on the boundary of each set $\calN_k$. Define $\Theta=\bigcup_{k=1}^K \partial\calN_k$, which is the union of finitely many hyperplanes. By \cite{Leobacher_Steinicke:2022}, piecewise Lipschitz continuous function of this type (i.e., piecewise Lipschitz with polyhedral intersections) is also Lipschitz continuous on the entire domain. As a result, $\Qh_\tau^N(\xb;\zb)$ is Lipschitz continuous on $\calX$. Since $f(\xb,\xib^i)$ is Lipschitz continuous with Lipschitz constant $L=\max_{\pib\in\Pit}\norms{\Tb^\tp\pib}<\infty$ for all $i\in\{1,\dots,N\}$,   $\Qh_\tau^N(\xb;\zb)$ is also Lipschitz continuous with Lipschitz constant $L$.

\subsection{Proof of Theorem~\ref{thm:convergence_as_est_quantile}} \label{appdx:thm_convergence_est_quantile_pf}

Recall that we would like to show that for any $\xb\in\calX$, the estimated quantile $\Qh_\tau^N(\xb;\zb)$, the minimizer of \eqref{eqn:CQM_sample}, converges to the true conditional quantile $Q_\tau(\xb;\zb)$, the minimizer of \eqref{eqn:cond_quantile_min_characterization}, which is unique by Assumption \ref{assumption:regularity_cond_1}. 

\begin{itemize}[leftmargin=5mm]
    \item Note that for every $u\in\R$,
    \begin{align*}
    \E_\Prob\big[\rho_\tau\big(f(\xb,\xib)-u\big)\big] 
    &=\tau \E_\Prob\Big\{ \big( f(\xb,\xib)-u\big) \one\big(u < f(\xb,\xib)\big) \Big\} + (1-\tau) \E_\Prob\Big\{ \big( u-f(\xb,\xib)\big) \one\big(u \geq f(\xb,\xib)\big) \Big\} \\
    &\leq \tau \E_\Prob\big|f(\xb,\xib)-u\big| + (1-\tau) \E_\Prob\big|f(\xb,\xib)-u\big| \\
    & <\infty
    \end{align*}
    since $f(\xb,\xib)$ is integrable by Assumption \ref{assumption:regularity_cond_1}. Moreover, the optimal solution to \eqref{eqn:cond_quantile_min_characterization}, i.e., $Q_\tau(\xb;\zb)$ always exists.

    \item For any $y\in\R$, the function $\rho_\tau(y-u)$ is equicontinuous in $u$. That is, for any $u\in\R$ and $\epsilon>0$, there exists $\delta>0$ such that $|\rho_\tau(y-u)-\rho_\tau(y-u')| \leq \varepsilon$ for all $y\in\R$ and $u'$ with $|u-u'|\leq\delta$. Indeed, for any $y\in\R$, by picking $\delta=\varepsilon$, we have
    $$|\rho_\tau(y-u)-\rho_\tau(y-u')| \leq |u-u'| \leq \varepsilon.$$

    \item Since $\rho_\tau(y-u) \geq 0$ for any $y\in\R$ and $u\in\R$, we have $\liminf_{|u|\rightarrow\infty} \inf_{y\in\R} \rho_\tau(y-u) \geq 0$. Also, by assumption~(a) in Theorem~\ref{thm:convergence_as_est_quantile}, there exists a compact set $B_{\zb}\subseteq\Xi$ such that $\Prob(\xib\in B_{\zb}\mid \Zb=\zb)>0$. Therefore, since $f(\xb,\cdot)$ is continuous, there exists a compact set $P_{\zb}\subset\R$ such that $\Prob(f(\xb,\xib)\in P_{\zb}\mid\Zb=\zb)>0$. Moreover, by the coerciveness of $\rho_\tau(\cdot)$, we have $\lim_{|u|\rightarrow\infty} \inf_{y\in P_{\zb}} \rho_\tau(y-u)=\infty$. 
\end{itemize}

Together with assumptions~(b) and (c) in Theorem~\ref{thm:convergence_as_est_quantile}, it follows from Lemma~EC.6 of \cite{Bertsimas_Kallus:2020} that $\Qh^N_{\tau}(\xb;\zb)\as  Q_\tau(\xb;\zb)$ as $N\rightarrow\infty$ for any $\xb\in\calX$.

\subsection{Proof of Lemma~\ref{lem:unif_convergence_as_obj}} \label{appdx:lem_unif_convergence_obj}

Consider a sample path of data such that $\Qh_\tau^N(\xb;\zb)\rightarrow Q_\tau(\xb;\zb)$ for any $\xb\in\calX$. Let $\{\xb^N\}\subseteq \calX$ be a sequence such that $\xb^N\rightarrow \xbbar\in\calX$ as $N\rightarrow\infty$. Note that
\begin{align}
\Big|\Qh_\tau^N(\xb^N;\zb) - Q_\tau(\xbbar;\zb) \Big| &\leq \Big|\Qh_\tau^N(\xb^N;\zb) - \Qh_\tau^N(\xbbar;\zb) \Big| + \Big| \Qh_\tau^N(\xbbar;\zb) - Q_\tau(\xbbar;\zb) \Big| \nonumber \\
&\leq L_2\norm{\xb^N-\xbbar}^\beta+ \Big| \Qh_\tau^N(\xbbar;\zb) - Q_\tau(\xbbar;\zb) \Big|, \label{eqn_pf:lem:unif_convergence_as_obj_1}
\end{align}
where the last inequality follows from Assumption~\ref{assumption:Holder_cont}(b). Since $\xb^N\rightarrow\xbbar$, the first term in \eqref{eqn_pf:lem:unif_convergence_as_obj_1} converges to zero. Also, since $\Qh_\tau^N(\xb;\zb)\rightarrow Q_\tau(\xb;\zb)$ for any $\xb\in\calX$, the second term in \eqref{eqn_pf:lem:unif_convergence_as_obj_1} also converges to zero. Therefore, we have $\Qh_\tau^N(\xb^N;\zb) \rightarrow Q_\tau(\xbbar;\zb)$ as $N\rightarrow\infty$.

Now, we prove the desired uniform convergence via a contradiction argument. Suppose, on the contrary, that $\sup_{\xb\in \calX} \big|\Qh_\tau^N(\xb;\zb) - Q_\tau(\xb;\zb) \big| \not\rightarrow 0$. Then, there exist $\varepsilon>0$ and a sequence $\{\xb^N\}\subseteq\calX$ such that
\begin{equation} \label{eqn_pf:lem:unif_convergence_as_obj_2}
\Big|\Qh_\tau^N(\xb^N;\zb) - Q_\tau(\xb^N;\zb) \Big| \geq \varepsilon
\end{equation}
infinitely often. Pick a sequence such that \eqref{eqn_pf:lem:unif_convergence_as_obj_2} always hold, and then a convergent subsequence $\{\xb^{N_k}\}_{k\in\N}$ such that $\xb^{N_k}\rightarrow\xbbar\in\calX$. Note that such a convergent subsequence can be constructed since $\calX$ is compact by Assumption~\ref{assumption:compact_X}. Thus, we have
$$\varepsilon \leq \Big| \Qh_\tau^{N_k}(\xb^{N_k};\zb) - Q_\tau(\xb^{N_k};\zb) \Big| \leq  \Big| \Qh_\tau^{N_k}(\xb^{N_k};\zb) - Q_\tau(\xbbar;\zb) \Big| + \Big| Q_\tau(\xbbar;\zb) - Q_\tau(\xb^{N_k};\zb) \Big|.$$
From the argument in the first part of the proof, since $\Qh_\tau^N(\xb^N;\zb) \rightarrow Q_\tau(\xbbar;\zb)$ as $N\rightarrow\infty$, there exists $k_1\in\N$ such that $\big| \Qh_\tau^{N_k}(\xb^{N_k};\zb) - Q_\tau(\xbbar;\zb) \big| \leq \varepsilon/4$ for all $k\geq k_1$. Also, by Assumption~\ref{assumption:Holder_cont}(a), there exists $k_2\in\N$ such that $\big| Q_\tau(\xbbar;\zb) - Q_\tau(\xb^{N_k};\zb) \big|\leq \varepsilon/4$ for all $k \geq k_2$. Therefore, for all $k\geq \max\{k_1,\,k_2\}$, we have $\varepsilon \leq \big| \Qh_\tau^{N_k}(\xb^{N_k};\zb) - Q_\tau(\xb^{N_k};\zb) \big| \leq  \varepsilon/2$, which gives the desired contradiction. Finally, by Theorem~\ref{thm:convergence_as_est_quantile}, we have $\Qh_\tau^N(\xb;\zb)\as Q_\tau(\xb;\zb)$ for any $\xb\in\calX$. Therefore, we conclude that $\sup_{\xb\in \calX} \big|\Qh_\tau^N(\xb;\zb) - Q_\tau(\xb;\zb) \big| \as 0$.

\subsection{Proof of Theorem \ref{thm:convergence_as_opt_val_sol}} \label{appdx:thm:convergence_as_opt_val_sol}


For part~(i), note that it follows from Lemma~\ref{lem:unif_convergence_as_obj} that
\begin{equation} \label{eqn_pf:thm:convergence_as_opt_val_sol_1}
    \delta^N := \Big| \min_{\xb\in\calX} \Qh_\tau^N(\xb;\zb) - \min_{\xb\in\calX} Q_\tau(\xb;\zb) \Big| \leq \sup_{\xb\in\calX} \Big|\Qh_\tau^N(\xb;\zb) - Q_\tau(\xb;\zb)\Big| \as 0.
\end{equation}
Next, we prove part~(ii). Note that $\sup_{\xb^N\in\calXhat^N(\zb)}  \Qh_\tau^N(\xb^N;\zb) =  \min_{\xb\in\calX}  \Qh_\tau^N(\xb;\zb) $ since $\calXhat^N(\zb)$ is the set of optimal solutions to problem~\eqref{eqn:CQM_sample}. Therefore, we have
\begin{align*}
&\quad\,\,\Bigg| \sup_{\xb^N\in\calXhat^N(\zb)} Q_\tau(\xb^N;\zb) - v^\star(\zb) \Bigg| \\
&=\Bigg| \sup_{\xb^N\in\calXhat^N(\zb)} Q_\tau(\xb^N;\zb) - \min_{\xb\in\calX}  Q_\tau(\xb;\zb) \Bigg| \\
& \leq \Bigg| \sup_{\xb^N\in\calXhat^N(\zb)} Q_\tau(\xb^N;\zb) -  \sup_{\xb^N\in\calXhat^N(\zb)} \Qh_\tau^N(\xb^N;\zb)\Bigg| + \bigg| \min_{\xb\in\calX} \Qh_\tau^N(\xb;\zb) - \min_{\xb\in\calX} Q_\tau(\xb;\zb) \bigg| \\
&\leq \sup_{\xb\in\calX} \Big|\Qh_\tau^N(\xb;\zb) - Q_\tau(\xb;\zb)\Big| + \delta^N \as 0,
\end{align*}
where the convergence follows from \eqref{eqn_pf:thm:convergence_as_opt_val_sol_1}. This proves part~(ii). Finally, we prove part~(iii). Suppose, on the contrary, that $D(\calXhat^N(\zb),\calX^\star(\zb))\not\as 0$. Consider a sample path of data such that $D(\calXhat^N(\zb),\calX^\star(\zb))\not\rightarrow 0$. Then, there exists a sequence $\{\xb^{N_k}\}_{k\in\N}\subseteq\calX$ such that $\xb^{N_k}\in\calXhat^{N_k}(\zb)$, $\xb^{N_k}\rightarrow\xbbar\in\calX$, and $\inf_{\xb\in\calX^\star(\zb)} \norms{\xb^{N_k}-\xb} = \eta > 0$. Note that $\xbbar\not\in\calX^\star(\zb)$ since
$$\inf_{\xb\in\calX^\star(\zb)} \norms{\xbbar - \xb} \geq \inf_{\xb\in\calX^\star(\zb)} \norms{\xb^{N_k}-\xb} - \norms{\xb^{N_k} - \xbbar} \rightarrow \eta > 0.$$
Define $\varepsilon := Q_\tau(\xbbar;\zb) - v^\star(\zb) > 0$. By Assumption~\ref{assumption:Holder_cont}(a), the quantile function $Q_\tau(\xb;\zb)$ is continuous in $\xb$. Therefore, there exists $k_0\in\N$ such that we have $\big| Q_\tau(\xb^{N_k};\zb) -  Q_\tau(\xbbar;\zb)\big| \leq \varepsilon/2$ for all $k\geq k_0$. Then, for $k\geq k_0$, we have
\begin{subequations}
\begin{align}
C_k &:=\sup_{\xb\in\calX} \Big|\Qh_\tau^{N_k}(\xb;\zb) - Q_\tau(\xb;\zb)\Big| + \delta^{N_k} \nonumber  \\
&\geq \Big|  \Qh_\tau^{N_k}(\xb^{N_k};\zb) - Q_\tau(\xb^{N_k};\zb)\Big| + \Big|  \Qh_\tau^{N_k}(\xb^{N_k};\zb) - \min_{x\in\calX} Q_\tau(\xb^{N_k};\zb) \Big|  \label{eqn_pf:thm:convergence_as_opt_val_sol_2a} \\
&\geq \Big| Q_\tau(\xb^{N_k};\zb) -  \min_{x\in\calX} Q_\tau(\xb^{N_k};\zb) \Big| \label{eqn_pf:thm:convergence_as_opt_val_sol_2b} \\
&\geq  Q_\tau(\xb^{N_k};\zb) -  \min_{x\in\calX} Q_\tau(\xb^{N_k};\zb) \label{eqn_pf:thm:convergence_as_opt_val_sol_2c} \\
&\geq Q_\tau(\xbbar;\zb) -  \min_{x\in\calX} Q_\tau(\xb^{N_k};\zb) - \frac{\varepsilon}{2} \label{eqn_pf:thm:convergence_as_opt_val_sol_2d}\\
&= \frac{\varepsilon}{2}. \label{eqn_pf:thm:convergence_as_opt_val_sol_2e}
\end{align}
\end{subequations}
Here, \eqref{eqn_pf:thm:convergence_as_opt_val_sol_2a} follows from $\xb^{N_k}\in\calX$ and the definition of $\delta^{N_k}$ in \eqref{eqn_pf:thm:convergence_as_opt_val_sol_1}; \eqref{eqn_pf:thm:convergence_as_opt_val_sol_2b} follows from the triangular inequality; \eqref{eqn_pf:thm:convergence_as_opt_val_sol_2d} follows from  $\big| Q_\tau(\xb^{N_k};\zb) -  Q_\tau(\xbbar;\zb)\big| \leq \varepsilon/2$ for all $k\geq k_0$; \eqref{eqn_pf:thm:convergence_as_opt_val_sol_2e} follows from the definition of $\varepsilon$. On the other hand, it follows from \eqref{eqn_pf:thm:convergence_as_opt_val_sol_1} that $C_k\as 0$ as $k\rightarrow\infty$. This gives the desired contradiction, implying that $D(\calXhat^N(\zb),\calX^\star(\zb))\as 0$ as $N\rightarrow\infty$.

\subsection{Proof of Theorem~\ref{thm:convergence_prob_est_quantile}} \label{appdx:thm:convergence_prob_est_quantile}

To lay the foundation of our analysis, we first define a few notation. Denote $\Ch_N(u,\xb;\zb)= \E_{\Probh_N(\zb)}\big[\rho_\tau\big(f(\xb,\xib) - u\big)\big]$ as the objective function of \eqref{eqn:cond_quantile_min_characterization_est}. Let $\mh^N(\xb;\zb)$ and $\calUhat^N(\xb;\zb)$ be the optimal value and the set of optimal solutions to problem~\eqref{eqn:cond_quantile_min_characterization_est}. Similarly, denote $C(u,\xb;\zb)=E_{\Prob(\zb)}\big[\rho_\tau\big(f(\xb,\xib) - u\big)\big]$ as the objective function of \eqref{eqn:cond_quantile_min_characterization}. Let $h^\star(\xb;\zb)$ and $\calU^\star(\xb;\zb)$ be the optimal value and the set of optimal solutions to problem~\eqref{eqn:cond_quantile_min_characterization}. By Assumption~\ref{assumption:regularity_cond_1}, we have $\calU^\star(\xb;\zb)=\{Q_\tau(\xb;\zb)\}$. In the following, we use ``$\Pro$'' to denote the probability measure induced by the data set $\calD_N$. We divide the proof into three steps. 

\textit{Step 1.} We show that there exists a compact set $\calU^\infty(\xb;\zb)\subset\R$ such that $Q_\tau(\xb;\zb)\in \calU^\infty(\xb;\zb)$ and $\Pro\big(\calUhat^N(\xb;\zb) \subseteq \calU^\infty(\xb;\zb)\big) \rightarrow 1$. Since $Q_\tau(\xb;\zb)$ is finite, if we find a compact set $\calU^\infty(\xb;\zb)\subset\R$ such that $\Pro\big(\calUhat^N(\xb;\zb) \subseteq \calU^\infty(\xb;\zb)\big) \rightarrow 1$, we can enlarge the set $\calUhat^N(\xb;\zb)$ to include $\calU^\infty(\xb;\zb)\subset\R$. Thus, it suffices to show that there exists a compact set $\calU^\infty(\xb;\zb)\subset\R$ such that $\Pro\big(\calUhat^N(\xb;\zb) \subseteq \calU^\infty(\xb;\zb)\big) \rightarrow 1$. Suppose, on the contrary, that such a compact set $\calU^\infty(\xb;\zb)$ does not exists.  Consider $u^\star\in\calU^\star(\xb;\zb)$ and a fixed $\varepsilon>0$. Since $\min_{u\in\R} \Ch_N(u,\xb;\zb)\leq \Ch_N(u^\star,\xb;\zb)$, we have
\begin{align*}
    \Pro\bigg(\min_{u\in\R} \Ch_N(u,\xb;\xib) \leq h^\star(\xb;\zb)+\varepsilon\bigg)&\geq \Pro\Big(\Ch_N(u^\star,\xb;\zb)\leq C(u^\star,\xb;\zb)+\varepsilon\Big) \\
    &\geq\Pro\Big(\Big|\Ch_N(u^\star,\xb;\zb)-C(u^\star,\xb;\zb)\Big|\leq\varepsilon\Big)\rightarrow1
\end{align*}
as $N\rightarrow\infty$, where the convergence follows from assumption~(b) in Theorem~\ref{thm:convergence_prob_est_quantile} that $\Ch_N(u,\xb;\zb)\convp C(u,\xb;\zb)$ for any $u\in\R$  and $\xb\in\calX$. Since $\varepsilon>0$ is arbitrary, we have $\Pro\big(\min_{u\in\R} \Ch_N(u,\xb;\zb)\leq h^\star(\xb;\zb)\big)\rightarrow1$.

To arrive at a contradiction, it suffices to show that $\Pro\big(\min_{u\in\R} \Ch_N(u,\xb;\zb)\geq h^\star(\xb;\zb)+1\big)\not\rightarrow 0$ as $N\rightarrow\infty$. Since $\calUhat^{N}(\xb;\zb)$ is the set of optimal solutions to \eqref{eqn:cond_quantile_min_characterization_est}, we have
\begin{equation} \label{eqn_pf:thm:convergence_prob_est_quantile_1.5}
    \Pro\bigg(\min_{u\in\R}\Ch_{N}(u,\xb;\zb)\geq h^\star(\xb;\zb)+1\bigg)=\Pro\bigg(\max_{u\in\calUhat_{N}(\xb;\zb)} \Ch_{N}(u,\xb;\zb)\geq h^\star(\xb;\zb)+1\bigg).
\end{equation}
If follows from the definition of $\Ch_{N}(u,\xb;\zb)$ that for any $u^N\in\calUhat^N(\xb;\zb)$, we have
\begin{align}
    &\quad\,\,\Ch_{N}(u^{N},\xb;\zb) \nonumber \\
    &=\E_{\Probh_{N}(\zb)}\Big[\rho_\tau\big(f(\xb,\xib) - u^N\big)\Big] \nonumber\\
    &=\E_{\Probh_{N}(\zb)}(\one_{B_{\zb}})\,\E_{\Probh_{N}(\zb)}\Big[\rho_\tau\big(f(\xb,\xib) - u^N\big)\,\,\Big|\,\, \xib\in B_{\zb}\Big]+ \E_{\Probh_{N}(\zb)}(\one_{B^c_{\zb}})\,\E_{\Probh_{N}(\zb)}\Big[\rho_\tau\big(f(\xb,\xib) - u^N\big)\,\,\Big|\,\, \xib\in B^c_{\zb}\Big] \nonumber\\
    &\geq\E_{\Probh_{N}(\zb)}(\one_{B_{\zb}})\,\E_{\Probh_{N}(\zb)}\Big[\rho_\tau\big(f(\xb,\xib) - u^N\big)\,\,\Big|\,\, \xib\in B_{\zb}\Big],  \label{eqn_pf:thm:convergence_prob_est_quantile_1}
\end{align}
where $B^c_{\zb}$ is the complement of the set $B_{\zb}$ and the last inequality follows from the fact that the check loss function $\rho_\tau$ is non-negative. Therefore, for any $\delta\in\big(0,\E_{\Prob(\zb)}(\one_{B_{\zb}})\big)$, where we recall by assumption~(a)  in Theorem~\ref{thm:convergence_prob_est_quantile} that $\E_{\Prob(\zb)}(\one_{B_{\zb}})>0$, we have 
\allowdisplaybreaks
\begin{subequations} \label{eqn_pf:thm:convergence_prob_est_quantile_2}
\begin{align}
    &\quad\,\,\Pro\bigg(\max_{u\in\calUhat^{N}(\xb;\zb)} \Ch_{N}(u,\xb;\zb)\geq v^\star+1\bigg) \nonumber \\
    &\geq\Pro\bigg(\max_{u\in\calUhat^{N}(\xb;\zb)}  \E_{\Probh_{N}(\zb)}(\one_{B_{\zb}})\,\E_{\Probh_{N}(\zb)}\Big[\rho_\tau\big(f(\xb,\xib) - u\big)\,\,\Big|\,\, \xib\in B_{\zb}\Big]\geq h^\star(\xb;\zb)+1\bigg)  \label{eqn_pf:thm:convergence_prob_est_quantile_2a}\\
    &\geq\Pro\bigg(\bigg\{ \E_{\Probh_{N}(\zb)}(\one_{B_{\zb}})\cdot\max_{u\in\calUhat^{N}(\xb;\zb)} \E_{\Probh_{N}(\zb)}\Big[\rho_\tau\big(f(\xb,\xib) - u\big)\,\,\Big|\,\, \xib\in B_{\zb}\Big] \geq h^\star(\xb;\zb)+1\bigg\} \nonumber \\
    &\hspace{30mm}\cap \bigg\{\max_{u\in\calUhat^{N}(\xb;\zb)} \E_{\Probh_{N}(\zb)}\Big[\rho_\tau\big(f(\xb,\xib) - u\big)\,\,\Big|\,\, \xib\in B_{\zb}\Big] \geq \frac{h^\star(\xb;\zb)+1}{ \E_{\Prob(\zb)}(\one_{B_{\zb}})-\delta}\bigg\}\bigg) \label{eqn_pf:thm:convergence_prob_est_quantile_2b} \\
    &\geq \Pro\bigg(\underbrace{\Big\{  \E_{\Probh_N(\zb)}(\one_{B_{\zb}})\geq \E_{\Prob(\zb)}(\one_{B_{\zb}})-\delta\Big\}}_{=:A_{1,N}} \nonumber \\
    & \hspace{30mm}\cap \underbrace{\bigg\{\max_{u\in\calUhat^{N}(\xb;\zb)} \E_{\Probh_{N}(\zb)}\Big[\rho_\tau\big(f(\xb,\xib) - u\big)\,\,\Big|\,\, \xib\in B_{\zb}\Big] \geq \frac{h^\star(\xb;\zb)+1}{ \E_{\Prob(\zb)}(\one_{B_{\zb}})-\delta}\bigg\}}_{=:A_{2,N}}\bigg).  \label{eqn_pf:thm:convergence_prob_est_quantile_2c}
\end{align}
\end{subequations}
Here, \eqref{eqn_pf:thm:convergence_prob_est_quantile_2a} follows from \eqref{eqn_pf:thm:convergence_prob_est_quantile_1}; \eqref{eqn_pf:thm:convergence_prob_est_quantile_2b} follows by intersecting the original event in \eqref{eqn_pf:thm:convergence_prob_est_quantile_2a} with an additional event $A_{2,N}$; \eqref{eqn_pf:thm:convergence_prob_est_quantile_2c} follows from replacing the original event in \eqref{eqn_pf:thm:convergence_prob_est_quantile_2a} with a subset $A_{1,N}$ given the event $A_{2,N}$.

We now analyze $\Pro(A_{1,N})$ and $\Pro(A_{2,N})$. First, we have
\begin{align}
    \Pro(A_{1,N})&=\Pro\Big(\E_{\Probh_N(\zb)}(\one_{B_{\zb}})\geq \E_{\Prob(\zb)}(\one_{B_{\zb}})-\delta\Big) \nonumber \\
    &\geq \Pro\Big(\Big|\E_{\Probh_N(\zb)}(\one_{B_{\zb}})-\E_{\Prob(\zb)}(\one_{B_{\zb}})\Big|\leq \delta\Big) \rightarrow 1  \label{eqn_pf:thm:convergence_prob_est_quantile_3}
\end{align}
as $N\rightarrow\infty$, where the convergence follows from assumption~(c) in Theorem~\ref{thm:convergence_prob_est_quantile}. Next, it follows from assumption~(a) in Theorem~\ref{thm:convergence_prob_est_quantile} that $\lim_{|u|\rightarrow\infty} \inf_{\xib\in B_{\zb}} \rho_\tau\big(f(\xb,\xib)-u\big)=\infty$. In other words, there exists $\rbar>0$ such that $ \rho_\tau\big(f(\xb,\xib)-u\big)\geq [h^\star(\xb;\zb)+1]/[\E_{\Prob(\zb)}(\one_{B_{\zb}})-\delta]$ for all $\xib\in B_{\zb}$ and $|u|>\rbar$. Denote $\B_{\rbar}=\{u\in\R\mid |u|\leq \rbar\}$ as the one-dimensional norm ball with radius $\rbar$. Note that if $\calUhat^N(\xb;\zb)\cap \B^c_{\rbar}\ne\emptyset$, i.e., there exists $\ubar\in\calUhat^N(\xb;\zb)$ such that $|\ubar|>\rbar$, then
$$\max_{u\in\calUhat^{N}(\xb;\zb)} \E_{\Probh_{N}(\zb)}\Big[\rho_\tau\big(f(\xb,\xib) - u\big)\,\,\Big|\,\, \xib\in B_{\zb}\Big]\geq \E_{\Probh_{N}(\zb)}\Big[\rho_\tau\big(f(\xb,\xib) - \ubar\big)\,\,\Big|\,\, \xib\in B_{\zb}\Big]\geq \frac{h^\star(\xb;\zb)+1}{ \E_{\Prob(\zb)}(\one_{B_{\zb}})-\delta}.$$
Therefore, we have
\begin{align}
    \Pro(A_{2,N})&=\Pro\Bigg(\max_{u\in\calUhat^{N}(\xb;\zb)} \E_{\Probh_{N}(\zb)}\Big[\rho_\tau\big(f(\xb,\xib) - u\big)\,\,\Big|\,\, \xib\in B_{\zb}\Big]\geq \frac{h^\star(\xb;\zb)+1}{ \E_{\Prob(\zb)}(\one_{B_{\zb}})-\delta}\Bigg) \nonumber \\
    &\geq \Pro\Big(\calUhat^N(\xb;\zb)\cap \B^c_{\rbar}\ne\emptyset\Big). \label{eqn_pf:thm:convergence_prob_est_quantile_4}
\end{align}
Recall the initial assumption that there does not exist a compact set $\calU^\infty(\xb;\zb)\subset\R$ such that $\Pro\big(\calUhat^N(\xb;\zb) \subseteq \calU^\infty(\xb;\zb)\big) \rightarrow 1$. This implies that there exists $N_{\rbar}\in\N$ and $\varepsilon>0$ such that $\Pro(\calUhat^N(\xb;\zb)\not\subseteq \B_{\rbar})\geq\varepsilon$ for all $N\geq N_{\rbar}$, i.e., we have $\Pro\big(\calUhat^N(\xb;\zb)\cap \B^c_{\rbar}\ne\emptyset\big)\geq\varepsilon$ for all $N\geq N_{\rbar}$. Hence, we have
\begin{subequations} \label{eqn_pf:thm:convergence_prob_est_quantile_5}
\begin{align}
    \Pro\bigg(\max_{u\in\calUhat^{N}(\xb;\zb)} \Ch_{N}(u,\xb;\zb)\geq v^\star+1\bigg) &\geq \Pro\big(A_{1,N}\cap A_{2,N}\big) \label{eqn_pf:thm:convergence_prob_est_quantile_5a}\\
    &\geq 1-\Prob\big(A^c_{1,N}\big)-\Prob\big(A^c_{2,N}\big) \label{eqn_pf:thm:convergence_prob_est_quantile_5b}\\
    &= \varepsilon - \Prob\big(A^c_{1,N}\big) \label{eqn_pf:thm:convergence_prob_est_quantile_5c}\\
    &\rightarrow\varepsilon \label{eqn_pf:thm:convergence_prob_est_quantile_5d}
\end{align}
\end{subequations}
as $N\rightarrow\infty$. Here, \eqref{eqn_pf:thm:convergence_prob_est_quantile_5a} follows from \eqref{eqn_pf:thm:convergence_prob_est_quantile_2}, \eqref{eqn_pf:thm:convergence_prob_est_quantile_5b} follows from the union bound, \eqref{eqn_pf:thm:convergence_prob_est_quantile_5c} follows from \eqref{eqn_pf:thm:convergence_prob_est_quantile_4}, and \eqref{eqn_pf:thm:convergence_prob_est_quantile_5d} follows from \eqref{eqn_pf:thm:convergence_prob_est_quantile_3}. It then follows from \eqref{eqn_pf:thm:convergence_prob_est_quantile_1.5} that
$$\Pro\bigg(\min_{u\in\R}\Ch_{N}(u,\xb;\zb)\geq h^\star(\xb;\zb)+1\bigg)=\Pro\bigg(\max_{u\in\calUhat_{N}(\xb;\zb)} \Ch_{N}(u,\xb;\zb)\geq h^\star(\xb;\zb)+1\bigg)\rightarrow\varepsilon >0,$$
which gives the desired contradiction.

\textit{Step~2.} We show that $\mh^N(\xb;\zb)\convp m^\star(\xb;\zb)$ as $N\rightarrow\infty$. Note that the check loss function $\rho_\tau\big(f(\xb,\xib)-u\big)$ is Lipschitz continuous in $u$ with Lipschitz constant $L=\max\{1-\tau,\tau\}$. Therefore, both $\Ch_N(\cdot,\xb;\zb)$ and $C(\cdot,\xb;\zb)$ are Lipschitz continuous with Lipschitz constant $L$, which implies that $\big|\Ch_N(\cdot,\xb;\zb)-C(\cdot,\xb;\zb)\big|$ is Lipschitz continuous with Lipschitz constant $2L$. Now, let $d_{\calU}<\infty$ be the diameter of $\calU^\infty(\xb;\zb)$. Consider any fixed $\varepsilon>0$. Since $\calU^\infty(\xb;\zb)$ is compact, there exists a finite cover of $\calU^\infty(\xb;\zb)$ by norm balls of radius $\varepsilon/(4L)$.  Denote such a finite cover as $\calU^{K}:=\{u_1,\dots,u_{K}\}$, where $K\leq 4\rho d_\calU L/\varepsilon$ for some constant $\rho>0$. Then, we have
\begin{subequations} \label{eqn_pf:thm:convergence_prob_est_quantile_6}
\begin{align}
    &\quad\,\,\Pro\big(\big|\mh^N(\xb;\zb)-m^\star(\xb;\zb)\big|>\varepsilon\big) \nonumber\\
    &=\Pro\Big(\big|\mh^N(\xb;\zb)-m^\star(\xb;\zb)\big|>\varepsilon,\,\, \calUhat^N(\xb;\zb)\not\subseteq \calU^\infty(\xb;\zb)\Big) \nonumber \\
    &\qquad +\Pro\Big(\big|\mh^N(\xb;\zb)-m^\star(\xb;\zb)\big|>\varepsilon\,\,\Big|\,\, \calUhat^N(\xb;\zb)\subseteq \calU^\infty(\xb;\zb)\Big) \Pro\Big(\calUhat^N(\xb;\zb)\subseteq \calU^\infty(\xb;\zb)\Big)   \label{eqn_pf:thm:convergence_prob_est_quantile_6a}\\
    &\leq \Pro\Big(\calUhat^N(\xb;\zb)\not\subseteq \calU^\infty(\xb;\zb)\Big) \nonumber \\
    &\qquad + \Pro\Bigg( \sup_{u\in\calU^\infty(\xb;\zb)} \bigg|  \E_{\Probh_N(\zb)}\Big[\rho_\tau\big(f(\xb,\xib)-u\big)\Big] - \E_{\Prob(\zb)}\Big[\rho_\tau\big(f(\xb,\xib)-u\big)\Big] \bigg| >\varepsilon \,\,\Bigg|\,\, \calUhat^N(\xb;\zb)\subseteq \calU^\infty(\xb;\zb)\Bigg) \label{eqn_pf:thm:convergence_prob_est_quantile_6b}\\
    &\leq \Pro\Big(\calUhat^N(\xb;\zb)\not\subseteq \calU^\infty(\xb;\zb)\Big) \nonumber \\
    &\qquad +\Bigg(\frac{4\rho d_\calU L}{\varepsilon}\Bigg)\cdot \sup_{u\in\calU^{K}}\Pro\Bigg( \Big|  \E_{\Probh_N(\zb)}\Big[\rho_\tau\big(f(\xb,\xib)-u\big)\Big] - \E_{\Prob(\zb)}\Big[\rho_\tau\big(f(\xb,\xib)-u\big)\Big] \Big| > \frac{\varepsilon}{2} \Bigg) \label{eqn_pf:thm:convergence_prob_est_quantile_6c}\\
    &\rightarrow 0 \label{eqn_pf:thm:convergence_prob_est_quantile_6d}
\end{align}
\end{subequations}
Here, \eqref{eqn_pf:thm:convergence_prob_est_quantile_6b} follows from the fact that 
$$\big|\mh^N(\xb;\zb)-m^\star(\xb;\zb)\big| \leq \sup_{u\in\calU^\infty(\xb;\zb)} \bigg|  \E_{\Probh_N(\zb)}\Big[\rho_\tau\big(f(\xb,\xib)-u\big)\Big] - \E_{\Prob(\zb)}\Big[\rho_\tau\big(f(\xb,\xib)-u\big)\Big] \bigg|$$
given the event $\calUhat^N(\xb;\zb)\subseteq \calU^\infty(\xb;\zb)$; \eqref{eqn_pf:thm:convergence_prob_est_quantile_6c} follows from Lemma~8 of \cite{Bertsimas_McCord:2019}; the convergence in \eqref{eqn_pf:thm:convergence_prob_est_quantile_6d} follows from assumption~(b) in Theorem~\ref{thm:convergence_prob_est_quantile}  and the facts that $\Pro\big(\calUhat^N(\xb;\zb) \subseteq \calU^\infty(\xb;\zb)\big) \rightarrow 1$ and $\calU^{K}$ is finite. This shows that $\mh^N(\xb;\zb)\convp m^\star(\xb;\zb)$ as $N\rightarrow\infty$.

\textit{Step 3.} We show that $\Qh_\tau^N(\xb;\zb)\convp Q_\tau(\xb;\zb)$ as $N\rightarrow\infty$. Since $\Qh_\tau^N(\xb;\zb)$ and $Q_\tau(\xb;\zb)$ are optimal solutions to problems~\eqref{eqn:cond_quantile_min_characterization_est} and \eqref{eqn:cond_quantile_min_characterization}, respectively, it suffices to show that $D\big(\calUhat^N(\xb;\zb),\calU^\star(\xb;\zb)\big)\convp 0$. Suppose, on the contrary, that $D\big(\calUhat^N(\xb;\zb),\calU^\star(\xb;\zb)\big)\not\convp 0$. Then, there exists $\delta>0$, $\varepsilon>0$, and a subsequence $\{N_k\}_{k\in\N}\subseteq\N$ such that $\Pro\big(D\big(\calUhat^{N_k}(\xb;\zb),\calU^\star(\xb;\zb)\big)>\delta\big)\geq \varepsilon$ for all $k\in\N$. Note that if $D\big(\calUhat^{N_k}(\xb;\zb),\calU^\star(\xb;\zb)\big)>\delta$, we have $\sup_{u\in\calUhat^{N_k}(\xb;\zb)} C(u,\xb;\zb)-m^\star(\xb;\zb)\geq \inf_{u\in\calU^\star_\delta(\xb;\zb)} C(u,\xb;\zb)-m^\star(\xb;\zb)=:\kappa_\delta>0$, where $\calU^\star_\delta(\xb;\zb)=\{u\in\R\mid \inf_{u^\star\in\calU^\star(\xb;\zb)}\norms{u-u^\star}\geq \delta\}$. This implies that $\Pro\big(\sup_{u\in\calUhat^{N_k}(\xb;\zb)} C(u,\xb;\zb)-m^\star(\xb;\zb)\geq\kappa_\delta\big)\geq \varepsilon$. On the other hand, we have
\begin{subequations} \label{eqn_pf:thm:convergence_prob_est_quantile_7}
\begin{align}
    &\quad\,\,\Pro\bigg(\sup_{u\in\calUhat^{N_k}(\xb;\zb)} C(u,\xb;\zb)-m^\star(\xb;\zb)\geq\kappa_\delta\bigg) \nonumber \\
    &\leq \Pro\bigg(\sup_{u\in\calUhat^{N_k}(\xb;\zb)} C(u,\xb;\zb) - \sup_{u\in\calU^\star(\xb;\zb)} \Ch_{N_k}(u,\xb;\zb)\geq\frac{\kappa_\delta}{2}\bigg) + \Pro\bigg(\sup_{u\in\calU^\star(\xb;\zb)} \Ch_{N_k}(u,\xb;\zb)-m^\star(\xb;\zb)\geq\frac{\kappa_\delta}{2}\bigg) \label{eqn_pf:thm:convergence_prob_est_quantile_7a}\\
    &\leq  \Pro\bigg(\sup_{u\in\calUhat^{N_k}(\xb;\zb)} C(u,\xb;\zb)-\mh^{N_k}(\xb;\zb) \geq \frac{\kappa_\delta}{2}\bigg) + \Pro\bigg(\sup_{u\in\calU^\star(\xb;\zb)} \Ch_{N_k}(u,\xb;\zb)-m^\star(\xb;\zb)\geq\frac{\kappa_\delta}{2}\bigg) \label{eqn_pf:thm:convergence_prob_est_quantile_7b} \\
    &\leq  \Pro\bigg(\sup_{u\in\calUhat^{N_k}(\xb;\zb)} \big|C(u,\xb;\zb)- \Ch_{N_k}(u,\xb;\zb)\big|\geq\frac{\kappa_\delta}{2}\bigg) + \Pro\bigg(\sup_{u\in\calU^\star(\xb;\zb)} \big|\Ch_{N_k}(u,\xb;\zb)-C(u,\xb;\zb)\big|\geq\frac{\kappa_\delta}{2}\bigg) \label{eqn_pf:thm:convergence_prob_est_quantile_7c}\\
    &\rightarrow0 \label{eqn_pf:thm:convergence_prob_est_quantile_7d}
\end{align}    
\end{subequations}
as $k\rightarrow\infty$. Here, \eqref{eqn_pf:thm:convergence_prob_est_quantile_7a} follows from the union bound; \eqref{eqn_pf:thm:convergence_prob_est_quantile_7b} follows from the fact that $\sup_{u\in\calU^\star(\xb;\zb)} \Ch_{N_k}(u,\xb;\zb)\leq \sup_{u\in\R} \Ch_{N_k}(u,\xb;\zb)=\mh^{N_k}(\xb;\zb)$;  \eqref{eqn_pf:thm:convergence_prob_est_quantile_7c} follows from the facts that $\sup_{u\in\calUhat^{N_k}(\xb;\zb)} \Ch_{N_k}(u,\xb;\zb)=\mh^{N_k}(\xb;\zb)$ and $\sup_{u\in\calU^\star(\xb;\zb)} C(u,\xb;\zb)=m^\star(\xb;\zb)$; the convergence in \eqref{eqn_pf:thm:convergence_prob_est_quantile_7d} follows from a similar argument as in \eqref{eqn_pf:thm:convergence_prob_est_quantile_6}. This gives the desired contradiction, showing that  $\Qh_\tau^N(\xb;\zb)\convp Q_\tau(\xb;\zb)$ as $N\rightarrow\infty$. This completes the proof.

\subsection{Proof of Lemma~\ref{lem:unif_convergence_prob_obj}} \label{appdx:lem:unif_convergence_prob_obj}

By Assumption~\ref{assumption:compact_X}, the feasible set $\calX$ is compact. Thus, there exists a finite cover of $\calX$ by norm balls of radius $\min\{\nu,1\}$ for any $\nu>0$. Denote the finite cover as $\calX^K=\{\xb_k\}_{k=1}^K$ with $K\leq (\rho d_\calX/\min\{\nu,1\})^n$ for some $\rho>0$ that depends on the norm used, where $d_\calX$ is the diameter of $\calX$. Define the function $g^N(\xb;\zb)= \big|\Qh^N_\tau(\xb';\zb)-Q_\tau(\xb';\zb)\big|$. Note that for any $\{\xb',\xb''\}\subseteq\calX$ such that $\norms{\xb'-\xb''}\leq 1$, we have
\begin{subequations} \label{eqn_pf:lem:unif_convergence_prob_obj_1}
\begin{align}
    \big|g^N(\xb';\zb)-g^N(\xb'';\zb)\big|&=\bigg| \Big|\Qh^N_\tau(\xb';\zb)-Q_\tau(\xb';\zb)\Big| - \Big|\Qh^N_\tau(\xb'';\zb)-Q_\tau(\xb'';\zb)\Big| \bigg| \\
    &\leq \bigg| \Big[\Qh^N_\tau(\xb';\zb)-Q_\tau(\xb';\zb)\Big] - \Big[\Qh^N_\tau(\xb'';\zb)-Q_\tau(\xb'';\zb)\Big] \bigg| \label{eqn_pf:lem:unif_convergence_prob_obj_1a}\\
    &\leq \Big|\Qh^N_\tau(\xb';\zb)-\Qh^N_\tau(\xb'';\zb) \Big|  + \Big| Q_\tau(\xb';\zb) -Q_\tau(\xb'';\zb)\Big| \label{eqn_pf:lem:unif_convergence_prob_obj_1b}\\
    &\leq L_2\norms{\xb'-\xb''}^\beta + L_1\norms{\xb'-\xb''}^\alpha  \label{eqn_pf:lem:unif_convergence_prob_obj_1c}\\
    &\leq (L_1+L_2) \norms{\xb'-\xb''}^{\min\{\alpha,\beta\}}, \label{eqn_pf:lem:unif_convergence_prob_obj_1d}
\end{align}
\end{subequations}
where \eqref{eqn_pf:lem:unif_convergence_prob_obj_1a} follows from the reverse triangle inequality; \eqref{eqn_pf:lem:unif_convergence_prob_obj_1b} follows from the triangle inequality; \eqref{eqn_pf:lem:unif_convergence_prob_obj_1c} follows from Assumption~\ref{assumption:Holder_cont}; \eqref{eqn_pf:lem:unif_convergence_prob_obj_1d} follows from the fact that if $t\in[0,1]$, we have $t^{\gamma_1}>t^{\gamma_2}$ for any $0<\gamma_1<\gamma_2$. Now, note that for any $\xb\in\calX$, there exists some $k'\in\{1,\dots,K\}$ such that $\norms{\xb-\xb_{k'}}\leq \min\{\nu,1\}$. Therefore,
$$g^N(\xb;\zb)\leq g^N(\xb_{k'};\zb)+\big|g^N(\xb_{k'};\zb)-g^N(\xb;\zb)\big|\leq  g^N(\xb_{k'};\zb)+(L_1+L_2) \min\{\nu,1\}^{\min\{\alpha,\beta\}},$$
where the last inequality follows from \eqref{eqn_pf:lem:unif_convergence_prob_obj_1}. Therefore, we have
\begin{equation}  \label{eqn_pf:lem:unif_convergence_prob_obj_2}
    \sup_{\xb\in\calX} g^N(\xb;\zb) \leq \sup_{k\in\{1,\dots,K\}} g^N(\xb_{k};\zb) + (L_1+L_2) \min\{\nu,1\}^{\min\{\alpha,\beta\}}.
\end{equation}
Therefore, using~\eqref{eqn_pf:lem:unif_convergence_prob_obj_2}, for any $\varepsilon>0$, we have
\begin{equation}
    \Pro\bigg(\sup_{\xb\in\calX} g^N(\xb;\zb)>\varepsilon \bigg)\leq \Pro\bigg(\sup_{k\in\{1,\dots,K\}} g^N(\xb_{k};\zb) + (L_1+L_2) \min\{\nu,1\}^{\min\{\alpha,\beta\}}>\varepsilon\bigg).
\end{equation}
Now, let us pick $\nu=\big\{\varepsilon/[\phi(L_1+L_2)]\big\}^{1/\min\{\alpha,\beta\}}$ for any $\phi>\max\{2,\varepsilon/(L_1+L_2)\}$. Note that by construction, we have $\nu<1$, and hence, $(L_1+L_2)\min\{\nu,1\}^{\min\{\alpha,\beta\}}\leq \varepsilon/2$. Therefore, we obtain
\begin{align*}
    \Pro\bigg(\sup_{\xb\in\calX} g^N(\xb;\zb)>\varepsilon \bigg)&\leq \Pro\bigg(\sup_{k\in\{1,\dots,K\}} g^N(\xb_{k};\zb) + \frac{\varepsilon}{2}>\varepsilon\bigg)\\
    &=\Pro\bigg(\sup_{k\in\{1,\dots,K\}} g^N(\xb_{k};\zb) >\frac{\varepsilon}{2}\bigg)\\
    &\leq \Bigg[\rho d_\calX \bigg(\frac{\phi(L_1+L_2)}{\varepsilon}\bigg)^{\frac{1}{\min\{\alpha,\beta\}}}\Bigg]^{n}\sup_{\xb\in\calX^K}\Prob\bigg(g^N(\xb;\zb)>\frac{\varepsilon}{2}\bigg)\\
    &\rightarrow 0
\end{align*}
where the last inequality follows from the union bound; the convergence follows from the facts that $\calX^K$ is finite and $g^N(\xb;\zb)\convp 0$ as $N\rightarrow\infty$, which is shown in Theorem~\ref{thm:convergence_prob_est_quantile}. This shows that $\Qh_\tau^N(\xb;\zb)\convp Q_\tau(\xb;\zb)$ uniformly over $\xb\in\calX$.

\subsection{Proof of Theorem~\ref{thm:convergence_prob_opt_val_sol}} \label{appdx:thm:convergence_prob_opt_val_sol}

For part~(i), it follows from Lemma~\ref{lem:unif_convergence_prob_obj} that for any $\varepsilon>0$, we have
\begin{align*}
    \Pro\Big(\big|\vhat^N(\zb)-v^\star\big|>\varepsilon\Big)&=\Pro\Bigg(\bigg|\min_{\xb\in\calX} \Qh^N_\tau(\xb;\zb)-\min_{\xb\in\calX}Q_\tau(\xb;\zb)\bigg|>\varepsilon \Bigg) \\
    &\leq\Pro\Bigg(\sup_{\xb\in\calX}\Big|\Qh^N_\tau(\xb;\zb)-Q_\tau(\xb;\zb)\Big|>\varepsilon \Bigg)\rightarrow 0.
\end{align*}
This shows that $\vhat^N(\zb) \convp  v^\star(\zb)$ as $N\rightarrow\infty$. Next, we prove part~(ii). Following a similar argument in \eqref{eqn_pf:thm:convergence_prob_est_quantile_7}, for any $\varepsilon>0$, we have
%
\begin{align}
    &\quad\,\,\Pro\bigg(\sup_{\xb\in\calXhat^{N}(\zb)} Q_\tau(\xb;\zb)-v^\star(\zb)\geq\varepsilon\bigg) \nonumber \\
    &\leq \Pro\bigg(\sup_{\xb\in\calXhat^{N}(\zb)} Q_\tau(\xb;\zb) - \sup_{\xb\in\calX^\star(\zb)} \Qh^{N}_\tau(\xb;\zb)\geq\frac{\varepsilon}{2}\bigg) + \Pro\bigg(\sup_{\xb\in\calX^\star(\xb)} Q_\tau^N(\xb;\zb)-v^\star(\zb)\geq\frac{\varepsilon}{2}\bigg)\nonumber \\
    &\leq  \Pro\bigg(\sup_{\xb\in\calXhat^{N}(\xb)} Q_\tau(\xb;\zb)-\vhat^{N}(\zb) \geq \frac{\varepsilon}{2}\bigg) + \Pro\bigg(\sup_{\xb\in\calX^\star(\zb)} \Qh^{N}_\tau(\xb;\zb)-v^\star(\zb)\geq\frac{\varepsilon}{2}\bigg) \nonumber\\
    &\leq  \Pro\bigg(\sup_{\xb\in\calXhat^{N}(\zb)} \big|Q_\tau(\xb;\zb)- \Qh^{N}_\tau(\xb;\zb)\big|\geq\frac{\varepsilon}{2}\bigg) + \Pro\bigg(\sup_{\xb\in\calX^\star(\zb)} \big|\Qh^{N}_\tau(\xb;\zb)-Q_\tau(\xb;\zb)\big|\geq\frac{\varepsilon}{2}\bigg) \nonumber\\
    &\rightarrow0,  \label{eqn_pf:thm:convergence_prob_opt_val_sol_1}
\end{align}    
%
where the convergence again follows from Lemma~\ref{lem:unif_convergence_prob_obj}. This shows that  $\sup_{\xb\in\calXhat^N(\zb)} Q_\tau(\xb^N;\zb)  \convp v^\star(\zb)$ as $N\rightarrow\infty$. Finally, we prove part~(iii).  Suppose, on the contrary, that $D\big(\calXhat^N(\zb),\calX^\star(\zb)\big)\not\convp 0$. Then, there exists $\delta>0$, $\varepsilon>0$, and a subsequence $\{N_k\}_{k\in\N}\subseteq\N$ such that $\Pro\big(D\big(\calXhat^{N}(\zb),\calX^\star(\zb)\big)>\delta\big)\geq\varepsilon$ for all $k\in\N$. Note that if $D\big(\calXhat^{N}(\zb),\calX^\star(\zb)\big)>\delta$, we have $\sup_{\xb\in\calXhat^{N}(\zb)} Q_\tau(\xb;\zb)-v^\star(\zb)\geq \inf_{\xb\in\calX^\star_\delta(\zb)} Q_\tau(\xb;\zb)-v^\star(\zb)=:\kappa_\delta>0$, where $\calX^\star_\delta(\zb)=\{\xb\in\calX\mid \inf_{\xb^\star\in\calX^\star(\zb)}\norms{\xb-\xb^\star}\geq \delta\}$. This implies that $\Pro\big(\sup_{\xb\in\calXhat^{N}(\zb)} Q_\tau(\xb;\zb)-v^\star(\zb)\geq\kappa_\delta\big)\geq\varepsilon$, which contradicts with \eqref{eqn_pf:thm:convergence_prob_opt_val_sol_1}. This shows that $D(\calXhat^N(\zb),\calX^\star(\zb))=\sup_{\xb^N\in\calXhat^N(\zb)}\inf_{\xb^\star\in\calX^\star(\zb)}\norms{\xb^\star-\xb^N}\convp 0$ as $N\rightarrow\infty$.

\subsection{Proof of Proposition \ref{prop:MILP_quan_min_LP_form}} \label{appdx:prop_MILP_quan_min_LP_form}

We first show that for any optimal solution to \eqref{model:MILP_quan_min}, there is a feasible solution to \eqref{model:MILP_quan_min_LP_form} with the same objective value. Let $(\xb^\star,t^\star,\vb^\star)$ be an optimal solution to \eqref{model:MILP_quan_min}. Then, there exists $\yb^*_i\in\calY(\xb^\star,\xib^i)$ such that $\qb^\tp\yb^\star_i = f(\xb^\star,\xib^i)\leq M(1-v^\star_i)+t^\star$ for $i\in\{1,\dots,N\}$. Thus, this shows that $(\xb^\star,t^\star,\vb^\star,\yb^\star)$ is a feasible solution to \eqref{model:MILP_quan_min_LP_form} with the same objective value $t^\star$. Next, we show that for any optimal solution to \eqref{model:MILP_quan_min_LP_form}, there is a feasible solution to \eqref{model:MILP_quan_min} with the same objective value. Let $(\xb^\star,t^\star,\vb^\star,\yb^\star)$ be an optimal solution to \eqref{model:MILP_quan_min_LP_form}. Since $\yb^\star_i\in\calY(\xb^\star,\xib^i)$, we have $f(\xb,\xib^i)\leq\qb^\tp\yb^\star_i$. It follows that, by constraints \eqref{model:MILP_quan_min_LP_form_con1}, $f(\xb^\star,\xib^i)\leq\qb^\tp\yb^\star_i\leq M(1-v^\star_i)+t^\star$. Thus, $(\xb^\star,t^\star,\vb^\star)$ is feasible to \eqref{model:MILP_quan_min} with the same objective value $t^\star$. This also implies that $(\xb^\star,t^\star,\vb^\star)$ is an optimal solution to \eqref{model:MILP_quan_min}.

\subsection{Proof of Theorem \ref{thm:SiCG_convergence}} \label{appdx:thm_SiCG_convergence}

The proof follows a similar argument in \cite{Tsang_et_al:2023} with modifications due to the stochasticity in the exploration step. First, we adopt two fundamental results, Propositions 1 and 2 from \cite{Tsang_et_al:2023}, which provide guarantees on the optimality gap with the employed backtracking routine. Following the same proof logic of Proposition 1 in \cite{Tsang_et_al:2023}, we have that $L^\ell$ is a valid lower bound a.s. at any iteration. Therefore, the quantity $(\Ubar-L^\ell)/\Ubar$ is a valid measure of the optimality gap (since $L^\ell$ and $\Ubar$ are lower and upper bounds of $\widehat{\upsilon}^N(\zb)$). Moreover, adopting from Proposition 2 in \cite{Tsang_et_al:2023}, if SiCG reaches the exploitation step at iteration $j$, then the optimality gap $(\Ubar-L^\ell)/\Ubar$ is upper bounded by  $(1-\epst)^{-1} \prod_{k=\ell}^j(1-\varepsilon_{MP}^k)^{-1}-1$ a.s., which is related to the parameter $\epst$ (governing the trade-off between exploitation and exploration) and the MIP gap imposed in the master problems.

Next, to show the a.s. finite convergence of SiCG, we first claim that at each exploration step, the current dual set $\Pih$ is enlarged by one new element (i.e., a new dual vertex). Define $\calV=\{\vb\in\{0,1\}^N\mid \sum_{i=1}^N w_iv_i \geq \tau\}$. Indeed, if at iteration $j$, we have $f(\xb^j,\xib^i)=\max_{\pib\in\Pih}\{(\hb-\Tb\xb^j-\Cb\xib^i)^\tp\pib\}$ for all $i\in\{1,\dots,N\}$, then, by our construction of $U^j$, we have
\begin{align}
    U^j = t^j &\geq \max_{i=1,\dots,N}\Bigg\{\max_{\pib\in\Pih} \Big\{(\hb-\Tb\xb^j-\Cb\xib^i)^\tp\pib \Big\} -M(1-v_i^j)\Bigg\} \label{eqn:pf_thm_SiCG_convergence_1} \\
    &\geq \min_{\vb\in\calV} \Bigg\{ \max_{i=1,\dots,N} \bigg\{\max_{\pib\in\Pih} \Big\{(\hb-\Tb\xb^j-\Cb\xib^i)^\tp\pib\Big\} -M(1-v_i) \bigg\} \Bigg\} \label{eqn:pf_thm_SiCG_convergence_2}\\
    &= \Qh^N_\tau\Big(\max_{\pib\in\Pih} \big\{(\hb-\Tb\xb^j-\Cb\xib)^\tp\pib \big\} \,\Big|\, \Zb=\zb \Big) \label{eqn:pf_thm_SiCG_convergence_3}\\
    &= \Qh^N_\tau(\xb^j;\zb), \nonumber
\end{align}
where the inequality in \eqref{eqn:pf_thm_SiCG_convergence_1} follows from constraints \eqref{model:MP_quan_min_SiCG_con1}, \eqref{eqn:pf_thm_SiCG_convergence_2} comes from $\vb^j\in\calV$, \eqref{eqn:pf_thm_SiCG_convergence_3} is based on the fact that the $\tau$th quantile of $\max_{\pib\in\Pih} \big\{(\hb-\Tb\xb^j-\Cb\xib)^\tp\pib \big\}$ is the optimal value of \eqref{eqn:pf_thm_SiCG_convergence_2}. Hence, this immediately implies that
$$\Ubar-U^j \leq  \Ubar-\Qh^N_\tau(\xb^j;\zb)\leq 0.$$
Thus, in such a case, the algorithm proceeds to the exploitation step. As a result, if, at iteration $j$, the algorithm proceeds to the exploration step (and thus, not the exploitation step), then there exists $i'\in\{1,\dots,N\}$ such that $$f(\xb^j,\xib^{i'})=\max_{\pib\in\Pi}\{(\hb-\Tb\xb^j-\Cb\xib^{i'})^\tp\pib\}>\max_{\pib\in\Pih}\{(\hb-\Tb\xb^j-\Cb\xib^{i'})^\tp\pib\},$$
which implies that $\Pih$ does not include all the vertices of the dual feasible set $\Pi$. By our assumption, the dual feasible set consists of a finite number of extreme points, and we always randomly select one of the vertices in the exploration step. Therefore, at each exploration step, a new dual variable will be found and the set $\Pih$ is enlarged. 

Finally, we prove the a.s. finite convergence of SiCG. Suppose, on the contrary, that SiCG does not terminate in a finite number of iterations. Note that if SiCG does not terminate, at iteration, it must proceed to either the exploration or the exploitation step. However, as we have just shown, the set $\Pih$ is enlarged by a vertex at each exploration step, and $\Pi$ consists of a finite number of vertices only, it implies that SiCG proceeds to the exploitation step infinitely many times. This corresponds to the situation that the master problem with a given set of vertices $\Pit$ is solved infinitely many times (and SiCG always proceeds to the exploitation step). However, as mentioned earlier, the optimality gap $(\Ubar-L^\ell)/\Ubar$ is a.s. upper bounded by $(1-\epst)^{-1} (1-\varepsilon_{MP}^\ell)^{-1}-1$. Moreover, $\varepsilon_{MP}^j$ is decreased at each exploitation step, and hence, $\varepsilon_{MP}^j$ converges to zero. Thus, the upper bound $(1-\epst)^{-1} (1-\varepsilon_{MP}^\ell)^{-1}-1$ converges to $(1-\epst)^{-1}-1 =\epst/(1-\epst)$. By our assumption that $\epst<\varepsilon/(1+\varepsilon)$, this shows that after a sufficiently large number of iterations, the termination condition will be satisfied. Therefore, this contradicts that SiCG does not terminate in a finite number of iterations.

\section{SiCG for the Appointment Scheduling Problem} \label{appdx:SiCG_ASP}

In this section, we present the SiCG implemented in our experiments. First, note that for a given $\xb$, the objective value $f(\xb,\xib)$ in \eqref{model:ASP} can be easily computed without solving the optimization problem. Specifically, the optimal solutions are $u_{i-1}=\max\{x_{i-1}-s_{i-1}-w_{i-1},0\}$ and $w_i = \max\{s_{i-1}+w_{i-1}-x_{i-1},0\}$ for all $i\in\{2,\dots,n+1\}$. As a result, evaluating $\Qh^N_\tau(\xb^j;\zb)$ in the subproblem step reduces to direct computations. Moreover, note that the dual of $f(\xb,\xib)$ in \eqref{model:ASP} reads as follows:
\begin{subequations} 
\begin{align}
f(\xb,\sb)=\underset{\yb}{\text{maximize}\,} \quad
&  \sum_{i=1}^n (s_i-x_i ) y_i     \label{model:appt_sch_dual_obj} \\    
\text{subject to} \quad
&  -\cu_n \leq y_n \leq \co,  \label{model:appt_sch_dual_con1} \\    
&  -\cu_{i-1}\leq y_{i-1} \leq y_i+\cw_i,\quad\forall i\in\{2,\dots,n\}. \label{model:2nd_LP_form_dual_con2}   
\end{align} \label{model:ASP_dual}%
\end{subequations} 
As pointed out in \cite{Jiang_et_al:2017}, there is a one-to-one correspondence between the extreme point of the dual feasible set and a partition of $\{1,\dots,n+1\}$. In particular, for a given interval $\{k,\dots,j\}$ in the partition, we set
$$y_i = \begin{cases} -\cu_j + \sum_{\ell=i+1}^j \cw_\ell &\text{if } 1\leq i\leq j \leq n, \\ \co +\sum_{\ell=i+1}^n \cw_\ell &\text{if } 1\leq i \leq n,\, j = n+1,\end{cases}$$
for $i\in\{k,\dots,j\}$. As a result, to generate a dual feasible solution, it suffices to generate a random partition of $\{1,\dots,n+1\}$. Thus, we modify the exploration step as in Algorithm \ref{algo:SiCG_exploration_ASP}. We start by searching for some potential dual solutions. First, we record, from Step 2 of SiCG, the scenario \texttt{idx} that yields the true quantile $\Qh^N_\tau(\xb^j;\zb)$. We then evaluate $f(\xb^j;\xib^\texttt{idx})$ for a dual solution $\pib_\texttt{idx}$, which serves as a potential new scenario to be added into the current scenario set $\Pih$. Moreover, we also record, from Step 1 of SiCG, the scenario \texttt{idx'} that yields the quantile for $\max_{\pib\in\Pih}\{ (\hb-\Tb\xb^j-\Cb\xib)^\tp\pib\}$. Again, we evaluate $f(\xb^j;\xib^\texttt{idx'})$ for a dual solution $\pib_\texttt{idx'}$, which serves as another potential new scenario to be added into the current scenario set $\Pih$. If these two dual solutions are already included in $\Pih$, we start generating dual solutions until we find a new one. 

\LinesNumbered  
\IncMargin{1em}
\begin{algorithm}[t]
\DontPrintSemicolon  
\SetKwInOut{Initialization}{Initialization}
\Initialization{Scenario index \texttt{idx} that yields the true quantile (in Step 2 of SiCG); Scenario index \texttt{idx'} that yields the objective value of the master problem (in Step 1 of SiCG)}
Compute the dual solution $\pib_\texttt{idx}$ associated with scenario \texttt{idx}. \tcc*[r]{Trial} 
If $\pib_\texttt{idx}\not\in\Pih$, then \textbf{return} $\pib_\texttt{idx}$. \\
Compute the dual solution $\pib_\texttt{idx'}$ associated with scenario \texttt{idx'}. \tcc*[r]{Trial} 
If $\pib_\texttt{idx'}\not\in\Pih$, then \textbf{return} $\pib_\texttt{idx'}$. \\
\While{\textup{\texttt{true}}}{
    Generate a random partition of $\{1,\dots,n+1\}$ and obtain the associated dual solution $\pib$. \tcc*[r]{Generate dual solutions} 
    If $\pib\not\in\Pih$, then \textbf{return} $\pib$.
}
\BlankLine
\caption{Exploration Step in SiCG for ASP}\label{algo:SiCG_exploration_ASP}
\end{algorithm}\DecMargin{1em}

\section{Additional Computational Results} \label{appdx:add_comp_results}

\subsection{Additional Results on Computational Time} \label{appdx:expt_compt_time}

\begin{figure}[t] 
\centering  
\includegraphics[scale=0.75]{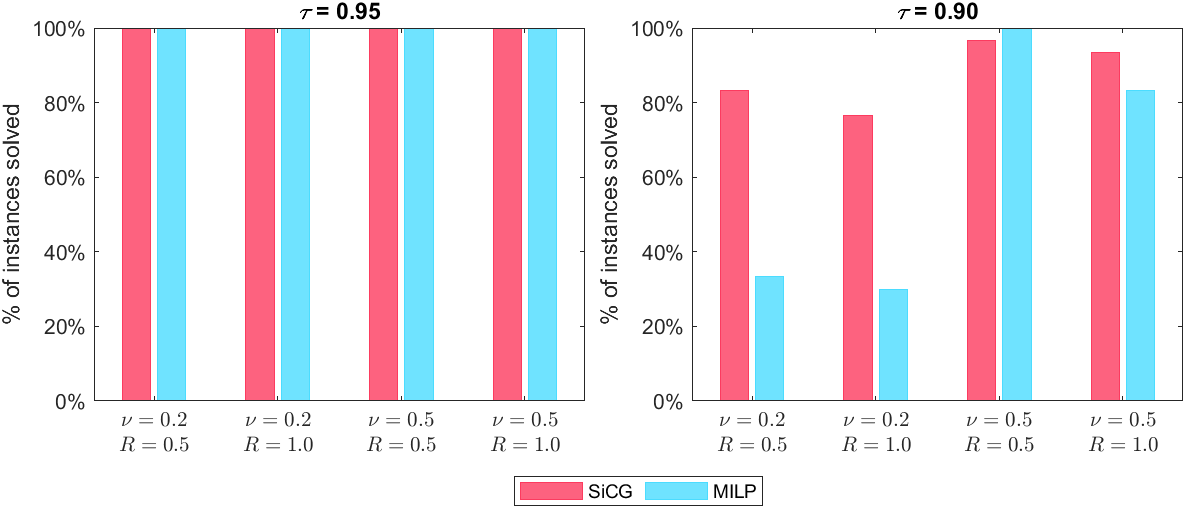}
\caption{Percentage of instances (over $30$ replications) solved within $1$ hour time limit using SiCG and MILP for $\tau\in\{0.95,0.90\}$ when $n=6$ and $N=500$.} \label{fig:num_solved_n8}
\end{figure}

In this section, we present similar computational time results when there are $n=8$ appointments. First, we note that when $N=1000$, the number of instances solved within $1$ hour for both SiCG and MILP are similar when $\tau=0.95$, and most of the instances become challenging when $\tau=0.90$. Therefore, we focus on $N\in\{200,500\}$ in our discussion. Figure \ref{fig:num_solved_n8} shows the percentage of instances solved within $1$ hour time limit. When $\tau=0.90$, we observe that the number of instances solved by SiCG is generally larger. In particular, when $\nu=0.2$ (more computationally challenging instances), SiCG can solve around $80\%$ of the instances but MILP can only solve less than $40\%$ of them.

Tables \ref{table:sol_time_n8_tau95}--\ref{table:sol_time_n8_tau90} present the lower quantile, median, and upper quantile of the solution times (over $30$ generated instances) under different number of scenarios and settings under $\tau=0.95$ and $\tau=0.90$, respectively. We observe that when $N=200$, SiCG may take slightly longer time to solve due to its iterative nature. However, when $N=500$, we could observe improvements in computational efficiency, especially when $\tau=0.90$. For instance, when $\nu=0.5$ and $R=1.0$, while SiCG requires less than $1500$ seconds to solve $75\%$ of the instances, MILP takes around $2600$ seconds. These results further demonstrate that our SiCG could have a better computational performance than MILP.
\begin{table}[t]\centering  \small
\ra{1.0}  
\caption{Solution time (in seconds) when $n=8$ and $\tau=0.95$ under different number of scenarios $N\in\{200,500\}$. Q1, Q2, and Q3, are the lower quantile, median, and upper quantile of solution times over $30$ generated instances.} \label{table:sol_time_n8_tau95}
\begin{tabular}{ll||rrr|rrr} \Xhline{1.0pt}
\multicolumn{2}{c||}{$N = 200$} & \multicolumn{3}{c|}{SiCG} & \multicolumn{3}{c}{MILP} \\
\multicolumn{2}{c||}{}          & Q1     & Q2     & Q3     & Q1     & Q2     & Q3     \\  \hline
$\nu = 0.2$    & $R = 0.5$    & 3.40   & 4.52   & 6.40   & 2.32   & 2.89   & 3.39   \\
$\nu = 0.2$    & $R = 1.0$    & 3.58   & 5.08   & 7.00   & 3.06   & 3.65   & 4.33   \\
$\nu = 0.5$    & $R = 0.5$    & 2.29   & 2.89   & 4.58   & 1.58   & 1.77   & 1.97   \\
$\nu = 0.5$    & $R = 1.0$    & 2.49   & 3.02   & 4.03   & 2.40   & 2.61   & 3.21   \\ \Xhline{1.0pt}
\multicolumn{2}{c||}{$N = 500$} & \multicolumn{3}{c|}{SiCG} & \multicolumn{3}{c}{MILP} \\
\multicolumn{2}{c||}{}          & Q1     & Q2     & Q3     & Q1     & Q2     & Q3     \\  \hline
$\nu = 0.2$    & $R = 0.5$    & 56.20  & 93.51  & 134.18 & 104.49 & 115.10 & 141.95 \\
$\nu = 0.2$    & $R = 1.0$    & 55.67  & 127.23 & 216.37 & 111.20 & 130.40 & 161.49 \\
$\nu = 0.5$    & $R = 0.5$    & 16.42  & 23.63  & 33.15  & 13.70  & 17.44  & 45.37  \\
$\nu = 0.5$    & $R = 1.0$    & 28.42  & 40.47  & 67.61  & 37.57  & 76.36  & 92.17  \\
\Xhline{1.0pt}
\end{tabular}
\end{table}

\begin{table}[t]\centering  \small
\ra{1.0}  
\caption{Solution time (in seconds) when $n=8$ and $\tau=0.90$ under different number of scenarios $N\in\{200,500\}$. Q1, Q2, and Q3, are the lower quantile, median, and upper quantile of solution times over $30$ generated instances.} \label{table:sol_time_n8_tau90}
\begin{tabular}{ll||rrr|rrr} \Xhline{1.0pt}
\multicolumn{2}{c||}{$N = 200$} & \multicolumn{3}{c|}{SiCG}   & \multicolumn{3}{c}{MILP}    \\
\multicolumn{2}{c||}{}          & Q1     & Q2      & Q3      & Q1      & Q2      & Q3      \\ \hline
$\nu = 0.2$    & $R = 0.5$    & 9.19   & 11.61   & 16.66   & 6.33    & 8.20    & 13.80   \\
$\nu = 0.2$    & $R = 1.0$    & 10.07  & 19.18   & 35.44   & 9.03    & 11.55   & 20.32   \\
$\nu = 0.5$    & $R = 0.5$    & 5.42   & 7.51    & 9.90    & 4.62    & 6.01    & 7.74    \\
$\nu = 0.5$    & $R = 1.0$    & 5.75   & 10.06   & 17.51   & 7.46    & 9.44    & 12.08   \\ \Xhline{1.0pt}
\multicolumn{2}{c||}{$N = 500$} & \multicolumn{3}{c|}{SiCG}   & \multicolumn{3}{c}{MILP}    \\
\multicolumn{2}{c||}{}          & Q1     & Q2      & Q3      & Q1      & Q2      & Q3      \\ \hline
$\nu = 0.2$    & $R = 0.5$    & 594.90 & 975.29  & 1918.81 & 2294.71 & $>$3600 & $>$3600 \\
$\nu = 0.2$    & $R = 1.0$    & 816.26 & 1122.13 & 3413.37 & 2836.38 & $>$3600 & $>$3600 \\
$\nu = 0.5$    & $R = 0.5$    & 167.58 & 274.20  & 396.83  & 303.92  & 587.83  & 782.15  \\
$\nu = 0.5$    & $R = 1.0$    & 407.25 & 655.35  & 1436.33 & 993.53  & 1502.39 & 2619.30 \\
\Xhline{1.0pt}
\end{tabular}
\end{table}

\subsection{Effect of Sub-Sampling} \label{appdx:expt_subsample}

In this section, we discuss the effect of sub-sampling on our experiments in Sections \ref{subsec:opt_schedule}--\ref{subsec:OS_performance}. In this experiment, we follow the experiment description in Section \ref{subsec:opt_schedule} to generate a set of $N=10,000$ pairs historical data. Then, we randomly generate $20$ sets of sub-samples of size $N'\in\{10,20,50,100,200,500,1000\}$ and solve the corresponding CSO or SAA model for an optimal solution. Based on the optimal solution, we obtain the out-of-sample performance using different sub-samples of different sizes. As in Section \ref{subsec:OS_performance}, we focus on the out-of-sample 95\% percentiles that may represent a risk-averse decision-maker. For brevity, since results for other bandwidths and settings are similar, we present the results for CSO with bandwidth $h=1.0$ and SAA under the setting $\nu=0.2$ and $R=0.5$.

\begin{figure}[t] 
\centering  
\includegraphics[scale=0.65]{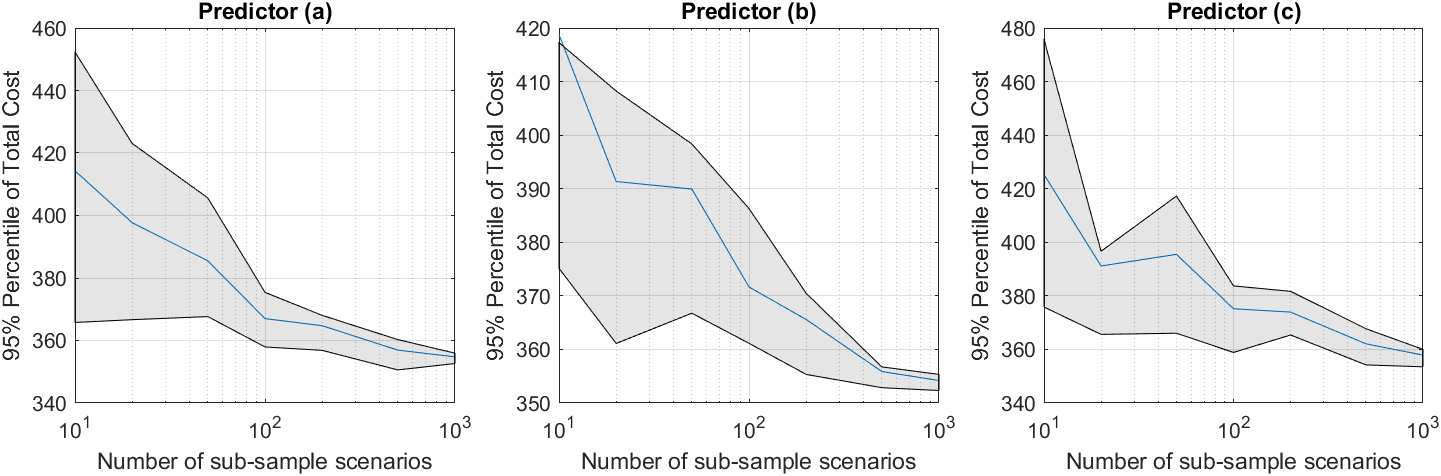}
\caption{Out-of-sample 95\% percentiles of the total cost over $20$ optimal schedules for different sub-sampling size $N'$ for CSO when $\nu=0.2$ and $R=0.5$. The solid line shows the mean over $20$ replications and the shaded region depicts the lower and upper quantiles.} \label{fig:subsample_size_CSO}
\end{figure}

\begin{figure}[t] 
\centering  
\includegraphics[scale=0.65]{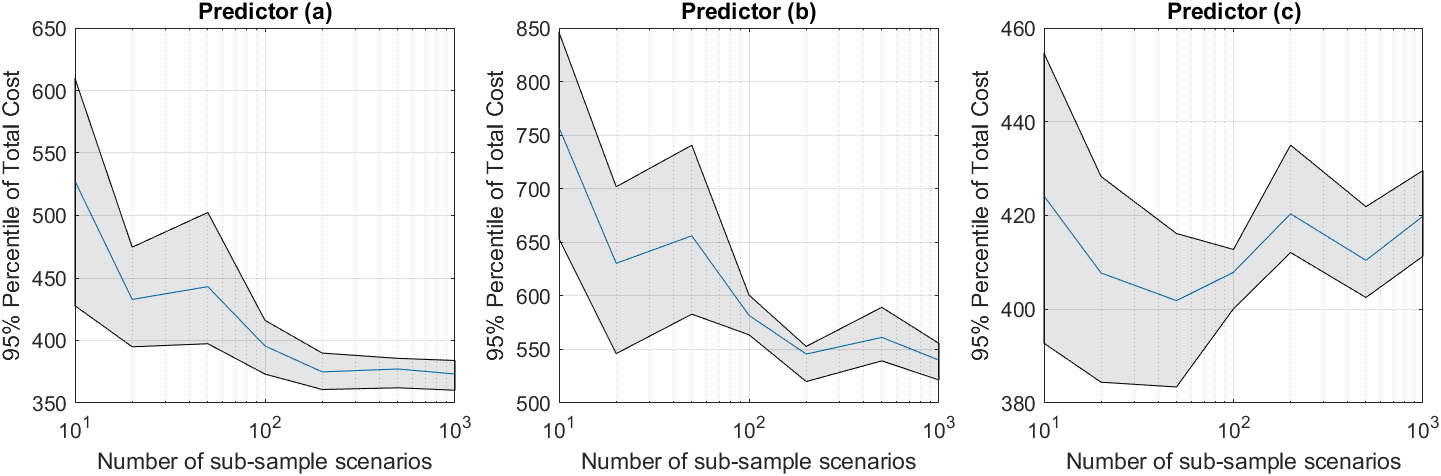}
\caption{Out-of-sample 95\% percentiles of the total cost over $20$ optimal schedules for different sub-sampling size $N'$ for SAA when $\nu=0.2$ and $R=0.5$. The solid line shows the mean over $20$ replications and the shaded region depicts the lower and upper quantiles.} \label{fig:subsample_size_SAA}
\end{figure}

Figures \ref{fig:subsample_size_CSO}--\ref{fig:subsample_size_SAA} show the corresponding out-of-sample 95\% percentiles of total cost for CSO and SAA, respectively. The solid line gives the mean over $20$ replications and the shaded region depicts the lower and upper quantiles. From these figures, we observe that the average out-of-sample performance is improving, mainly because of the fact that  more data points are included in the optimization model. Moreover, we observe that such an improvement would stabilized when $N'$ is large, say when $N'\in\{500,1000\}$. In particular, the shaded region, which represents the variation among different random sub-samples, becomes narrower when $N'$ reaches $500$ and $1,000$. These results demonstrate that using a random sub-sample of size $1,000$ would be sufficient to derive managerial insights for our experiments. 

\subsection{Effect of Bandwidth} \label{appdx:expt_bandwidth}

In this section, we examine the effect of bandwidth $h$ on the out-of-sample performance. In this experiment, we follow the experiment description in Section \ref{subsec:opt_schedule} to generate a set of $N=10,000$ pairs historical data, where we generate a sub-sample of size $N'=1,000$ to obtain an optimal solution by solving CSO with bandwidth $h$. As we have shown in \ref{appdx:expt_subsample}, performances of different random sub-samples when $N'=1,000$ are relatively stable. Next, we repeat this procedure to generate $20$ different optimal solutions for each bandwidth $h$. To examine the effect of the choice of bandwidth, we consider $h\in\{0.2,0.4,\dots,2.2\}$. Moreover, we solve CSO with both quantile and expectation objectives (since they may give different results). For CSO with a quantile objective, we focus on the out-of-sample 95\% percentile of the total cost while for CSO with an expectation objective, we focus on the out-of-sample average total cost. Note that such a choice is reasonable as the out-of-sample metric matches with the objective about which the decision-maker is concerned. For brevity, since results for other settings are similar, we present the results under the setting $\nu=0.2$ and $R=0.5$.

\begin{figure}[t] 
\centering  
\includegraphics[scale=0.65]{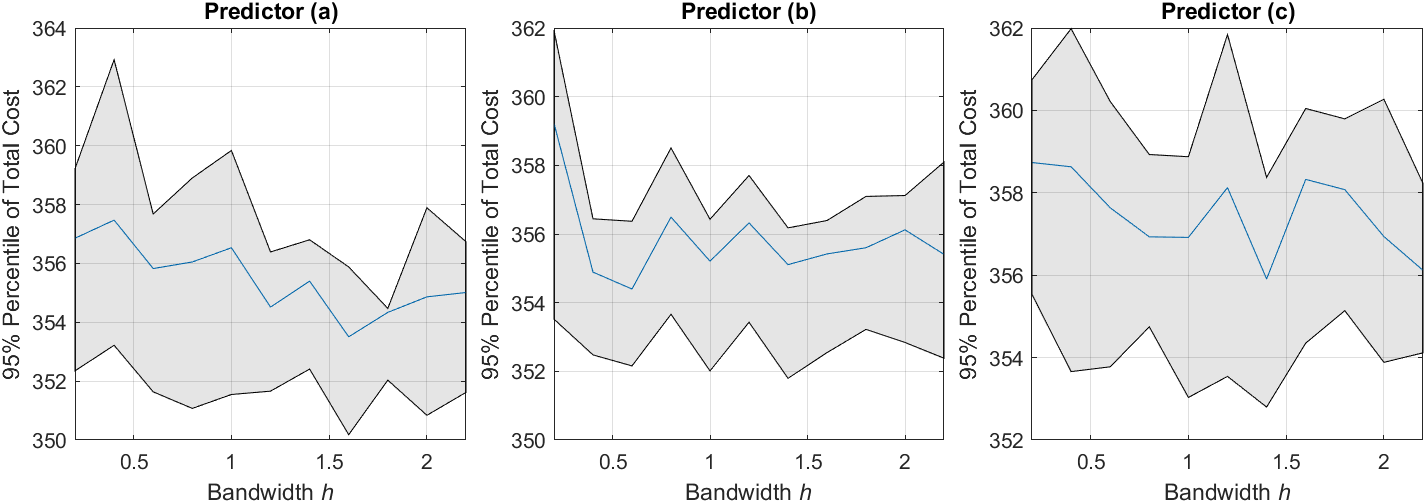}
\caption{Out-of-sample 95\% percentiles of the total cost over $20$ optimal schedules for different bandwidth $h$ for CSO under a quantile objective when $\nu=0.2$ and $R=0.5$. The solid line shows the mean over $20$ replications and the shaded region depicts the lower and upper quantiles.} \label{fig:bandwidth_CSO_quantile}
\end{figure}

\begin{figure}[t] 
\centering  
\includegraphics[scale=0.65]{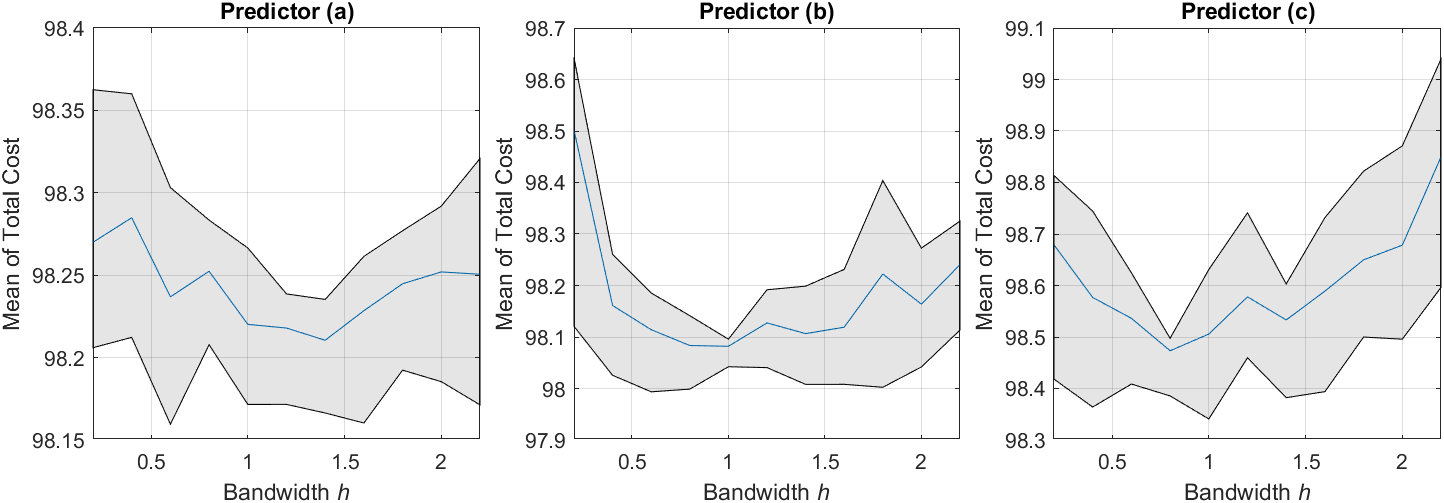}
\caption{Out-of-sample average total cost over $20$ optimal schedules for different bandwidth $h$ for CSO under an expectation objective when $\nu=0.2$ and $R=0.5$. The solid line shows the mean over $20$ replications and the shaded region depicts the lower and upper quantiles.} \label{fig:bandwidth_CSO_expectation}
\end{figure}

Figures \ref{fig:bandwidth_CSO_quantile}--\ref{fig:bandwidth_CSO_expectation} show the corresponding out-of-sample performance for CSO with a quantile and an expectation objective respectively. The solid line gives the mean over $20$ replications and the shaded region depicts the lower and upper quantiles. First, from Figure \ref{fig:bandwidth_CSO_quantile} for CSO with a quantile objective, we observe that the out-of-sample 95\% percentile does not change significantly within the tested range of bandwidths. Indeed, the out-of-sample 95\% percentiles fluctuates around some values, and a specific trend could not be observed. On the other hand, from Figure \ref{fig:bandwidth_CSO_expectation} for CSO with an expectation objective, we notice that the average out-of-sample total cost exhibits a slight decreasing trend at the beginning, followed by a slight increasing trend. Despite the changes in actual magnitude are relatively small (in the order of $0.1$), a choice of value $h=1.0$ seems to yield a better out-of-sample performance in general. Therefore, in the experiments, we choose the bandwidth for CSO models as $h=1.0$. We remark that, in practice, decision-makers could adopt methods such as cross-validation to calibrate the choice of bandwidth. 

\subsection{Additional Results on Optimal Schedules} \label{appdx:expt_opt_schedule}

\begin{figure}[t] 
\centering    
\includegraphics[scale=0.62]{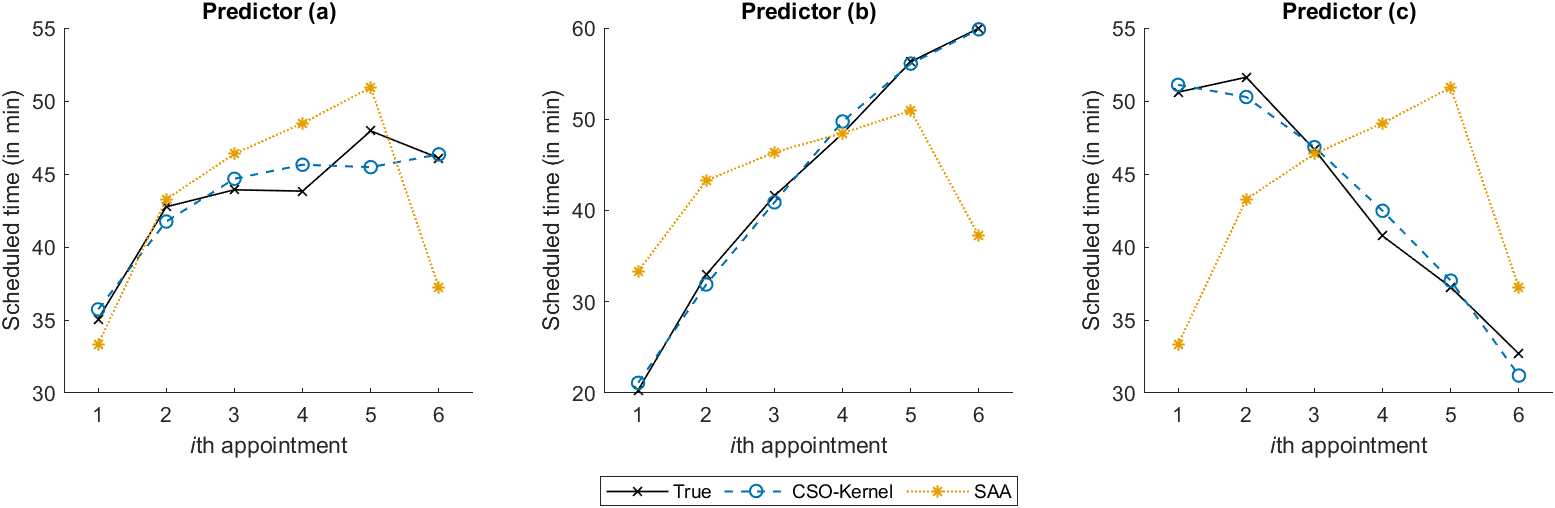}
\caption{Optimal appointment schedules from three different models when $\nu=0.2$ and $R=1.0$ under a quantile objective.} \label{fig:opt_schedule_2_quantile}
\end{figure}

\begin{figure}[t] 
\centering  
\includegraphics[scale=0.62]{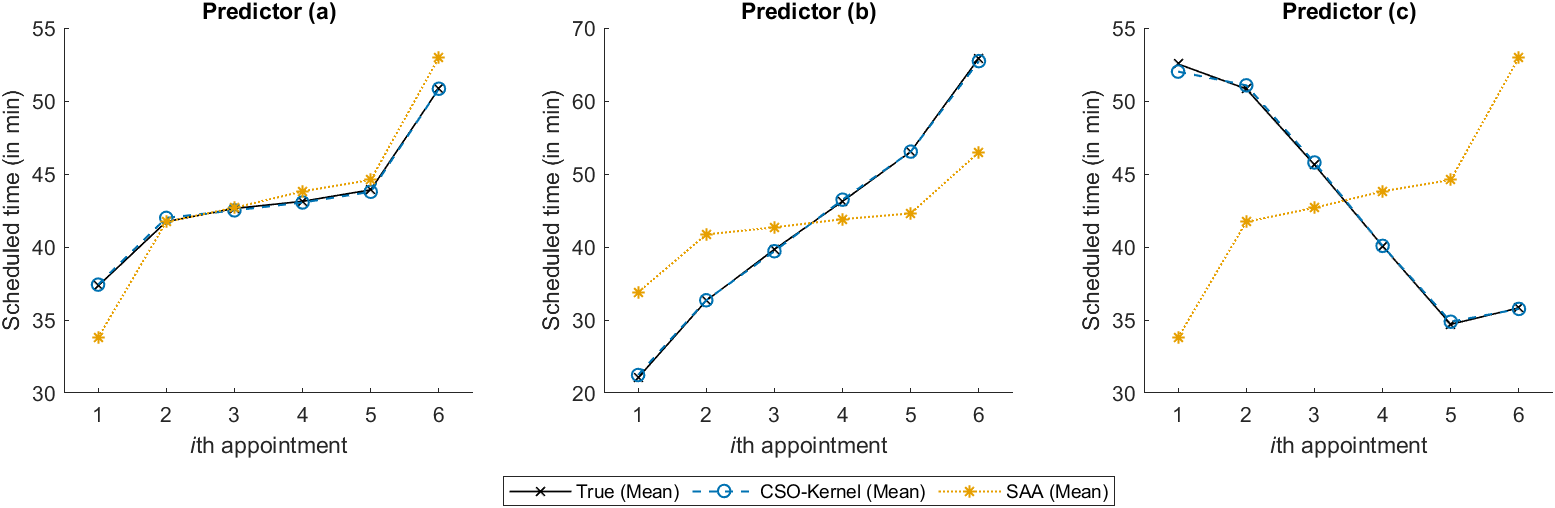}
\caption{Optimal appointment schedules from three different models when $\nu=0.2$ and $R=1.0$ under an expectation objective.} \label{fig:opt_schedule_2_mean}
\end{figure}

In this section, we present results on optimal schedules under the remaining settings that are not covered in Section \ref{subsec:opt_schedule}. First, in Figures \ref{fig:opt_schedule_2_quantile}--\ref{fig:opt_schedule_2_mean}, we plot the optimal schedules under quantile and expectation objectives respectively when $\nu=0.2$ and $R=1.0$. Similar to the observations in Section~\ref{subsec:opt_schedule}, the optimal schedule from CSO is close to the one obtained from simulation data under true distribution. While CSO is able to capture the predictor information, SAA could not and gives an optimal solution significantly different from CSO. Again, using a quantile objective provides a more risk-averse schedule that allocates more time to each appointment than using an expectation objective.

\begin{figure}[t] 
\centering  
\includegraphics[scale=0.62]{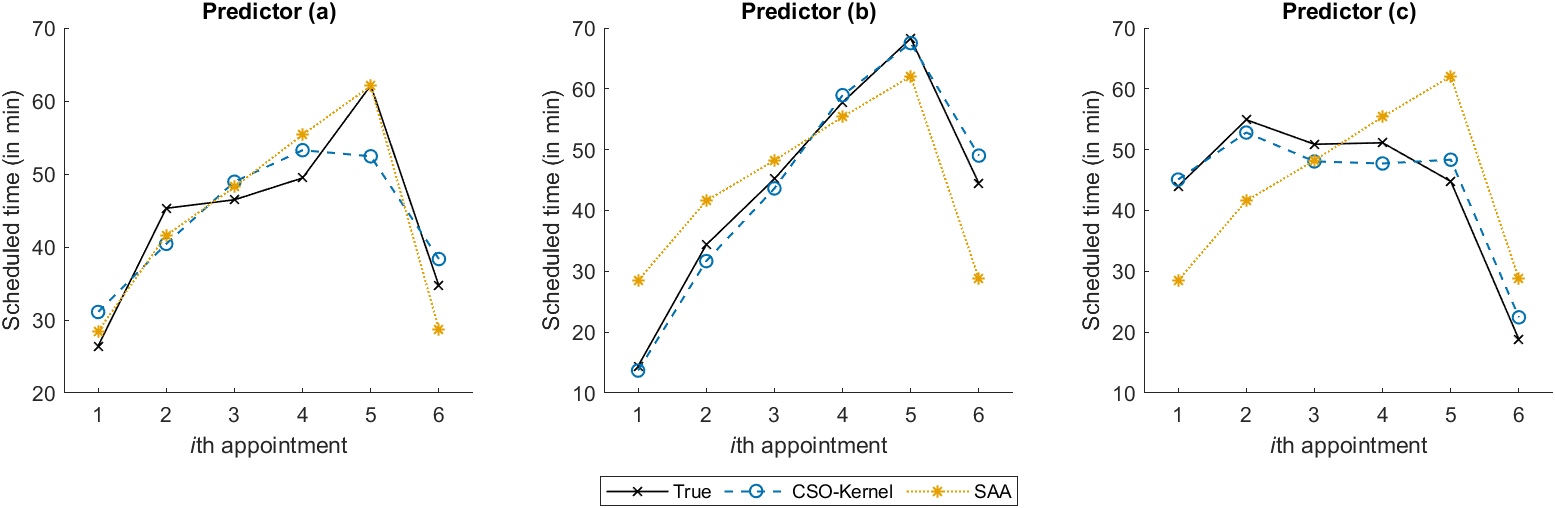}
\caption{Optimal appointment schedules from three different models when $\nu=0.5$ and $R=0.5$ under a quantile objective.} \label{fig:opt_schedule_3_quantile}
\end{figure}

\begin{figure}[t] 
\centering  
\includegraphics[scale=0.62]{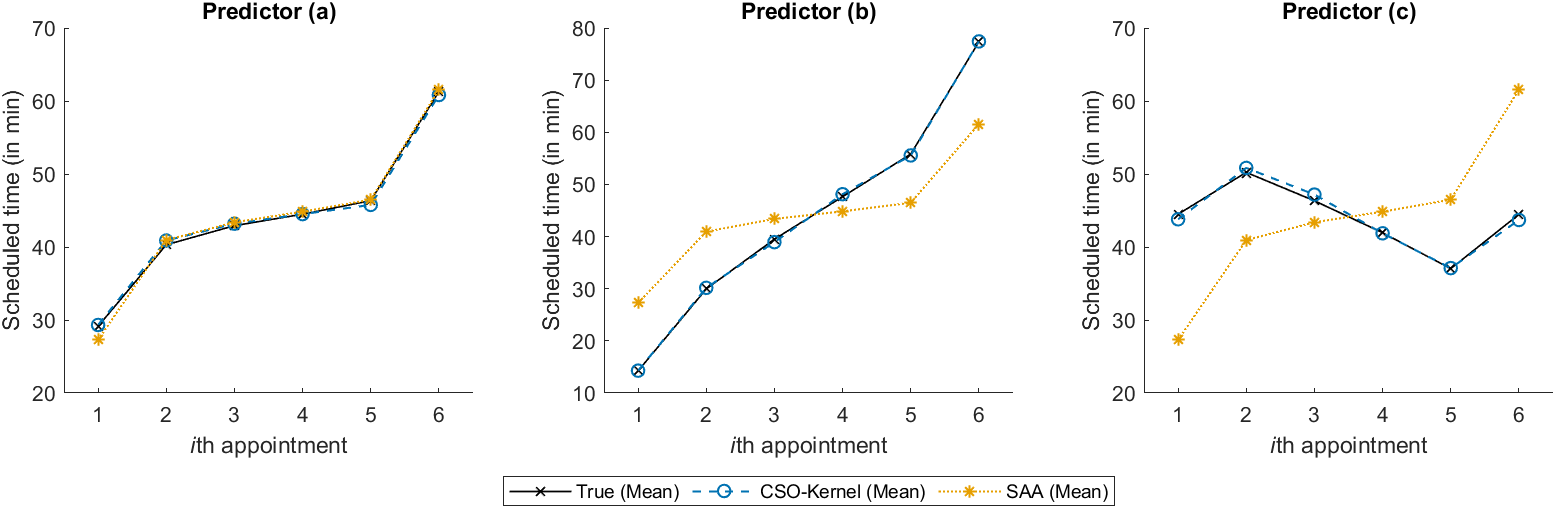}
\caption{Optimal appointment schedules from three different models when $\nu=0.5$ and $R=0.5$ under an expectation objective.} \label{fig:opt_schedule_3_mean}
\end{figure}

\begin{figure}[t] 
\centering  
\includegraphics[scale=0.62]{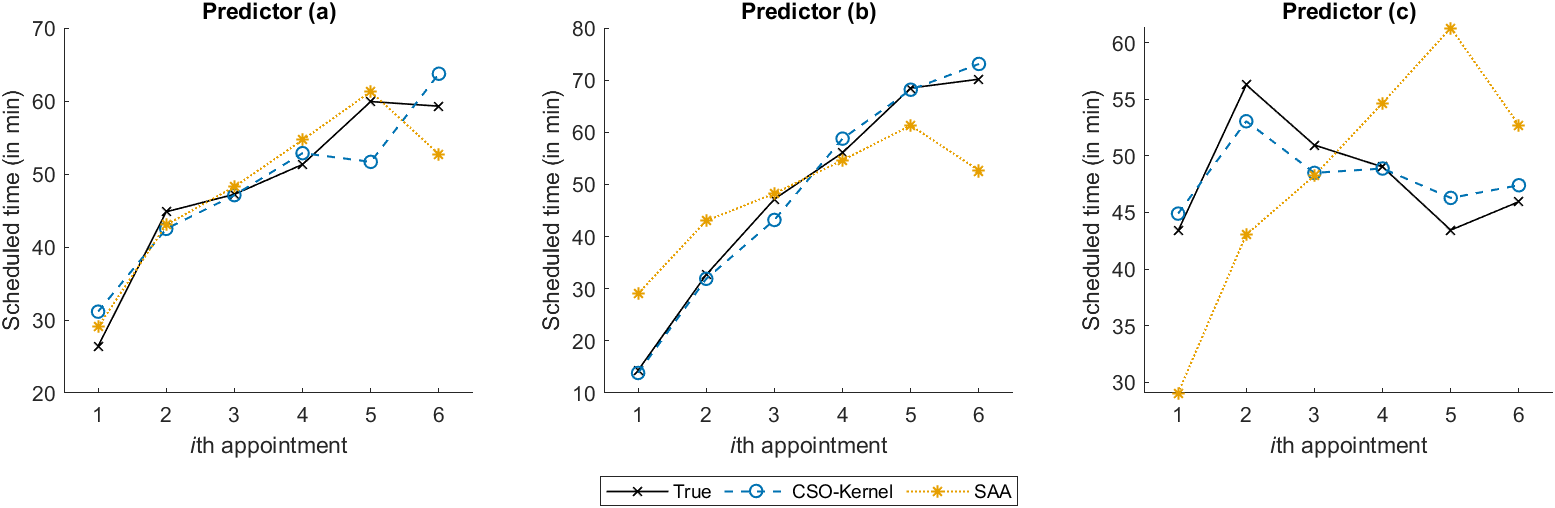}
\caption{Optimal appointment schedules from three different models when $\nu=0.5$ and $R=1.0$ under a quantile objective.} \label{fig:opt_schedule_4_quantile}
\end{figure}

\begin{figure}[t] 
\centering  
\includegraphics[scale=0.62]{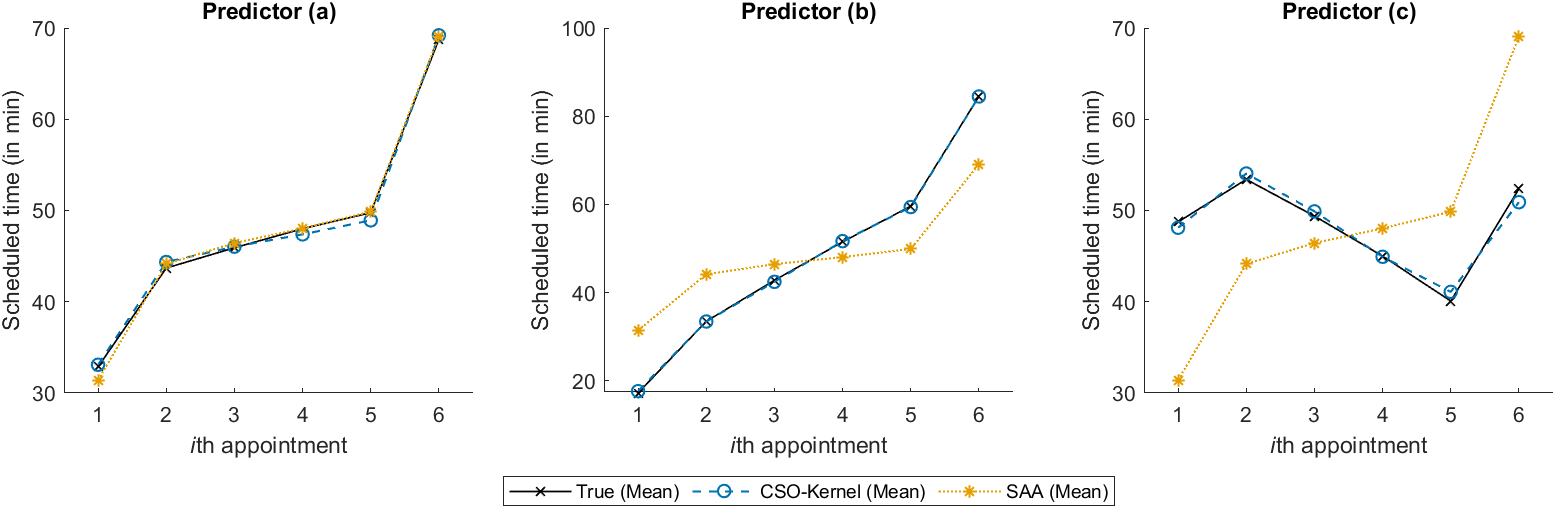}
\caption{Optimal appointment schedules from three different models when $\nu=0.5$ and $R=1.0$ under an expectation objective.} \label{fig:opt_schedule_4_mean}
\end{figure}

Next, Figures \ref{fig:opt_schedule_3_quantile}--\ref{fig:opt_schedule_3_mean} show the results on the optimal schedules when $\nu=0.5$ and $R=0.5$. Note that in this setting, the individual service duration variance is relatively larger. As a result, both optimal schedules from CSO and SAA are similar. However, when the underlying distribution exhibits an increasing or decreasing pattern as in predictors (b) or (c), the optimal schedules from CSO and SAA become different. Nevertheless, the time allocated to each appointment using a quantile objective is still longer than the one using an expectation objective. We have similar observations for the results when $\nu=0.5$ and $R=1.0$ in Figures \ref{fig:opt_schedule_4_quantile}--\ref{fig:opt_schedule_4_mean}.

\subsection{Additional Results on Out-of-Sample Performances} \label{appdx:expt_OS}

\begin{figure}[t] 
\centering  
\includegraphics[scale=0.65]{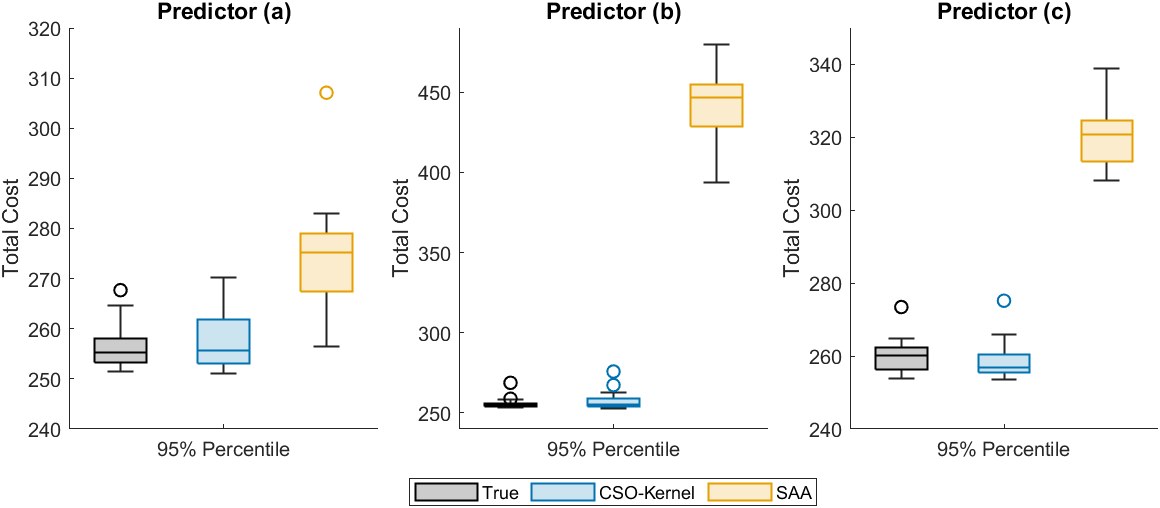}
\caption{Out-of-sample 95\% percentiles of the total cost over $20$ optimal schedules under a quantile objective from three different models when $\nu=0.2$ and $R=1.0$.} \label{fig:OS_95Per_2}
\end{figure}

\begin{figure}[t] 
\centering  
\includegraphics[scale=0.65]{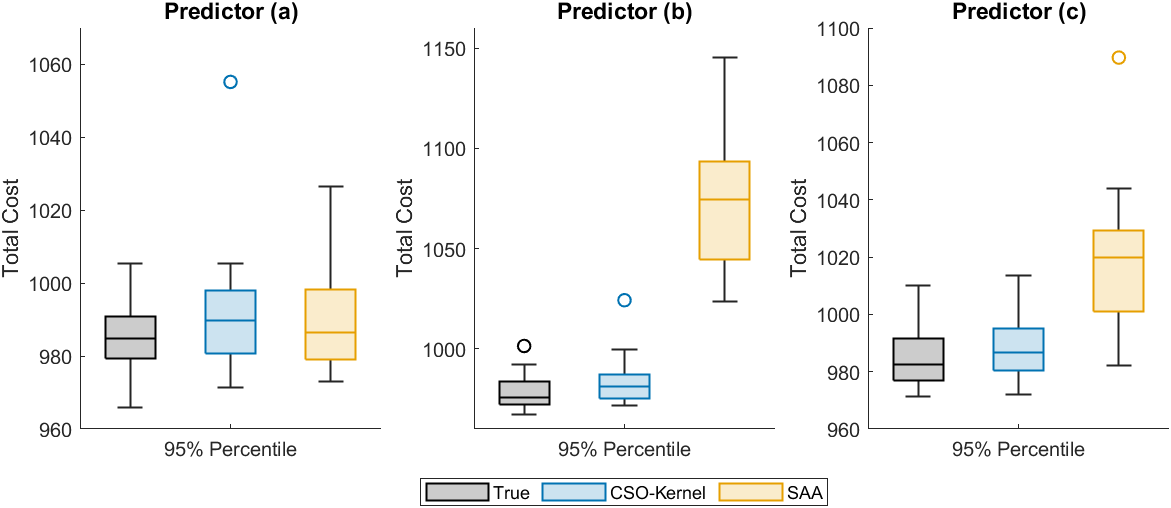}
\caption{Out-of-sample 95\% percentiles of the total cost over $20$ optimal schedules under a quantile objective from three different models when $\nu=0.5$ and $R=0.5$.} \label{fig:OS_95Per_3}
\end{figure}

\begin{figure}[t] 
\centering  
\includegraphics[scale=0.65]{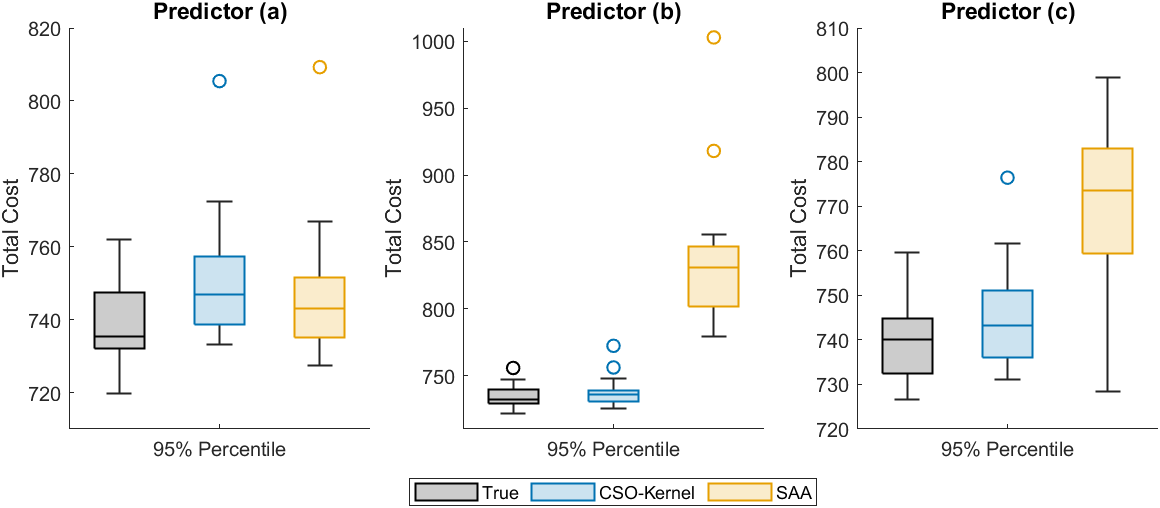}
\caption{Out-of-sample 95\% percentiles of the total cost over $20$ optimal schedules under a quantile objective from three different models when $\nu=0.5$ and $R=1.0$.} \label{fig:OS_95Per_4}
\end{figure}

In this section, we provide additional results on out-of-sample performances. First, we present the figures for out-of-sample 95\% percentiles of the total cost discussed in the beginning of Section \ref{subsec:OS_performance}. Figures \ref{fig:OS_95Per_2}--\ref{fig:OS_95Per_4} show the corresponding plots for the remaining settings: $(\nu,R)\in\{(0.2,1.0),(0.5,0.5),(0.5,1.0)\}$. As discussed in Section \ref{subsec:OS_performance}, we observe that CSO could perform better than SAA in general and such an improvement could be significant under predictors (b) and (c). Moreover, we note that when $\nu$ is relatively large, e.g., $\nu=0.5$, the performance between CSO and SAA could be similar under predictor (a). This could be explained by the similarity of the optimal schedules discussed in \ref{appdx:expt_opt_schedule}.

\begin{figure}[t] 
\centering  
\includegraphics[scale=0.65]{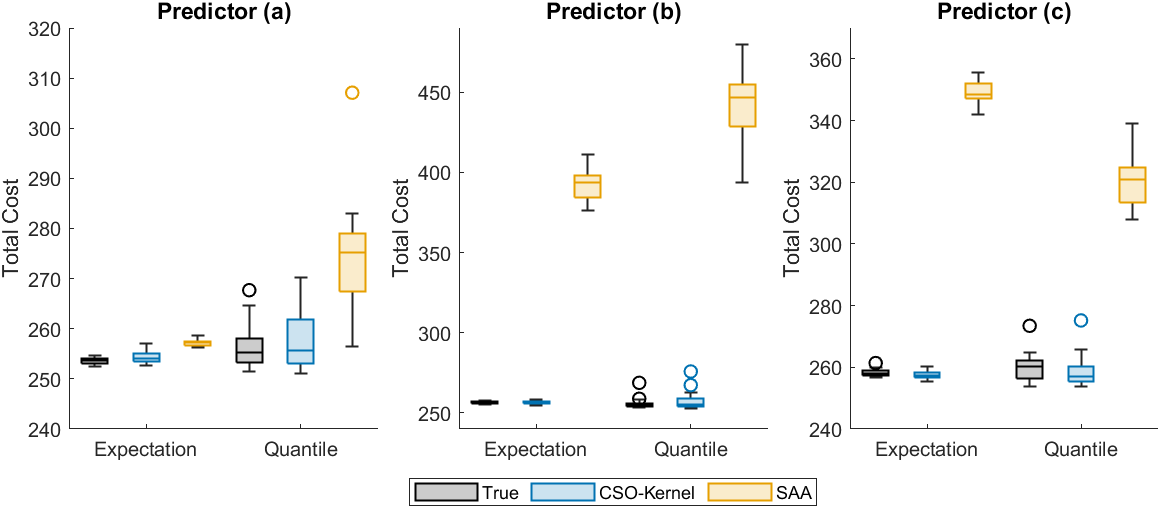}
\caption{Out-of-sample 95\% percentiles of the total cost over $20$ optimal schedules under an expectation or quantile objective from three different models when $\nu=0.2$ and $R=1.0$.} \label{fig:OS_with_mean_95Per_2}
\end{figure}

\begin{figure}[t] 
\centering  
\includegraphics[scale=0.65]{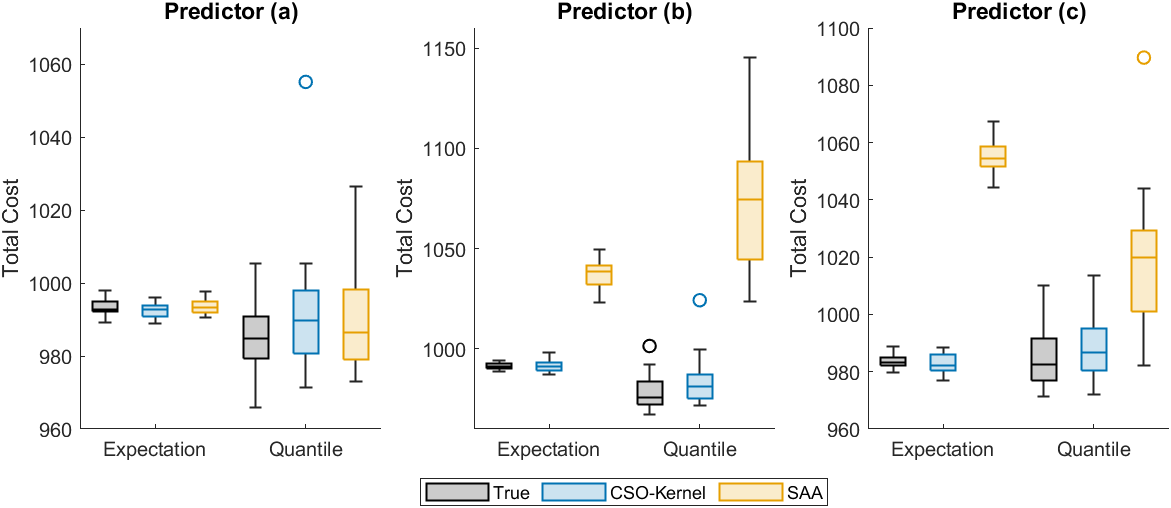}
\caption{Out-of-sample 95\% percentiles of the total cost over $20$ optimal schedules under an expectation or quantile objective from three different models when $\nu=0.5$ and $R=0.5$.} \label{fig:OS_with_mean_95Per_3}
\end{figure}

\begin{figure}[t] 
\centering  
\includegraphics[scale=0.65]{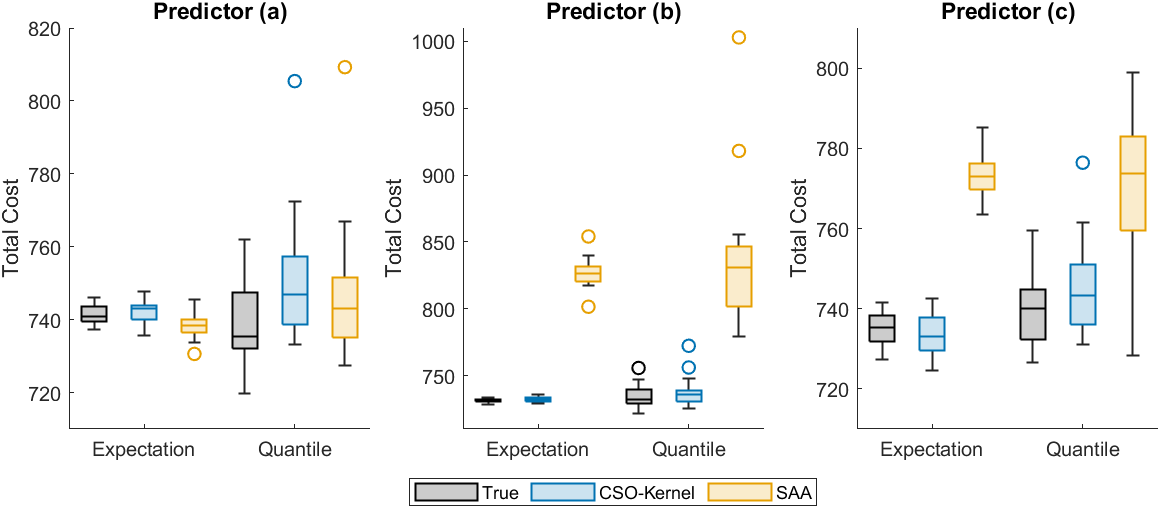}
\caption{Out-of-sample 95\% percentiles of the total cost over $20$ optimal schedules under an expecatation or quantile objective from three different models when $\nu=0.5$ and $R=1.0$.} \label{fig:OS_with_mean_95Per_4}
\end{figure}

Next, we present the figures for out-of-sample 95\% percentiles of the total cost for both quantile and expectation objectives discussed in the middle of Section \ref{subsec:OS_performance}.  Figures \ref{fig:OS_with_mean_95Per_2}--\ref{fig:OS_with_mean_95Per_4} show the corresponding plots for the remaining settings: $(\nu,R)\in\{(0.2,1.0),(0.5,0.5),(0.5,1.0)\}$. First, in Figure \ref{fig:OS_with_mean_95Per_2}, we observe that the performance using a quantile objective performs similar to the one using an expectation objective, with occasional improvements. This is reasonable since $R$ becomes larger, and thus the length of the day is longer. As a result, both models have a larger extent of flexibility for scheduling, in the sense that both models could assign more times to each appointment to protect against possible huge waiting costs. This could be observed from the optimal schedules shown in \ref{appdx:expt_opt_schedule}. Finally, we note that observations from Figures \ref{fig:OS_with_mean_95Per_3} and \ref{fig:OS_with_mean_95Per_4} are similar to those from Figure \ref{fig:OS_with_mean_95Per_1} and \ref{fig:OS_with_mean_95Per_2}, i.e., when $\nu=0.2$. 

Next, we present additional results under mis-specified distributions as discussed at the end of Section \ref{subsec:OS_performance}. Tables \ref{table:OS_misspec_95Per_2}--\ref{table:OS_misspec_95Per_4} present the corresponding mean and standard deviation for the remaining settings: $(\nu,R)\in\{(0.2,1.0),(0.5,0.5),(0.5,1.0)\}$. First, when comparing CSO with SAA, we note that the percentage increase in 95\% percentiles from CSO to SAA is reduced when the distribution is perturbed. As explained in Section \ref{subsec:OS_performance}, this is because the predictor becomes less informative under such a circumstance. However, CSO could still perform better than SAA in many cases, especially under predictor (c). Second, when comparing CSO with the one with an expectation objective, we observe that the percentage increase in the 95\% percentiles is increased when the distribution is perturbed, especially under Set II. This demonstrates the advantages of a risk-averse schedule over a risk-neutral one, in the sense that the risk-neutral schedule is less robust against possible mis-specified distributions.
\begin{table}[t]\centering  \small
\ra{1.0}  
\caption{Mean and standard deviation of the percentage change in the 95\% percentile of total cost when changing the implementation of the CSO solution based on a quantile objective to the SAA solution or the CSO solution based on an expectation objective over $20$ optimal solutions under $\nu=0.2$ and $R=1.0$. (The `Original' column corresponds to the out-of-sample scenarios generated under the true distribution as for the in-sample scenarios).} \label{table:OS_misspec_95Per_2}
\begin{tabular}{l||rr|rr|rr} \Xhline{1.0pt}
Predictor (a)       & \multicolumn{2}{c|}{Original} & \multicolumn{2}{c|}{Set I} & \multicolumn{2}{c}{Set II} \\ 
                    & Mean          & Std           & Mean         & Std         & Mean         & Std         \\ \hline
CSO v.s. SAA        & 6.4\%         & 4.2\%         & -0.4\%       & 1.4\%       & -1.3\%       & 2.1\%       \\
CSO v.s. CSO (Mean) & -1.3\%        & 2.1\%         & 0.7\%        & 1.0\%       & 0.5\%        & 1.5\%       \\ \Xhline{1.0pt}
Predictor (b)       & \multicolumn{2}{c|}{Original} & \multicolumn{2}{c|}{Set I} & \multicolumn{2}{c}{Set II} \\
                    & Mean          & Std           & Mean         & Std         & Mean         & Std         \\ \hline
CSO v.s. SAA        & 71.0\%        & 11.1\%        & 25.9\%       & 4.2\%       & -1.2\%       & 2.8\%       \\
CSO v.s. CSO (Mean) & -0.3\%        & 2.1\%         & 1.1\%        & 1.0\%       & 0.5\%        & 1.1\%       \\ \Xhline{1.0pt}
Predictor (c)       & \multicolumn{2}{c|}{Original} & \multicolumn{2}{c|}{Set I} & \multicolumn{2}{c}{Set II} \\
                    & Mean          & Std           & Mean         & Std         & Mean         & Std         \\ \hline
CSO v.s. SAA        & 23.9\%        & 4.6\%         & 8.2\%        & 2.0\%       & 11.0\%       & 1.9\%       \\
CSO v.s. CSO (Mean) & -0.5\%        & 1.9\%         & 0.6\%        & 0.7\%       & 0.4\%        & 1.3\%       \\ \hline
\Xhline{1.0pt}
\end{tabular}
\end{table}
\begin{table}[t]\centering  \small
\ra{1.0}  
\caption{Mean and standard deviation of the percentage change in the 95\% percentile of total cost when changing the implementation of the CSO solution based on a quantile objective to the SAA solution or the CSO solution based on an expectation objective over $20$ optimal solutions under $\nu=0.5$ and $R=0.5$. (The `Original' column corresponds to the out-of-sample scenarios generated under the true distribution as for the in-sample scenarios).} \label{table:OS_misspec_95Per_3}
\begin{tabular}{l||rr|rr|rr} \Xhline{1.0pt}
Predictor (a)       & \multicolumn{2}{c|}{Original} & \multicolumn{2}{c|}{Set I} & \multicolumn{2}{c}{Set II} \\ \
                    & Mean           & Std          & Mean         & Std         & Mean         & Std         \\ \hline
CSO v.s. SAA        & -0.1\%         & 2.2\%        & -0.4\%       & 1.2\%       & -0.4\%       & 2.4\%       \\
CSO v.s. CSO (Mean) & 0.1\%          & 1.7\%        & 0.7\%        & 0.9\%       & 2.5\%        & 1.9\%       \\ \Xhline{1.0pt}
Predictor (b)       & \multicolumn{2}{c|}{Original} & \multicolumn{2}{c|}{Set I} & \multicolumn{2}{c}{Set II} \\
                    & Mean           & Std          & Mean         & Std         & Mean         & Std         \\ \hline
CSO v.s. SAA        & 8.9\%          & 4.1\%        & 4.6\%        & 1.9\%       & -0.1\%       & 1.7\%       \\
CSO v.s. CSO (Mean) & 0.8\%          & 1.3\%        & 0.9\%        & 0.6\%       & 2.9\%        & 1.5\%       \\ \Xhline{1.0pt}
Predictor (c)       & \multicolumn{2}{c|}{Original} & \multicolumn{2}{c|}{Set I} & \multicolumn{2}{c}{Set II} \\
                    & Mean           & Std          & Mean         & Std         & Mean         & Std         \\ \hline
CSO v.s. SAA        & 3.3\%          & 2.6\%        & 1.2\%        & 1.3\%       & 5.8\%        & 3.3\%       \\
CSO v.s. CSO (Mean) & -0.5\%         & 1.0\%        & 0.5\%        & 0.7\%       & 2.0\%        & 2.0\%       \\ 
\Xhline{1.0pt}
\end{tabular}
\end{table}
\begin{table}[t]\centering  \small
\ra{1.0}  
\caption{Mean and standard deviation of the percentage change in the 95\% percentile of total cost when changing the implementation of the CSO solution based on a quantile objective to the SAA solution or the CSO solution based on an expectation objective over $20$ optimal solutions under $\nu=0.5$ and $R=1.0$. (The `Original' column corresponds to the out-of-sample scenarios generated under the true distribution as for the in-sample scenarios).} \label{table:OS_misspec_95Per_4}
\begin{tabular}{l||rr|rr|rr} \Xhline{1.0pt}
Predictor (a)       & \multicolumn{2}{c|}{Original} & \multicolumn{2}{c|}{Set I} & \multicolumn{2}{c}{Set II} \\ 
                    & Mean           & Std          & Mean         & Std         & Mean         & Std         \\ \hline
CSO v.s. SAA        & -0.4\%         & 3.6\%        & -0.3\%       & 1.5\%       & -0.7\%       & 2.5\%       \\
CSO v.s. CSO (Mean) & -1.0\%         & 2.3\%        & 0.0\%        & 1.0\%       & -0.2\%       & 2.2\%       \\ \Xhline{1.0pt}
Predictor (b)       & \multicolumn{2}{c|}{Original} & \multicolumn{2}{c|}{Set I} & \multicolumn{2}{c}{Set II} \\
                    & Mean           & Std          & Mean         & Std         & Mean         & Std         \\ \hline
CSO v.s. SAA        & 13.4\%         & 7.2\%        & 6.2\%        & 3.3\%       & 0.4\%        & 2.8\%       \\
CSO v.s. CSO (Mean) & -0.7\%         & 1.4\%        & 0.0\%        & 0.8\%       & 0.3\%        & 1.7\%       \\ \Xhline{1.0pt}
Predictor (c)       & \multicolumn{2}{c|}{Original} & \multicolumn{2}{c|}{Set I} & \multicolumn{2}{c}{Set II} \\
                    & Mean           & Std          & Mean         & Std         & Mean         & Std         \\ \hline
CSO v.s. SAA        & 3.5\%          & 3.4\%        & 1.3\%        & 1.0\%       & 6.5\%        & 3.3\%       \\
CSO v.s. CSO (Mean) & -1.6\%         & 1.5\%        & -0.7\%       & 0.6\%       & -0.8\%       & 2.6\%       \\ 
\Xhline{1.0pt}
\end{tabular}
\end{table}

\clearpage
\bibliographystyle{elsarticle-harv}
\bibliography{references}

\begin{thebibliography}{78}
\expandafter\ifx\csname natexlab\endcsname\relax\def\natexlab#1{#1}\fi
\providecommand{\url}[1]{\texttt{#1}}
\providecommand{\href}[2]{#2}
\providecommand{\path}[1]{#1}
\providecommand{\DOIprefix}{doi:}
\providecommand{\ArXivprefix}{arXiv:}
\providecommand{\URLprefix}{URL: }
\providecommand{\Pubmedprefix}{pmid:}
\providecommand{\doi}[1]{\href{http://dx.doi.org/#1}{\path{#1}}}
\providecommand{\Pubmed}[1]{\href{pmid:#1}{\path{#1}}}
\providecommand{\bibinfo}[2]{#2}
\ifx\xfnm\relax \def\xfnm[#1]{\unskip,\space#1}\fi
\bibitem[{Babat et~al.(2018)Babat, Vera and Zuluaga}]{Babat_et_al:2018}
\bibinfo{author}{Babat, O.}, \bibinfo{author}{Vera, J.C.},
  \bibinfo{author}{Zuluaga, L.F.}, \bibinfo{year}{2018}.
\newblock \bibinfo{title}{Computing near-optimal value-at-risk portfolios using
  integer programming techniques}.
\newblock \bibinfo{journal}{European Journal of Operational Research}
  \bibinfo{volume}{266}, \bibinfo{pages}{304--315}.
\bibitem[{Bai et~al.(2021)Bai, Sun and Zheng}]{Bai_et_al:2021}
\bibinfo{author}{Bai, X.}, \bibinfo{author}{Sun, J.}, \bibinfo{author}{Zheng,
  X.}, \bibinfo{year}{2021}.
\newblock \bibinfo{title}{An augmented {L}agrangian decomposition method for
  chance-constrained optimization problems}.
\newblock \bibinfo{journal}{INFORMS Journal on Computing} \bibinfo{volume}{33},
  \bibinfo{pages}{1056--1069}.
\bibitem[{Ban and Rudin(2019)}]{Ban_Rudin:2019}
\bibinfo{author}{Ban, G.Y.}, \bibinfo{author}{Rudin, C.}, \bibinfo{year}{2019}.
\newblock \bibinfo{title}{The big data newsvendor: Practical insights from
  machine learning}.
\newblock \bibinfo{journal}{Operations Research} \bibinfo{volume}{67},
  \bibinfo{pages}{90--108}.
\bibitem[{Benati and Rizzi(2007)}]{Benati_Rizzi:2007}
\bibinfo{author}{Benati, S.}, \bibinfo{author}{Rizzi, R.},
  \bibinfo{year}{2007}.
\newblock \bibinfo{title}{A mixed integer linear programming formulation of the
  optimal mean/value-at-risk portfolio problem}.
\newblock \bibinfo{journal}{European Journal of Operational Research}
  \bibinfo{volume}{176}, \bibinfo{pages}{423--434}.
\bibitem[{Bertsimas and Kallus(2020)}]{Bertsimas_Kallus:2020}
\bibinfo{author}{Bertsimas, D.}, \bibinfo{author}{Kallus, N.},
  \bibinfo{year}{2020}.
\newblock \bibinfo{title}{From predictive to prescriptive analytics}.
\newblock \bibinfo{journal}{Management Science} \bibinfo{volume}{66},
  \bibinfo{pages}{1025--1044}.
\bibitem[{Bertsimas and Koduri(2022)}]{Bertsimas_Koduri:2022}
\bibinfo{author}{Bertsimas, D.}, \bibinfo{author}{Koduri, N.},
  \bibinfo{year}{2022}.
\newblock \bibinfo{title}{Data-driven optimization: A reproducing kernel
  {Hilbert} space approach}.
\newblock \bibinfo{journal}{Operations Research} \bibinfo{volume}{70},
  \bibinfo{pages}{454--471}.
\bibitem[{Bertsimas and McCord(2019)}]{Bertsimas_McCord:2019}
\bibinfo{author}{Bertsimas, D.}, \bibinfo{author}{McCord, C.},
  \bibinfo{year}{2019}.
\newblock \bibinfo{title}{From predictions to prescriptions in multistage
  optimization problems}.
\newblock \bibinfo{journal}{arXiv preprint arXiv:1904.11637} .
\bibitem[{Bertsimas et~al.(2023)Bertsimas, McCord and
  Sturt}]{Bertsimas_et_al:2019}
\bibinfo{author}{Bertsimas, D.}, \bibinfo{author}{McCord, C.},
  \bibinfo{author}{Sturt, B.}, \bibinfo{year}{2023}.
\newblock \bibinfo{title}{Dynamic optimization with side information}.
\newblock \bibinfo{journal}{European Journal of Operational Research}
  \bibinfo{volume}{304}, \bibinfo{pages}{634--651}.
\bibitem[{Bertsimas and Van~Parys(2022)}]{Bertsimas_et_al:2021}
\bibinfo{author}{Bertsimas, D.}, \bibinfo{author}{Van~Parys, B.},
  \bibinfo{year}{2022}.
\newblock \bibinfo{title}{Bootstrap robust prescriptive analytics}.
\newblock \bibinfo{journal}{Mathematical Programming} \bibinfo{volume}{195},
  \bibinfo{pages}{39--78}.
\bibitem[{Beyerlein(2014)}]{Beyerlein:2014}
\bibinfo{author}{Beyerlein, A.}, \bibinfo{year}{2014}.
\newblock \bibinfo{title}{Quantile regression—opportunities and challenges
  from a user's perspective}.
\newblock \bibinfo{journal}{American Journal of Epidemiology}
  \bibinfo{volume}{180}, \bibinfo{pages}{330--331}.
\bibitem[{Bouscary et~al.(2025)Bouscary, Zhang and Amin}]{Bouscary_et_al:2025}
\bibinfo{author}{Bouscary, M.}, \bibinfo{author}{Zhang, J.},
  \bibinfo{author}{Amin, S.}, \bibinfo{year}{2025}.
\newblock \bibinfo{title}{Reducing contextual stochastic bilevel optimization
  via structured function approximation}.
\newblock \bibinfo{journal}{arXiv preprint arXiv:2503.19991} .
\bibitem[{Cattaruzza et~al.(2024)Cattaruzza, Labb{\'e}, Petris, Roland and
  Schmidt}]{Cattaruzza_et_al:2024}
\bibinfo{author}{Cattaruzza, D.}, \bibinfo{author}{Labb{\'e}, M.},
  \bibinfo{author}{Petris, M.}, \bibinfo{author}{Roland, M.},
  \bibinfo{author}{Schmidt, M.}, \bibinfo{year}{2024}.
\newblock \bibinfo{title}{Exact and heuristic solution techniques for
  mixed-integer quantile minimization problems}.
\newblock \bibinfo{journal}{INFORMS Journal on Computing} \bibinfo{volume}{36},
  \bibinfo{pages}{1084--1107}.
\bibitem[{{\c{C}}etinkaya and Thiele(2015)}]{Cetinkaya_Thiele:2015}
\bibinfo{author}{{\c{C}}etinkaya, E.}, \bibinfo{author}{Thiele, A.},
  \bibinfo{year}{2015}.
\newblock \bibinfo{title}{Data-driven portfolio management with quantile
  constraints}.
\newblock \bibinfo{journal}{OR Spectrum} \bibinfo{volume}{37},
  \bibinfo{pages}{761--786}.
\bibitem[{Chen et~al.(2006)Chen, Daskin, Shen and Uryasev}]{Chen_et_al:2006}
\bibinfo{author}{Chen, G.}, \bibinfo{author}{Daskin, M.S.},
  \bibinfo{author}{Shen, Z.J.M.}, \bibinfo{author}{Uryasev, S.},
  \bibinfo{year}{2006}.
\newblock \bibinfo{title}{The $\alpha$-reliable mean-excess regret model for
  stochastic facility location modeling}.
\newblock \bibinfo{journal}{Naval Research Logistics} \bibinfo{volume}{53},
  \bibinfo{pages}{617--626}.
\bibitem[{Chenreddy et~al.(2022)Chenreddy, Bandi and
  Delage}]{Chenreddy_et_al:2022}
\bibinfo{author}{Chenreddy, A.R.}, \bibinfo{author}{Bandi, N.},
  \bibinfo{author}{Delage, E.}, \bibinfo{year}{2022}.
\newblock \bibinfo{title}{Data-driven conditional robust optimization}.
\newblock \bibinfo{journal}{Advances in Neural Information Processing Systems}
  \bibinfo{volume}{35}, \bibinfo{pages}{9525--9537}.
\bibitem[{Costa and Iyengar(2023)}]{Costa_Iyengar:2023}
\bibinfo{author}{Costa, G.}, \bibinfo{author}{Iyengar, G.N.},
  \bibinfo{year}{2023}.
\newblock \bibinfo{title}{Distributionally robust end-to-end portfolio
  construction}.
\newblock \bibinfo{journal}{Quantitative Finance} \bibinfo{volume}{23},
  \bibinfo{pages}{1465--1482}.
\bibitem[{Deng and Sen(2022)}]{Deng_Sen:2022}
\bibinfo{author}{Deng, Y.}, \bibinfo{author}{Sen, S.}, \bibinfo{year}{2022}.
\newblock \bibinfo{title}{Predictive stochastic programming}.
\newblock \bibinfo{journal}{Computational Management Science}
  \bibinfo{volume}{19}, \bibinfo{pages}{65--98}.
\bibitem[{Donti et~al.(2017)Donti, Amos and Kolter}]{Donti_et_al:2017}
\bibinfo{author}{Donti, P.}, \bibinfo{author}{Amos, B.},
  \bibinfo{author}{Kolter, J.Z.}, \bibinfo{year}{2017}.
\newblock \bibinfo{title}{Task-based end-to-end model learning in stochastic
  optimization}.
\newblock \bibinfo{journal}{Advances in Neural Information Processing Systems}
  \bibinfo{volume}{30}.
\bibitem[{Duchi et~al.(2021)Duchi, Glynn and Namkoong}]{Duchi_et_al:2021}
\bibinfo{author}{Duchi, J.C.}, \bibinfo{author}{Glynn, P.W.},
  \bibinfo{author}{Namkoong, H.}, \bibinfo{year}{2021}.
\newblock \bibinfo{title}{Statistics of robust optimization: A generalized
  empirical likelihood approach}.
\newblock \bibinfo{journal}{Mathematics of Operations Research}
  \bibinfo{volume}{46}, \bibinfo{pages}{946--969}.
\bibitem[{Elmachtoub and Grigas(2022)}]{Elmachtoub_Grigas:2022}
\bibinfo{author}{Elmachtoub, A.N.}, \bibinfo{author}{Grigas, P.},
  \bibinfo{year}{2022}.
\newblock \bibinfo{title}{Smart “predict, then optimize”}.
\newblock \bibinfo{journal}{Management Science} \bibinfo{volume}{68},
  \bibinfo{pages}{9--26}.
\bibitem[{Feng et~al.(2015)Feng, W{\"a}chter and Staum}]{Feng_et_al:2015}
\bibinfo{author}{Feng, M.}, \bibinfo{author}{W{\"a}chter, A.},
  \bibinfo{author}{Staum, J.}, \bibinfo{year}{2015}.
\newblock \bibinfo{title}{Practical algorithms for value-at-risk portfolio
  optimization problems}.
\newblock \bibinfo{journal}{Quantitative Finance Letters} \bibinfo{volume}{3},
  \bibinfo{pages}{1--9}.
\bibitem[{Filippi et~al.(2020)Filippi, Guastaroba and
  Speranza}]{Filippi_et_al:2020}
\bibinfo{author}{Filippi, C.}, \bibinfo{author}{Guastaroba, G.},
  \bibinfo{author}{Speranza, M.G.}, \bibinfo{year}{2020}.
\newblock \bibinfo{title}{Conditional value-at-risk beyond finance: A survey}.
\newblock \bibinfo{journal}{International Transactions in Operational Research}
  \bibinfo{volume}{27}, \bibinfo{pages}{1277--1319}.
\bibitem[{Gaivoronski and Pflug(2005)}]{Gaivoronski_Pflug:2005}
\bibinfo{author}{Gaivoronski, A.A.}, \bibinfo{author}{Pflug, G.},
  \bibinfo{year}{2005}.
\newblock \bibinfo{title}{Value-at-risk in portfolio optimization: Properties
  and computational approach}.
\newblock \bibinfo{journal}{Journal of Risk} \bibinfo{volume}{7},
  \bibinfo{pages}{1--31}.
\bibitem[{Ge et~al.(2014)Ge, Wan, Wang and Zhang}]{Ge_et_al:2014}
\bibinfo{author}{Ge, D.}, \bibinfo{author}{Wan, G.}, \bibinfo{author}{Wang,
  Z.}, \bibinfo{author}{Zhang, J.}, \bibinfo{year}{2014}.
\newblock \bibinfo{title}{A note on appointment scheduling with piecewise
  linear cost functions}.
\newblock \bibinfo{journal}{Mathematics of Operations Research}
  \bibinfo{volume}{39}, \bibinfo{pages}{1244--1251}.
\bibitem[{Gneiting(2011)}]{Gneiting:2011}
\bibinfo{author}{Gneiting, T.}, \bibinfo{year}{2011}.
\newblock \bibinfo{title}{Making and evaluating point forecasts}.
\newblock \bibinfo{journal}{Journal of the American Statistical Association}
  \bibinfo{volume}{106}, \bibinfo{pages}{746--762}.
\bibitem[{Grass and Fischer(2016)}]{Grass_Fischer:2016}
\bibinfo{author}{Grass, E.}, \bibinfo{author}{Fischer, K.},
  \bibinfo{year}{2016}.
\newblock \bibinfo{title}{Two-stage stochastic programming in disaster
  management: A literature survey}.
\newblock \bibinfo{journal}{Surveys in Operations Research and Management
  Science} \bibinfo{volume}{21}, \bibinfo{pages}{85--100}.
\bibitem[{Grieco et~al.(2021)Grieco, Utley and Crowe}]{Grieco_et_al:2021}
\bibinfo{author}{Grieco, L.}, \bibinfo{author}{Utley, M.},
  \bibinfo{author}{Crowe, S.}, \bibinfo{year}{2021}.
\newblock \bibinfo{title}{Operational research applied to decisions in home
  health care: A systematic literature review}.
\newblock \bibinfo{journal}{Journal of the Operational Research Society}
  \bibinfo{volume}{72}, \bibinfo{pages}{1960--1991}.
\bibitem[{Henrion and R{\"o}misch(2004)}]{Henrion_Romisch:2004}
\bibinfo{author}{Henrion, R.}, \bibinfo{author}{R{\"o}misch, W.},
  \bibinfo{year}{2004}.
\newblock \bibinfo{title}{H{\"o}lder and {Lipschitz} stability of solution sets
  in programs with probabilistic constraints}.
\newblock \bibinfo{journal}{Mathematical Programming} \bibinfo{volume}{100},
  \bibinfo{pages}{589--611}.
\bibitem[{Hu et~al.(2023)Hu, Wang, Xie, Krause and Kuhn}]{Hu_et_al:2023}
\bibinfo{author}{Hu, Y.}, \bibinfo{author}{Wang, J.}, \bibinfo{author}{Xie,
  Y.}, \bibinfo{author}{Krause, A.}, \bibinfo{author}{Kuhn, D.},
  \bibinfo{year}{2023}.
\newblock \bibinfo{title}{Contextual stochastic bilevel optimization}.
\newblock \bibinfo{journal}{Advances in Neural Information Processing Systems}
  \bibinfo{volume}{36}, \bibinfo{pages}{78412--78434}.
\bibitem[{Ivanov and Naumov(2012)}]{Ivanov_Naumov:2012}
\bibinfo{author}{Ivanov, S.V.}, \bibinfo{author}{Naumov, A.},
  \bibinfo{year}{2012}.
\newblock \bibinfo{title}{Algorithm to optimize the quantile criterion for the
  polyhedral loss function and discrete distribution of random parameters}.
\newblock \bibinfo{journal}{Automation and Remote Control}
  \bibinfo{volume}{73}, \bibinfo{pages}{105--117}.
\bibitem[{Jiang et~al.(2019)Jiang, Ryu and Xu}]{Jiang_et_al:2019}
\bibinfo{author}{Jiang, R.}, \bibinfo{author}{Ryu, M.}, \bibinfo{author}{Xu,
  G.}, \bibinfo{year}{2019}.
\newblock \bibinfo{title}{Data-driven distributionally robust appointment
  scheduling over {Wasserstein} balls}.
\newblock \bibinfo{journal}{arXiv preprint arXiv:1907.03219} .
\bibitem[{Jiang et~al.(2017)Jiang, Shen and Zhang}]{Jiang_et_al:2017}
\bibinfo{author}{Jiang, R.}, \bibinfo{author}{Shen, S.},
  \bibinfo{author}{Zhang, Y.}, \bibinfo{year}{2017}.
\newblock \bibinfo{title}{Integer programming approaches for appointment
  scheduling with random no-shows and service durations}.
\newblock \bibinfo{journal}{Operations Research} \bibinfo{volume}{65},
  \bibinfo{pages}{1638--1656}.
\bibitem[{Jiang et~al.(2012)Jiang, Zhang, Li and Guan}]{Jiang_et_al:2012}
\bibinfo{author}{Jiang, R.}, \bibinfo{author}{Zhang, M.}, \bibinfo{author}{Li,
  G.}, \bibinfo{author}{Guan, Y.}, \bibinfo{year}{2012}.
\newblock \bibinfo{title}{Benders' decomposition for the two-stage security
  constrained robust unit commitment problem}, in: \bibinfo{booktitle}{IIE
  Annual Conference Proceedings}, \bibinfo{organization}{Institute of
  Industrial and Systems Engineers (IISE)}. p.~\bibinfo{pages}{1}.
\bibitem[{Jiang and Kou(2021)}]{Jiang_Kou:2021}
\bibinfo{author}{Jiang, W.}, \bibinfo{author}{Kou, S.}, \bibinfo{year}{2021}.
\newblock \bibinfo{title}{Simulating risk measures via asymptotic expansions
  for relative errors}.
\newblock \bibinfo{journal}{Mathematical Finance} \bibinfo{volume}{31},
  \bibinfo{pages}{907--942}.
\bibitem[{Kannan et~al.(2024)Kannan, Bayraksan and Luedtke}]{Kannan_et_al:2023}
\bibinfo{author}{Kannan, R.}, \bibinfo{author}{Bayraksan, G.},
  \bibinfo{author}{Luedtke, J.R.}, \bibinfo{year}{2024}.
\newblock \bibinfo{title}{Residuals-based distributionally robust optimization
  with covariate information}.
\newblock \bibinfo{journal}{Mathematical Programming} \bibinfo{volume}{207},
  \bibinfo{pages}{369--425}.
\bibitem[{Kannan et~al.(2025)Kannan, Bayraksan and Luedtke}]{Kannan_et_al:2025}
\bibinfo{author}{Kannan, R.}, \bibinfo{author}{Bayraksan, G.},
  \bibinfo{author}{Luedtke, J.R.}, \bibinfo{year}{2025}.
\newblock \bibinfo{title}{Data-driven sample average approximation with
  covariate information}.
\newblock \bibinfo{journal}{Operations Research} .
\bibitem[{Koenker and Hallock(2001)}]{Koenker_Hallock:2001}
\bibinfo{author}{Koenker, R.}, \bibinfo{author}{Hallock, K.F.},
  \bibinfo{year}{2001}.
\newblock \bibinfo{title}{Quantile regression}.
\newblock \bibinfo{journal}{Journal of Economic Perspectives}
  \bibinfo{volume}{15}, \bibinfo{pages}{143--156}.
\bibitem[{Laksaci et~al.(2009)Laksaci, Lemdani and
  Ould-Sa{\"\i}d}]{Laksaci_et_al:2009}
\bibinfo{author}{Laksaci, A.}, \bibinfo{author}{Lemdani, M.},
  \bibinfo{author}{Ould-Sa{\"\i}d, E.}, \bibinfo{year}{2009}.
\newblock \bibinfo{title}{A generalized {L1}-approach for a kernel estimator of
  conditional quantile with functional regressors: Consistency and asymptotic
  normality}.
\newblock \bibinfo{journal}{Statistics \& Probability Letters}
  \bibinfo{volume}{79}, \bibinfo{pages}{1065--1073}.
\bibitem[{Le et~al.(2021)Le, Nguyen, Yamada, Blanchet and
  Nguyen}]{Le_et_al:2021}
\bibinfo{author}{Le, T.}, \bibinfo{author}{Nguyen, T.},
  \bibinfo{author}{Yamada, M.}, \bibinfo{author}{Blanchet, J.},
  \bibinfo{author}{Nguyen, V.A.}, \bibinfo{year}{2021}.
\newblock \bibinfo{title}{Adversarial regression with doubly non-negative
  weighting matrices}.
\newblock \bibinfo{journal}{Advances in Neural Information Processing Systems}
  \bibinfo{volume}{34}.
\bibitem[{Leobacher and Steinicke(2022)}]{Leobacher_Steinicke:2022}
\bibinfo{author}{Leobacher, G.}, \bibinfo{author}{Steinicke, A.},
  \bibinfo{year}{2022}.
\newblock \bibinfo{title}{Exception sets of intrinsic and piecewise {Lipschitz}
  functions}.
\newblock \bibinfo{journal}{The Journal of Geometric Analysis}
  \bibinfo{volume}{32}.
\bibitem[{Li and Grossmann(2021)}]{Li_Grossmann:2021}
\bibinfo{author}{Li, C.}, \bibinfo{author}{Grossmann, I.E.},
  \bibinfo{year}{2021}.
\newblock \bibinfo{title}{A review of stochastic programming methods for
  optimization of process systems under uncertainty}.
\newblock \bibinfo{journal}{Frontiers in Chemical Engineering}
  \bibinfo{volume}{2}, \bibinfo{pages}{622241}.
\bibitem[{Li et~al.(2022)Li, Zhong and Brandeau}]{Li_et_al:2022}
\bibinfo{author}{Li, X.}, \bibinfo{author}{Zhong, H.},
  \bibinfo{author}{Brandeau, M.L.}, \bibinfo{year}{2022}.
\newblock \bibinfo{title}{Quantile {Markov} decision processes}.
\newblock \bibinfo{journal}{Operations Research} \bibinfo{volume}{70},
  \bibinfo{pages}{1428--1447}.
\bibitem[{Liu et~al.(2016)Liu, K{\"u}{\c{c}}{\"u}kyavuz and
  Luedtke}]{Liu_et_al:2016}
\bibinfo{author}{Liu, X.}, \bibinfo{author}{K{\"u}{\c{c}}{\"u}kyavuz, S.},
  \bibinfo{author}{Luedtke, J.}, \bibinfo{year}{2016}.
\newblock \bibinfo{title}{Decomposition algorithms for two-stage
  chance-constrained programs}.
\newblock \bibinfo{journal}{Mathematical Programming} \bibinfo{volume}{157},
  \bibinfo{pages}{219--243}.
\bibitem[{Luedtke(2014)}]{Luedtke:2014}
\bibinfo{author}{Luedtke, J.}, \bibinfo{year}{2014}.
\newblock \bibinfo{title}{A branch-and-cut decomposition algorithm for solving
  chance-constrained mathematical programs with finite support}.
\newblock \bibinfo{journal}{Mathematical Programming} \bibinfo{volume}{146},
  \bibinfo{pages}{219--244}.
\bibitem[{Mak et~al.(2015)Mak, Rong and Zhang}]{Mak_et_al:2015}
\bibinfo{author}{Mak, H.Y.}, \bibinfo{author}{Rong, Y.},
  \bibinfo{author}{Zhang, J.}, \bibinfo{year}{2015}.
\newblock \bibinfo{title}{Appointment scheduling with limited distributional
  information}.
\newblock \bibinfo{journal}{Management Science} \bibinfo{volume}{61},
  \bibinfo{pages}{316--334}.
\bibitem[{McKenzie et~al.(2023)McKenzie, Fung and Heaton}]{Mckenzie_et_al:2023}
\bibinfo{author}{McKenzie, D.}, \bibinfo{author}{Fung, S.W.},
  \bibinfo{author}{Heaton, H.}, \bibinfo{year}{2023}.
\newblock \bibinfo{title}{Differentiating through integer linear programs with
  quadratic regularization and {Davis-Yin} splitting}.
\newblock \bibinfo{journal}{arXiv preprint arXiv:2301.13395} .
\bibitem[{M{\'\i}nguez et~al.(2021)M{\'\i}nguez, van Ackooij and
  Garc{\'\i}a-Bertrand}]{Minguez_et_al:2021}
\bibinfo{author}{M{\'\i}nguez, R.}, \bibinfo{author}{van Ackooij, W.},
  \bibinfo{author}{Garc{\'\i}a-Bertrand, R.}, \bibinfo{year}{2021}.
\newblock \bibinfo{title}{Constraint generation for risk averse two-stage
  stochastic programs}.
\newblock \bibinfo{journal}{European Journal of Operational Research}
  \bibinfo{volume}{288}, \bibinfo{pages}{194--206}.
\bibitem[{Mu{\~n}oz et~al.(2022)Mu{\~n}oz, Pineda and
  Morales}]{Munoz_et_al:2022}
\bibinfo{author}{Mu{\~n}oz, M.A.}, \bibinfo{author}{Pineda, S.},
  \bibinfo{author}{Morales, J.M.}, \bibinfo{year}{2022}.
\newblock \bibinfo{title}{A bilevel framework for decision-making under
  uncertainty with contextual information}.
\newblock \bibinfo{journal}{Omega} \bibinfo{volume}{108},
  \bibinfo{pages}{102575}.
\bibitem[{Naumov and Bobylev(2012)}]{Naumov_Bobylev:2012}
\bibinfo{author}{Naumov, A.}, \bibinfo{author}{Bobylev, I.},
  \bibinfo{year}{2012}.
\newblock \bibinfo{title}{On the two-stage problem of linear stochastic
  programming with quantile criterion and discrete distribution of the random
  parameters}.
\newblock \bibinfo{journal}{Automation and Remote Control}
  \bibinfo{volume}{73}, \bibinfo{pages}{265--275}.
\bibitem[{Pavlikov et~al.(2017)Pavlikov, Veremyev and
  Pasiliao}]{Pavlikov_et_al:2017}
\bibinfo{author}{Pavlikov, K.}, \bibinfo{author}{Veremyev, A.},
  \bibinfo{author}{Pasiliao, E.L.}, \bibinfo{year}{2017}.
\newblock \bibinfo{title}{Optimization of value-at-risk: Computational aspects
  of {MIP} formulations}.
\newblock \bibinfo{journal}{Journal of the Operational Research Society} ,
  \bibinfo{pages}{1--15}.
\bibitem[{Pichler and Xu(2022)}]{Pichler_Xu:2022}
\bibinfo{author}{Pichler, A.}, \bibinfo{author}{Xu, H.}, \bibinfo{year}{2022}.
\newblock \bibinfo{title}{Quantitative stability analysis for minimax
  distributionally robust risk optimization}.
\newblock \bibinfo{journal}{Mathematical Programming} \bibinfo{volume}{191},
  \bibinfo{pages}{47--77}.
\bibitem[{Qi et~al.(2025)Qi, Grigas and Shen}]{Qi_et_al:2025}
\bibinfo{author}{Qi, M.}, \bibinfo{author}{Grigas, P.}, \bibinfo{author}{Shen,
  Z.J.}, \bibinfo{year}{2025}.
\newblock \bibinfo{title}{Integrated conditional estimation-optimization}.
\newblock \bibinfo{journal}{Operations Research} .
\bibitem[{Qiu et~al.(2014)Qiu, Ahmed, Dey and Wolsey}]{Qiu_et_al:2014}
\bibinfo{author}{Qiu, F.}, \bibinfo{author}{Ahmed, S.}, \bibinfo{author}{Dey,
  S.S.}, \bibinfo{author}{Wolsey, L.A.}, \bibinfo{year}{2014}.
\newblock \bibinfo{title}{Covering linear programming with violations}.
\newblock \bibinfo{journal}{INFORMS Journal on Computing} \bibinfo{volume}{26},
  \bibinfo{pages}{531--546}.
\bibitem[{Rachdi et~al.(2021)Rachdi, Laksaci and Al-Awadhi}]{Rachdi_et_al:2021}
\bibinfo{author}{Rachdi, M.}, \bibinfo{author}{Laksaci, A.},
  \bibinfo{author}{Al-Awadhi, F.A.}, \bibinfo{year}{2021}.
\newblock \bibinfo{title}{Parametric and nonparametric conditional quantile
  regression modeling for dependent spatial functional data}.
\newblock \bibinfo{journal}{Spatial Statistics} \bibinfo{volume}{43},
  \bibinfo{pages}{100498}.
\bibitem[{Rahimian and Pagnoncelli(2023)}]{Rahimian_Pagnoncelli:2023}
\bibinfo{author}{Rahimian, H.}, \bibinfo{author}{Pagnoncelli, B.},
  \bibinfo{year}{2023}.
\newblock \bibinfo{title}{Data-driven approximation of contextual
  chance-constrained stochastic programs}.
\newblock \bibinfo{journal}{SIAM Journal on Optimization} \bibinfo{volume}{33},
  \bibinfo{pages}{2248--2274}.
\bibitem[{Rahimian and Pagnoncelli(2024)}]{Rahimian_Pagnoncelli:2024}
\bibinfo{author}{Rahimian, H.}, \bibinfo{author}{Pagnoncelli, B.},
  \bibinfo{year}{2024}.
\newblock \bibinfo{title}{Contextual stochastic programs with expected-value
  constraints}.
\newblock \bibinfo{journal}{Optimization Online} .
\bibitem[{Rockafellar and Uryasev(2000)}]{Rockafellar_Uryasev:2000}
\bibinfo{author}{Rockafellar, R.T.}, \bibinfo{author}{Uryasev, S.},
  \bibinfo{year}{2000}.
\newblock \bibinfo{title}{Optimization of conditional value-at-risk}.
\newblock \bibinfo{journal}{Journal of Risk} \bibinfo{volume}{2},
  \bibinfo{pages}{21--42}.
\bibitem[{Sadana et~al.(2025)Sadana, Chenreddy, Delage, Forel, Frejinger and
  Vidal}]{Sadana_et_al:2025}
\bibinfo{author}{Sadana, U.}, \bibinfo{author}{Chenreddy, A.},
  \bibinfo{author}{Delage, E.}, \bibinfo{author}{Forel, A.},
  \bibinfo{author}{Frejinger, E.}, \bibinfo{author}{Vidal, T.},
  \bibinfo{year}{2025}.
\newblock \bibinfo{title}{A survey of contextual optimization methods for
  decision-making under uncertainty}.
\newblock \bibinfo{journal}{European Journal of Operational Research}
  \bibinfo{volume}{320}, \bibinfo{pages}{271--289}.
\bibitem[{Sadghiani and Motiian(2021)}]{Sadghiani_et_al:2021}
\bibinfo{author}{Sadghiani, N.S.}, \bibinfo{author}{Motiian, S.},
  \bibinfo{year}{2021}.
\newblock \bibinfo{title}{The contextual appointment scheduling problem}.
\newblock \bibinfo{journal}{arXiv preprint arXiv:2108.05531} .
\bibitem[{Sang et~al.(2021)Sang, Begen and Cao}]{Sang_et_al:2021}
\bibinfo{author}{Sang, P.}, \bibinfo{author}{Begen, M.A.},
  \bibinfo{author}{Cao, J.}, \bibinfo{year}{2021}.
\newblock \bibinfo{title}{Appointment scheduling with a quantile objective}.
\newblock \bibinfo{journal}{Computers \& Operations Research}
  \bibinfo{volume}{132}, \bibinfo{pages}{105295}.
\bibitem[{Shapiro et~al.(2014)Shapiro, Dentcheva and
  Ruszczy{\'n}ski}]{Shapiro_et_al:2014}
\bibinfo{author}{Shapiro, A.}, \bibinfo{author}{Dentcheva, D.},
  \bibinfo{author}{Ruszczy{\'n}ski, A.}, \bibinfo{year}{2014}.
\newblock \bibinfo{title}{Lectures on Stochastic Programming: Modeling and
  Theory}.
\newblock \bibinfo{edition}{Second} ed., \bibinfo{publisher}{SIAM}.
\bibitem[{Soleimani et~al.(2014)Soleimani, Seyyed-Esfahani and
  Kannan}]{Soleimani_et_al:2014}
\bibinfo{author}{Soleimani, H.}, \bibinfo{author}{Seyyed-Esfahani, M.},
  \bibinfo{author}{Kannan, G.}, \bibinfo{year}{2014}.
\newblock \bibinfo{title}{Incorporating risk measures in closed-loop supply
  chain network design}.
\newblock \bibinfo{journal}{International Journal of Production Research}
  \bibinfo{volume}{52}, \bibinfo{pages}{1843--1867}.
\bibitem[{Srivastava et~al.(2021)Srivastava, Wang, Hanasusanto and
  Ho}]{Srivastava_et_al:2021}
\bibinfo{author}{Srivastava, P.R.}, \bibinfo{author}{Wang, Y.},
  \bibinfo{author}{Hanasusanto, G.A.}, \bibinfo{author}{Ho, C.P.},
  \bibinfo{year}{2021}.
\newblock \bibinfo{title}{On data-driven prescriptive analytics with side
  information: A regularized nadaraya-watson approach}.
\newblock \bibinfo{journal}{arXiv preprint arXiv:2110.04855} .
\bibitem[{Sun et~al.(2023)Sun, Liu and Li}]{Sun_et_al:2023}
\bibinfo{author}{Sun, C.}, \bibinfo{author}{Liu, L.}, \bibinfo{author}{Li, X.},
  \bibinfo{year}{2023}.
\newblock \bibinfo{title}{Predict-then-calibrate: A new perspective of robust
  contextual {LP}}.
\newblock \bibinfo{journal}{Advances in Neural Information Processing Systems}
  \bibinfo{volume}{36}, \bibinfo{pages}{17713--17741}.
\bibitem[{Sun et~al.(2022)Sun, Shi, Wang, Tuan, Poor and Tao}]{Sun_et_al:2022}
\bibinfo{author}{Sun, H.}, \bibinfo{author}{Shi, Y.}, \bibinfo{author}{Wang,
  J.}, \bibinfo{author}{Tuan, H.D.}, \bibinfo{author}{Poor, H.V.},
  \bibinfo{author}{Tao, D.}, \bibinfo{year}{2022}.
\newblock \bibinfo{title}{Alternating differentiation for optimization layers}.
\newblock \bibinfo{journal}{arXiv preprint arXiv:2210.01802} .
\bibitem[{Sun and Xu(2016)}]{Sun_Xu:2016}
\bibinfo{author}{Sun, H.}, \bibinfo{author}{Xu, H.}, \bibinfo{year}{2016}.
\newblock \bibinfo{title}{Convergence analysis for distributionally robust
  optimization and equilibrium problems}.
\newblock \bibinfo{journal}{Mathematics of Operations Research}
  \bibinfo{volume}{41}, \bibinfo{pages}{377--401}.
\bibitem[{Tao et~al.(2025)Tao, Delage and Xu}]{Tao_et_al:2025}
\bibinfo{author}{Tao, Y.}, \bibinfo{author}{Delage, E.}, \bibinfo{author}{Xu,
  H.}, \bibinfo{year}{2025}.
\newblock \bibinfo{title}{Risk-averse decision making with contextual
  information: Model, sample average approximation, and kernelization}.
\newblock \bibinfo{journal}{arXiv preprint arXiv:2502.16607} .
\bibitem[{Thiele et~al.(2009)Thiele, Terry and Epelman}]{Thiele_et_al:2009}
\bibinfo{author}{Thiele, A.}, \bibinfo{author}{Terry, T.},
  \bibinfo{author}{Epelman, M.}, \bibinfo{year}{2009}.
\newblock \bibinfo{title}{Robust linear optimization with recourse}.
\newblock \bibinfo{type}{Technical Report}. Available in Optimization Online.
\bibitem[{Tsang and Shehadeh(2025)}]{Tsang_Shehadeh:2025}
\bibinfo{author}{Tsang, M.Y.}, \bibinfo{author}{Shehadeh, K.S.},
  \bibinfo{year}{2025}.
\newblock \bibinfo{title}{On the tradeoff between distributional belief and
  ambiguity: Conservatism, finite-sample guarantees, and asymptotic
  properties}.
\newblock \bibinfo{journal}{INFORMS Journal on Optimization} .
\bibitem[{Tsang et~al.(2023)Tsang, Shehadeh and Curtis}]{Tsang_et_al:2023}
\bibinfo{author}{Tsang, M.Y.}, \bibinfo{author}{Shehadeh, K.S.},
  \bibinfo{author}{Curtis, F.E.}, \bibinfo{year}{2023}.
\newblock \bibinfo{title}{An inexact column-and-constraint generation method to
  solve two-stage robust optimization problems}.
\newblock \bibinfo{journal}{Operations Research Letters} \bibinfo{volume}{51},
  \bibinfo{pages}{92--98}.
\bibitem[{Walk(2010)}]{Walk:2010}
\bibinfo{author}{Walk, H.}, \bibinfo{year}{2010}.
\newblock \bibinfo{title}{Strong laws of large numbers and nonparametric
  estimation}, in: \bibinfo{booktitle}{Recent Developments in Applied
  Probability and Statistics}. \bibinfo{publisher}{Springer}, pp.
  \bibinfo{pages}{183--214}.
\bibitem[{Wang et~al.(2023)Wang, Becker, Van~Parys and
  Stellato}]{Wang_et_al:2023}
\bibinfo{author}{Wang, I.}, \bibinfo{author}{Becker, C.},
  \bibinfo{author}{Van~Parys, B.}, \bibinfo{author}{Stellato, B.},
  \bibinfo{year}{2023}.
\newblock \bibinfo{title}{Learning decision-focused uncertainty sets in robust
  optimization}.
\newblock \bibinfo{journal}{arXiv preprint arXiv:2305.19225} .
\bibitem[{Wang et~al.(2018)Wang, Zhou, Song and Sherwood}]{Wang_et_al:2018}
\bibinfo{author}{Wang, L.}, \bibinfo{author}{Zhou, Y.}, \bibinfo{author}{Song,
  R.}, \bibinfo{author}{Sherwood, B.}, \bibinfo{year}{2018}.
\newblock \bibinfo{title}{Quantile-optimal treatment regimes}.
\newblock \bibinfo{journal}{Journal of the American Statistical Association}
  \bibinfo{volume}{113}, \bibinfo{pages}{1243--1254}.
\bibitem[{Wang et~al.(2020)Wang, Chen and Liu}]{Wang_et_al:2020}
\bibinfo{author}{Wang, S.}, \bibinfo{author}{Chen, Z.}, \bibinfo{author}{Liu,
  T.}, \bibinfo{year}{2020}.
\newblock \bibinfo{title}{Distributionally robust hub location}.
\newblock \bibinfo{journal}{Transportation Science} \bibinfo{volume}{54},
  \bibinfo{pages}{1189--1210}.
\bibitem[{Zeng et~al.(2014)Zeng, An and Kuznia}]{Zeng_et_al:2014}
\bibinfo{author}{Zeng, B.}, \bibinfo{author}{An, Y.}, \bibinfo{author}{Kuznia,
  L.}, \bibinfo{year}{2014}.
\newblock \bibinfo{title}{Chance constrained mixed integer program: Bilinear
  and linear formulations, and {Benders} decomposition}.
\newblock \bibinfo{journal}{arXiv preprint arXiv:1403.7875} .
\bibitem[{Zhenevskaya and Naumov(2018)}]{Zhenevskaya_Naumov:2018}
\bibinfo{author}{Zhenevskaya, I.D.}, \bibinfo{author}{Naumov, A.V.},
  \bibinfo{year}{2018}.
\newblock \bibinfo{title}{The decomposition method for two-stage stochastic
  linear programming problems with quantile criterion}.
\newblock \bibinfo{journal}{Automation and Remote Control}
  \bibinfo{volume}{79}, \bibinfo{pages}{229--240}.
\bibitem[{Zhou et~al.(2019)Zhou, Shahidehpour, Wei, Li, Sun and
  Chen}]{Zhou_et_al:2019}
\bibinfo{author}{Zhou, Y.}, \bibinfo{author}{Shahidehpour, M.},
  \bibinfo{author}{Wei, Z.}, \bibinfo{author}{Li, Z.}, \bibinfo{author}{Sun,
  G.}, \bibinfo{author}{Chen, S.}, \bibinfo{year}{2019}.
\newblock \bibinfo{title}{Distributionally robust unit commitment in
  coordinated electricity and district heating networks}.
\newblock \bibinfo{journal}{IEEE Transactions on Power Systems}
  \bibinfo{volume}{35}, \bibinfo{pages}{2155--2166}.
\bibitem[{Ziegel(2016)}]{Ziegel:2016}
\bibinfo{author}{Ziegel, J.F.}, \bibinfo{year}{2016}.
\newblock \bibinfo{title}{Coherence and elicitability}.
\newblock \bibinfo{journal}{Mathematical Finance} \bibinfo{volume}{26},
  \bibinfo{pages}{901--918}.

\end{thebibliography}

\end{document}